\magnification=\magstep1
\input amstex
\documentstyle{amsppt}
\pagewidth{30pc}\pageheight{48pc}
\input xy
\xyoption{all}

\topmatter
\title Galois extensions of structured ring spectra
\endtitle
\author John Rognes \endauthor
\date August 31st 2005 \enddate
\address Department of Mathematics, University of Oslo, Norway \endaddress
\email rognes\@math.uio.no \endemail
\abstract
We introduce the notion of a Galois extension of commutative
$S$-algebras ($E_\infty$ ring spectra), often localized with respect to
a fixed homology theory.  There are numerous examples, including some
involving Eilenberg--Mac\,Lane spectra of commutative rings, real and
complex topological $K$-theory, Lubin--Tate spectra and cochain
$S$-algebras.  We establish the main theorem of Galois theory in this
generality.  Its proof involves the notions of separable and {\'e}tale
extensions of commutative $S$-algebras, and the Goerss--Hopkins--Miller
theory for $E_\infty$ mapping spaces.  We show that the global sphere
spectrum~$S$ is separably closed, using Minkowski's discriminant
theorem, and we estimate the separable closure of its localization with
respect to each of the Morava $K$-theories.  We also define
Hopf--Galois extensions of commutative $S$-algebras, and study the
complex cobordism spectrum $MU$ as a common integral model for all of
the local Lubin--Tate Galois extensions.
\endabstract
\toc
\widestnumber\head{10}
\widestnumber\subhead{10.1}
\head 1. Introduction \endhead
\head 2. Galois extensions in algebra \endhead
\subhead 2.1. Galois extensions of fields \endsubhead
\subhead 2.2. Regular covering spaces \endsubhead
\subhead 2.3. Galois extensions of commutative rings \endsubhead
\head 3. Closed categories of structured module spectra \endhead
\subhead 3.1. Structured spectra \endsubhead
\subhead 3.2. Localized categories \endsubhead
\subhead 3.3. Dualizable spectra \endsubhead
\subhead 3.4. Stably dualizable groups \endsubhead
\subhead 3.5. The dualizing spectrum \endsubhead
\subhead 3.6. The norm map \endsubhead
\head 4. Galois extensions in topology \endhead
\subhead 4.1. Galois extensions of $E$-local commutative $S$-algebras
	\endsubhead
\subhead 4.2. The Eilenberg--Mac\,Lane embedding \endsubhead
\subhead 4.3. Faithful extensions \endsubhead
\head 5. Examples of Galois extensions \endhead
\subhead 5.1. Trivial extensions \endsubhead
\subhead 5.2. Eilenberg--Mac\,Lane spectra \endsubhead
\subhead 5.3. Real and complex topological $K$-theory \endsubhead
\subhead 5.4. The Morava change-of-rings theorem \endsubhead
\subhead 5.5. The $K(1)$-local case \endsubhead 
\subhead 5.6. Cochain $S$-algebras \endsubhead
\head 6. Dualizability and alternate characterizations \endhead
\subhead 6.1. Extended equivalences \endsubhead
\subhead 6.2. Dualizability \endsubhead
\subhead 6.3. Alternate characterizations \endsubhead
\subhead 6.4. The trace map and self-duality \endsubhead
\subhead 6.5. Smash invertible modules \endsubhead
\head 7. Galois theory I \endhead
\subhead 7.1. Base change for Galois extensions \endsubhead
\subhead 7.2. Fixed $S$-algebras \endsubhead
\head 8. Pro-Galois extensions and the Amitsur complex \endhead
\subhead 8.1. Pro-Galois extensions \endsubhead
\subhead 8.2. The Amitsur complex \endsubhead
\head 9. Separable and {\'e}tale extensions \endhead
\subhead 9.1. Separable extensions \endsubhead
\subhead 9.2. Symmetrically {\'e}tale extensions \endsubhead
\subhead 9.3. Smashing maps \endsubhead
\subhead 9.4. {\'E}tale extensions \endsubhead
\subhead 9.5. Henselian maps \endsubhead
\subhead 9.6. $I$-adic towers \endsubhead
\head 10. Mapping spaces of commutative $S$-algebras \endhead
\subhead 10.1. Obstruction theory \endsubhead
\subhead 10.2. Idempotents and connected $S$-algebras \endsubhead
\subhead 10.3. Separable closure \endsubhead
\head 11. Galois theory II \endhead
\subhead 11.1. Recovering the Galois group \endsubhead
\subhead 11.2. The brave new Galois correspondence \endsubhead
\head 12. Hopf--Galois extensions in topology \endhead
\subhead 12.1.  Hopf--Galois extensions of commutative $S$-algebras
	\endsubhead
\subhead 12.2. Complex cobordism \endsubhead
\head {\ } References \endhead
\endtoc
\endtopmatter


\define\ADer{\operatorname{\Cal A Der}}
\define\Alg{\operatorname{Alg}}
\define\Aut{\operatorname{Aut}}
\define\CDer{\operatorname{\Cal C Der}}
\define\C{\Bbb C}
\define\Ext{\operatorname{Ext}}
\define\E{\operatorname{E}}
\define\F{\Bbb F}
\define\Gal{\operatorname{Gal}}
\define\G{\Bbb G}
\define\Hom{\operatorname{Hom}}
\define\Map{\operatorname{Map}}
\define\Pic{\operatorname{Pic}}
\define\Q{\Bbb Q}
\define\SS{\Bbb S}
\define\Spec{\operatorname{Spec}}
\define\TAQ{T\!AQ}
\define\THH{T\!H\!H}
\define\Tor{\operatorname{Tor}}
\define\Tot{\operatorname{Tot}}
\define\W{\Bbb W}
\define\Z{\Bbb Z}
\define\colim{\operatornamewithlimits{colim}}

\define\holim{\operatornamewithlimits{holim}}
\define\llangle{\langle\!\langle}
\define\rrangle{\rangle\!\rangle}
\define\tmf{tm\!f}
\define\wprod{\mathop{{\prod}'}}

\hyphenation{topo-logy}
\hyphenation{homo-logy}

\document

\head 1. Introduction \endhead

The present paper is motivated by (1) the ``brave new rings'' paradigm
coined by Friedhelm Waldhausen, that structured ring spectra are an
unavoidable generalization of discrete rings, with arithmetic
properties captured by their algebraic $K$-theory, (2) the presumption
that algebraic $K$-theory will satisfy an extended form of the
{\'e}tale- and Galois descent foreseen by Dan Quillen, and (3) the
algebro-geometric perspective promulgated by Jack Morava, on how the
height-stratified moduli space of formal group laws influences stable
homotopy theory, by way of complex cobordism theory.

We here develop the arithmetic notion of a Galois extension of structured
ring spectra, viewed geometrically as an algebraic form of a regular
covering space, by always working intrinsically in a category of spectra,
rather than at the na{\"\i}ve level of coefficient groups.  The result
is a framework that well accommodates much recent work in stable homotopy
theory.  We hope that this study will eventually lead to a conceptual
understanding of objects like the algebraic $K$-theory of the sphere
spectrum, which by Waldhausen's stable parametrized $h$-cobordism theorem
bears on such seemingly unrelated geometric objects as the diffeomorphism
groups of manifolds, in much the same way that we now understand the
algebraic $K$-theory spectrum of the ring of integers.

\medskip

Let $E$ be any spectrum and $G$ a finite group.  We say that a map $A \to
B$ of $E$-local commutative $S$-algebras is an {\it $E$-local $G$-Galois
extension\/} if $G$ acts on $B$ through commutative $A$-algebra maps in
such a way that the two canonical maps
$$
i \: A \to B^{hG}
$$
and
$$
h \: B \wedge_A B \to \prod_G B
$$
induce isomorphisms in $E_*$-homology (Definition~4.1.3).  When $E =
S$ this means that the maps $i$ and $h$ are weak equivalences, and we
may talk of a {\it global\/} $G$-Galois extension.  In more detail,
the map $i$ is the standard inclusion into the homotopy fixed points
for the $G$-action on $B$ and $h$ is given in symbols by $h(b_1 \wedge
b_2) = \{g \mapsto b_1 \cdot g(b_2)\}$.  To make the definition homotopy
invariant we also assume that $A$ is a cofibrant commutative $S$-algebra
and that $B$ is a cofibrant commutative $A$-algebra.

There are many interesting examples of such ``brave new'' Galois
extensions.

\example{Examples 1.1}

(a)
The Eilenberg--Mac\,Lane functor $R \mapsto HR$ takes each $G$-Galois
extension $R \to T$ of commutative rings to a global $G$-Galois extension
$HR \to HT$ of commutative $S$-algebras (Proposition~4.2.1).

(b)
The complexification map $KO \to KU$ from real to complex topological
$K$-theory is a global $\Z/2$-Galois extension (Proposition~5.3.1).

(c)
For each rational prime $p$ and natural number $n$ the profinite
extended Morava stabilizer group $\G_n = \SS_n \rtimes \Gal$ acts on
the even periodic Lubin--Tate spectrum $E_n$, with $\pi_0(E_n) =
\W(\F_{p^n})[[u_1, \dots, u_{n-1}]]$, so that $L_{K(n)} S \to E_n$ is a
$K(n)$-local pro-$\G_n$-Galois extension (see Notation~3.2.2 and
Theorem~5.4.4(d)).

(d)
For most regular covering spaces $Y \to X$ the map of cochain
$H\F_p$-algebras $F(X_+, H\F_p) \to F(Y_+, H\F_p)$ is a Galois extension
(Proposition~5.6.3(a)).
\endexample

A map $A \to B$ of commutative $S$-algebras will be said to be {\it
faithful\/} if for each $A$-module $N$ with $N \wedge_A B \simeq *$ we
have $N \simeq *$ (Definition~4.3.1).  The map $A \to B$ is {\it
separable\/} if the multiplication map $\mu \: B \wedge_A B \to B$
admits a bimodule section up to homotopy (Definition~9.1.1).  A
commutative $S$-algebra $B$ is {\it connected,} in the sense of
algebraic geometry, if its space of idempotents $\Cal E(B)$ is weakly
equivalent to the two-point space $\{0,1\}$ (Definition~10.2.1).  There
are analogous definitions in each $E$-local context.

In commutative ring theory each Galois extension is faithful, but
it remains an open problem to decide whether each Galois extension
of commutative $S$-algebras is faithful (Question~4.3.6).  Rather
conveniently, a commutative $S$-algebra $B$ is connected if and only if
the ring $\pi_0(B)$ is connected (Proposition~10.2.2).

Here is our version of the Main Theorem of Galois theory for commutative
$S$-algebras.  The first two parts~(a) and~(b) of the theorem are
obtained by specializing Theorem~7.2.3 and Proposition~9.1.4 to the
case of a finite, discrete Galois group $G$.  The recovery in~(c) of
the Galois group is Theorem~11.1.1.  The converse part~(d) is the less
general part of Theorem~11.2.2.

\proclaim{Theorem 1.2}
Let $A \to B$ be a faithful $E$-local $G$-Galois extension.

(a)
For each subgroup $K \subset G$ the map $C = B^{hK} \to B$ is a faithful
$E$-local $K$-Galois extension, with $A \to C$ separable.

(b)
For each normal subgroup $K \subset G$ the map $A \to C = B^{hK}$ is a
faithful $E$-local $G/K$-Galois extension.

If furthermore $B$ is connected, then:

(c)
The Galois group $G$ is weakly equivalent to the mapping space $\Cal
C_A(B, B)$ of commutative $A$-algebra self-maps of $B$.

(d)
For each factorization $A \to C \to B$ of the $G$-Galois extension,
with $A \to C$ separable and $C \to B$ faithful, there is a subgroup $K
\subset G$ such that $C \simeq B^{hK}$ as an $A$-algebra over~$B$.
\endproclaim

In other words, for a faithful $E$-local $G$-Galois extension $A
\to B$ with $B$ connected there is a bijective contravariant Galois
correspondence $K \leftrightarrow C = B^{hK}$ between the subgroups of
$G$ and the weak equivalence classes of separable $A$-algebras mapping
faithfully to $B$.  The inverse correspondence takes $C$ to $K = \pi_0
\Cal C_C(B, B)$.

The main theorem fully describes the intermediate extensions in a
$G$-Galois extension $A \to B$, but what about the further extensions
of $B$?  We say that a connected $E$-local commutative $S$-algebra $A$ is
{\it separably closed\/} if there are no connected $E$-local $G$-Galois
extensions $A \to B$ for non-trivial groups $G$ (Definition~10.3.1).
The following fundamental example is a consequence of Minkowski's
discriminant theorem in number theory, and is proved as Theorem~10.3.3.

\proclaim{Theorem 1.3}
The global sphere spectrum $S$ is separably closed.
\endproclaim

The absence of localization is crucial for this result.  At the other
extreme the Morava $K(n)$-local category is maximally localized, for
each $p$ and~$n$.  Here the Lubin--Tate spectrum $E_n$ admits a
$K(n)$-local pro-$n\hat\Z$-Galois extension $E_n \to E_n^{nr}$, with
$$
\pi_0(E_n^{nr}) = \W(\bar\F_p)[[u_1, \dots, u_{n-1}]]
$$
given by adjoining all roots of unity of order prime to~$p$ (\S5.4.6).
We expect that each further $G$-Galois extension $E_n^{nr} \to B$ of such
a Landweber exact even periodic spectrum must again be Landweber exact
and even periodic, and such that $\pi_0(E_n^{nr}) \to \pi_0(B)$ will be
a $G$-Galois extension of commutative rings.  But $\W(\bar\F_p)[[u_1,
\dots, u_{n-1}]]$ is separably closed as a commutative ring, so such
a $\pi_0(B)$ cannot be connected, and $B$ the cannot be connected for
non-trivial groups $G$.  Therefore we expect:

\proclaim{Conjecture 1.4}
The extension $E_n^{nr}$ of the Lubin--Tate spectrum $E_n$ is
$K(n)$-locally separably closed.  In particular, the Galois group
$\G_n^{nr} = \SS_n \rtimes \hat\Z$ of $L_{K(n)} S \to E_n^{nr}$ is the
$K(n)$-local absolute Galois group of the $K(n)$-local sphere spectrum
$L_{K(n)}S$.
\endproclaim

Partial results supporting this conjecture have been obtained by Andy
Baker and Birgit Richter \cite{BR:r}, for global Galois extensions that
are furthermore assumed to be faithful and abelian.

The substantial supply of pro-Galois extensions in the $K(n)$-local
category, like $L_{K(n)}S \to E_n$, is not available in the
$E(n)$-local category (see \S5.5.4).  This draws extra attention to the
non-smashing Bousfield localizations, and thus to the distinction
between the whole category of modules over $L_{K(n)}S$ and its full
subcategory of $K(n)$-local modules.  A study of the sphere spectrum as
an algebro-geometric scheme- or stack-like object, that only involves
smashing localizations or only treats the whole module categories over
the various Bousfield localizations, does thus not capture these very
interesting examples of regular geometric covering spaces.

There are structured ring spectrum replacements for K{\"a}hler
differentials, called topological Hochschild homology (= $\THH$, see
Section~9.2) and topological Andr{\'e}--Quillen homology (= $\TAQ$, see
Section~9.4), in the context of associative and commutative
$S$-algebras, respectively.  These need not be $K(n)$-local when
applied to $K(n)$-local $S$-algebras (see Example~9.2.3).  Therefore
the notions of (formally) {\'e}tale extensions of associative or
commutative $S$-algebras will again give a richer theory when
considered within the $K(n)$-local subcategory, rather than in the
whole module category over $L_{K(n)}S$.  Thus also a study of the
algebraic geometry of the sphere spectrum with respect to the {\'e}tale
topology will become more substantial by taking these Bousfield local
subcategories into account.  This phenomenon differs from that which is
familiar in discrete algebraic geometry, since there all localizations
are, indeed, smashing.

It therefore appears to be better to think of the algebraic geometry of
the sphere spectrum as the ``$S$-algebraic stack'' of all Bousfield
$E$-local subcategories $\Cal M_{S,E}$ of spectra, for varying spectra
$E$, rather than the ``$S$-algebraic scheme'' of the Bousfield
$E$-local $S$-algebras $L_E S$ themselves.  The former stack maps to
the stack of module categories of the latter scheme, but it is the
former that carries the most interesting closed symmetric monoidal
structures.  See Definition~3.2.1 for the notations used here, and
Section~3.2, Chapter~9 and Section~12.2 for more on these
$S$-algebro-geometric ideas.

The (mono-)chromatic localizations $L_{K(n)}S$ of the sphere are of course
even more drastic than the $p$-localizations $S_{(p)}$, so that many of
the principal examples studied in this paper are of an even more local
nature than e.g.~local number fields.  But the arithmetic properties
of a global number field can usefully be studied by ad{\`e}lic means,
in terms of the system of local number fields that can be obtained from
it by the various completions that are available.  We are therefore also
interested in finding global models for the system of naturally occurring
$K(n)$-local Galois extensions of $L_{K(n)}S$, for varying $p$ and $n$.

The obvious candidate, given Quillen's discovery of the relation of
formal group law theory to complex cobordism, is the unit map $S \to
MU$ to the complex cobordism spectrum.  The following statement is
proved in Corollary~9.6.6, Proposition~12.2.1 and the discussion
surrounding diagram~(12.2.6).  In the second part, $S[BU]$ is the
commutative Hopf $S$-algebra $\Sigma^\infty BU_+$.  In summary, $MU$ is
very close to such a global model, up to formal thickenings by
Henselian maps.  This makes the author inclined to think of $S \to MU$
as a kind of (large) ramified global Galois extension, with $S[BU]$
playing the part of the functional dual of its imaginary Galois group.
To make good sense of this, we introduce the notion of a Hopf--Galois
extension of commutative $S$-algebras in Section~12.1.

\proclaim{Theorem 1.5}
For each prime $p$ and integer $n\ge1$ the $K(n)$-local pro-$\G_n$-Galois
extension $L_{K(n)}S \to E_n$ factors as the composite of the following
maps of commutative $S$-algebras
$$
L_{K(n)}S \to \hat L^{MU}_{K(n)}MU @>q>> \widehat{E(n)} \to E_n \,.
$$
Here the first map admits the global model $S \to MU$, by Bousfield
$K(n)$-localization in the category of $S$-modules and $K(n)$-nilpotent
completion in the category of $MU$-modules, respectively.  The second
map $q$ is a formal thickening, or more precisely, symmetrically (and
possibly commutatively) Henselian.  The third map is a finite Galois
extension (and can be avoided by passing to the even periodic version
$MUP$ of $MU$ and adjoining some roots of unity).

Furthermore, the global model $S \to MU$ is an $S[BU]$-Hopf--Galois
extension of commutative $S$-algebras, with coaction $\beta \: MU \to
MU \wedge S[BU]$ given by the Thom diagonal.  For each element $g \in
\G_n$ its Galois action on $E_n$ can be directly recovered from this
$S[BU]$-coaction, up to the adjunction of some roots of unity.
\endproclaim

Here are some more detailed references into the body of the paper.

Chapter~2 contains a review of the basic Galois theory for fields and for
commutative rings, together with some algebraic facts that we will need
for our examples.  We also make a comparison with the theory of regular
covering spaces, for the benefit of the topologically minded reader.

As hinted at above, we sometimes consider more general Galois groups
$G$ than finite (and profinite) groups.  For the initial theory, all
that is required is that the unreduced suspension spectrum $S[G] = L_E
\Sigma^\infty G_+$ admits a good Spanier--Whitehead dual in the
$E$-local stable homotopy category, i.e., that $G$ is {\it stably
dualizable} (Definition~3.4.1).  We review the basic properties of
stably dualizable groups and their actions on spectra in Chapter~3,
referring to the author's companion paper \cite{Rog:s} for most
proofs.  This chapter also contains a discussion of the various
categories of $E$-local $S$-modules and (commutative) $S$-algebras in
which we work.

The precise Definition~4.1.3 of a Galois extension of commutative
$S$-algebras is given in Chapter~4, followed by a discussion showing
that the Eilenberg--Mac\,Lane embedding from commutative rings preserves
and detects Galois extensions (Proposition~4.2.1).  We also consider the
elementary properties of faithful modules over structured ring spectra,
flatness being implicit in our homotopy invariant work.  We shall often
make use of how various algebro-geometric properties of $S$-algebras are
preserved by base change, or are detected by suitable forms of faithful
base change.

Chapter~5 is devoted to the many examples of Galois extensions
mentioned above, including all the intermediate $K(n)$-local Galois
extensions between $L_{K(n)} S$ and the maximal unramified extension
$E_n^{nr}$ of $E_n$.  We also go through the $K(1)$-local case of the
Lubin--Tate extensions in much detail, making explicit the close
analogy with the classification of abelian extensions of the $p$-adic
and rational fields $\Q_p$ and~$\Q$.  Finally we extend the example of
cochain algebras of regular covering spaces to cochain algebras of
principal $G$-bundles $P \to X$, for stably dualizable groups $G$.

Chapter~6 develops the formal consequences of the Galois conditions on
$A \to B$, including the basic fact that $B$ is a dualizable $A$-module
(Proposition~6.2.1), two useful alternate characterizations of (faithful)
Galois extensions (Propositions~6.3.1 and~6.3.2), and two further
characterizations of faithfulness (Proposition~6.3.3 and Lemma~6.5.4).
These let us prove in Chapter~7 that faithful Galois extensions are
preserved by arbitrary base change (Lemma~7.1.1) and are detected
by faithful and dualizable base change (Lemma~7.1.4(b)).  From these
results, in turn, the ``forward'' part of the Galois correspondence
(Theorem~7.2.3) follows rather formally, saying that for a faithful
$G$-Galois extension $A \to B$ the homotopy fixed point spectra $C =
B^{hK}$ give rise to $K$-Galois extensions $C \to B$ for subgroups $K
\subset G$, and to $G/K$-Galois extensions $A \to C$ when $K$ is normal.

When this much of the Galois correspondence has been established, we can
make sense of the notion of a pro-Galois extension, which we do somewhat
informally in Section~8.1.

The ``converse'' part of the Galois correspondence (Theorem~11.2.2) relies
on the possibility of recovering the Galois group $G$ in a $G$-Galois
extension $A \to B$ from the space $\Cal C_A(B, B)$ of commutative
$A$-algebra self-maps $B \to B$, or more generally, to recover the
subgroup $K$ from the mapping space $\Cal C_C(B, B)$, when $C = B^{hK}$
is a fixed $S$-algebra of $B$ (Proposition~11.2.1).  This is achieved in
Chapter~11, but relies on three preceding developments.

First of all, we use the commutative form of the Hopkins--Miller
theory, as developed by Paul Goerss and Mike Hopkins \cite{GH04}, to
study such mapping spaces.  We use an extension of their work, from
dealing with spaces of $E_\infty$ ring spectrum maps, or commutative
$S$-algebra maps, to spaces of commutative $A$-algebra maps.  This is
discussed in Section~10.1, where we also touch on the consequences for
this theory of working $E$-locally.  The main computational tool is the
Goerss--Hopkins spectral sequence~(10.1.4), whose $E_2$-term involves
suitable Andr{\'e}--Quillen cohomology groups, which fortunately vanish
in all relevant cases for the Galois extensions that we consider.

Second, the recovery of the Galois group $G$ from $\Cal C_A(B, B)$
only has a chance, judging from the discrete algebraic case, when $B$ is
connected in the geometric sense that it has no non-trivial idempotents.
For a commutative $S$-algebra $B$ there is a space $\Cal E(B)$ of
idempotents, which in turn is a commutative $B$-algebra mapping space of
the sort that can be studied by the Goerss--Hopkins spectral sequence.
So in Section~10.2 we treat connectivity in this geometric sense for
commutative $S$-algebras, reaching a convenient algebraic criterion in
Proposition~10.2.2.  This also lets us define separably closed commutative
$S$-algebras in Section~10.3.

Thirdly, not all commutative $A$-algebras $C$ mapping faithfully to $B$
occur in the Galois correspondence as fixed $S$-algebras $C = B^{hK}$.
As in the discrete algebraic case, the characteristic property is that
$C$ must be separable over $A$, and in Section~9.1 we develop the basic
theory of separable extensions of $S$-algebras.  As further
generalizations of separable maps we have the {\'e}tale maps, which we
discuss in three related contexts in Sections~9.2 through~9.4, leading
to the notions of symmetrically (=thh-){\'e}tale, smashing and
(commutatively) {\'e}tale maps of $S$-algebras, respectively.

Topological Hochschild homology $\THH$ controls the K{\"a}hler
differentials in the associative setting, while topological
Andr{\'e}--Quillen homology $\TAQ$ takes on the same role in the
purely commutative setting.  Our discussion here relies heavily on
the work of Maria Basterra \cite{Bas99} and Andrej Lazarev \cite{La01}.
There is much conceptual overlap between the triviality of the topological
Andr{\'e}--Quillen homology spectrum $\TAQ(B/A)$ for (formally)
{\'e}tale maps $A \to B$, and the vanishing of the Goerss--Hopkins
Andr{\'e}--Quillen cohomology groups $D^s_{B_*T}(B^A_*(B), \Omega^t B)$
for finite Galois extensions $A \to B$, but the direct connection is
not as well understood as might be desired.

The remainder of the text is concerned with the interpretation of $S
\to MU$ as a Hopf--Galois extension that provides an integral model, up
to Henselian maps, for all of the Lubin--Tate extensions $L_{K(n)}S \to
E_n$.  Thus we consider square-zero extensions, singular extensions and
Henselian maps as various forms of infinitesimal and formal thickenings
in Section~9.5.  We then obtain a good supply of relevant examples
in Section~9.6, using work of Baker and Lazarev on $I$-adic towers.
We have already cited Corollary~9.6.6 as relevant for part of Theorem~1.5.

The idea of Hopf--Galois extensions is to replace the action by the
Galois group~$G$ on a commutative $A$-algebra $B$ by a coaction by the
functional dual $DG_+ = F(G_+, S)$ of the Galois group, which is a
commutative Hopf $S$-algebra.  In the algebraic situation such
coactions have been useful, e.g.~to classify inseparable Galois
extensions of fields \cite{Cha71}.  In the absence of an actual Galois
group the condition that $i \: A \to B^{hG}$ is a weak equivalence must
be rewritten, by using a cosimplicial resolution for the coaction (the
Hopf cobar complex) in place of the homotopy fixed points.  This
rewriting can naturally go through a second cosimplicial resolution
associated to $A \to B$, which we know as the Amitsur complex.  We
discuss the Amitsur complex in Section~8.2, so as to have the
accompanying notion of (nilpotent) completion of $A$ along $B$
available in Chapters~9 and~10, and give the definitions of the Hopf
cobar complex and of Hopf--Galois extensions in Section~12.1.

To conclude the paper, in Section~12.2 we go through some of the details
of how the inseparable extension $S \to MU$ is an $S[BU]$-Hopf--Galois
extension, and how the Hopkins--Miller theory and the Lubin--Tate
deformation theory work together to show that the global $S[BU]$-coaction
on $MU$ captures the Morava stabilizer group action on $E_n$, at all
primes $p$ and chromatic heights $n$.

\subhead Acknowledgments \endsubhead

The study of idempotents in Chapter~10 got going during an Oberwolfach
hike with Neil Strickland, and the right way to use separability in
Chapter~11 was found in a discussion with Birgit Richter.  I am very
grateful to them, as well as to Andy Baker, Daniel G. Davis, Halvard
Fausk, Jack Morava and the referee, for a number of helpful or
encouraging comments.

Much of this work was done in the year 2000 and announced at various
conferences.  I apologize for the long delay in publication, which for
most of the time was due to the unresolved Question~4.3.6, on the
faithfulness of Galois extensions.

\head 2. Galois extensions in algebra \endhead

\subhead 2.1. Galois extensions of fields \endsubhead

We first recall the basics about Galois extensions of fields.  Let $G$
be a finite group acting effectively (only the unit element acts as the
identity) from the left by automorphisms on a field $E$, and let $F =
E^G$ be the fixed subfield.  Let
$$
j \: E\langle G\rangle \to \Hom_F(E, E)
$$
be the canonical associative ring homomorphism taking $e_1 g$ to the
homomorphism $e_2 \mapsto e_1 \cdot g(e_2)$, from the {\it twisted
group ring\/} of $G$ over $E$ to the $F$-module endomorphisms of $E$.
Then $j$ is an isomorphism, for by Dedekind's lemma $j$ is injective,
and $\dim_F(E)$ equals the order of $G$, so $j$ is also surjective by
a dimension count.  See \cite{Dr95, App.} for elementary proofs.  Let
$$
h \: E \otimes_F E \to \prod_G E
$$
be the canonical commutative ring homomorphism taking $e_1 \otimes e_2$
to the sequence $\{g \mapsto e_1 \cdot g(e_2)\}$, from the tensor product
of two copies of $E$ over $F$ to the product of $G$ copies of $E$.
Then also $h$ is an isomorphism, for it is the $E$-module dual of $j$,
by way of the identifications $\Hom_E(E \otimes_F E, E) \cong \Hom_F(E,
E)$ and $\Hom_E(\prod_G E, E) \cong E\langle G\rangle$ (using that $G$
is finite).

\subhead 2.2. Regular covering spaces \endsubhead

There is a parallel geometric theory of regular (= normal) covering spaces
\cite{Sp66, 2.6.7}, \cite{Ha02, 1.39}.  Let $G$ be a finite discrete group
acting from the right on a compact Hausdorff space $Y$.  Let $X = Y/G$
be the orbit space, and let $\pi \: Y \to X$ be the orbit projection.
There is a canonical map
$$
\xi \: Y \times G \to Y \times_X Y
$$
to the fiber product of $\pi$ with itself, taking $(y, g)$ to $(y, y
\cdot g)$.  This map is always surjective, by the definition of $X$ as
an orbit space, and it is injective if and only if $G$ acts freely on
$Y$.  So $\xi$ is a homeomorphism if and only if $Y \to X$ is a regular
covering space, with $G$ as its group of deck transformations, acting
freely and transitively on each fiber.  In general, the possible
failure of $\xi$ to be injective measures the extent to which $G$ does
not act freely on $Y$, which in turn can be interpreted as a measure of
to what extent $Y$ is {\it ramified\/} as a cover of $X$.  The theory
of Riemann surfaces provides numerous examples of the latter
phenomenon.

Dually, let $R = C(X)$ and $T = C(Y)$ be the rings of continuous (real or
complex) functions on $X$ and $Y$, respectively.  As usual the points of
$X$ can be recovered as the maximal ideals in $R$, and similarly for $Y$.
The group $G$ acts from the left on~$T$, by the formula $g(t) = g * t
\: y \mapsto t(y \cdot g)$, and the natural map $R \to T$ dual to~$\pi$
identifies $R$ with the invariant ring $T^G$, by the isomorphism $C(Y)^G
\cong C(Y/G)$.  The map $\xi$ above is dual to the canonical homomorphism
$$
h \: T \otimes_R T \to \prod_G T
$$
taking $t_1 \otimes t_2$ to the function $g \mapsto t_1 \cdot g(t_2)$,
considered as an element in the product $\prod_G T$.  Then $\xi$ is a
homeomorphism if and only if $h$ is an isomorphism, by the categorical
anti-equivalence between compact Hausdorff spaces and their function
rings.  The surjectivity of $\xi$ ensures that $h$ is always injective,
and in general the possible failure of $h$ to be surjective measures
the extent of ramification in $Y \to X$.

For a moment, let us also consider the more general case of a principal
$G$-bundle $\pi \: P \to X$ for a compact Hausdorff topological group $G$.
The map $\xi \: P \times G \to P \times_X P$ is a homeomorphism, now
with respect to the given topology on $G$.  Let $R = C(X)$, $T = C(P)$
and $H = C(G)$.  Then $H$ is a commutative Hopf algebra with coproduct
$\psi \: H \to H \otimes H$, if the map $H \to C(G \times G)$ dual to
the group multiplication $G \times G \to G$ factors through the canonical
map $H \otimes H \to C(G \times G)$.  Likewise $H$ coacts on $T$ from the
right by $\beta \: T \to T \otimes H$, if the map $T \to C(P \times G)$
induced by the group action $P \times G \to P$ factors through $T \otimes
H \to C(P \times G)$.  These factorizations can always be achieved by
using suitably completed tensor products, but we wish to refer to the
algebraic tensor products here.  Then the freeness of the group action
on $P$ is expressed by saying that the composite map
$$
h \: T \otimes_R T @>1\otimes\beta>> T \otimes_R T \otimes H
@>\mu\otimes1>> T \otimes H
$$
is an isomorphism.  We shall return to this dualized context in Chapter~12
on Hopf--Galois extensions.

\subhead 2.3. Galois extensions of commutative rings \endsubhead

Generalizing the two examples above, for finite Galois groups, Auslander
and Goldman \cite{AG60, App.} gave a definition of Galois extensions
of commutative rings as part of their study of separable algebras
over such rings.  Chase, Harrison and Rosenberg \cite{CHR65, \S1} found
several other equivalent definitions, and developed the Galois theory for
commutative rings to also encompass the fundamental Galois correspondence.
We now recall their basic results.

Let $R \to T$ be a homomorphism of commutative rings, making $T$ a
commutative $R$-algebra, and let $G$ be a finite group acting on $T$
from the left through $R$-algebra homomorphisms.  Let
$$
i \: R \to T^G
$$
be the inclusion into the fixed ring, let
$$
h \: T \otimes_R T \to \prod_G T
$$
be the commutative ring homomorphism that takes $t_1 \otimes t_2$ to
the sequence $\{g \mapsto t_1 \cdot g(t_2)\}$, and let
$$
j \: T \langle G\rangle \to \Hom_R(T, T)
$$
be the associative ring homomorphism that takes $t_1 g$ to the
$R$-module homomorphism $t_2 \mapsto t_1 \cdot g(t_2)$.  We give
$\prod_G T$ the pointwise product $(t_g)_g \cdot (t'_g)_g = (t_g
t'_g)_g$ and $T \langle G \rangle$ the twisted product $t_1 g_1 \cdot
t_2 g_2 = t_1 g_1(t_2) g_1 g_2$, using the left $G$-action $(g_1, t_2)
\mapsto g_1(t_2)$ on~$T$.

\definition{Definition 2.3.1}
Let $G$ act on $T$ over $R$, as above.  We say that $R \to T$ is a {\it
$G$-Galois extension of commutative rings\/} if both $i \: R \to T^G$
and $h \: T \otimes_R T \to \prod_G T$ are isomorphisms.
\enddefinition

Here we are following Greither \cite{Gre92, 0.1.5}.  Auslander and
Goldman \cite{AG60, p.~396} instead took the condition below on $i$, $j$
and $T$ to be the defining property, but Chase, Harrison and Rosenberg
\cite{CHR65, 1.3} proved that the two definitions are equivalent.

\proclaim{Proposition 2.3.2}
Let $G$ act on $T$ over $R$, as above.  Then $R \to T$ is a $G$-Galois
extension if and only if both $i \: R \to T^G$ and $j \: T\langle
G\rangle \to \Hom_R(T, T)$ are isomorphisms and $T$ is a finitely
generated projective $R$-module.
\endproclaim

The condition that $i$ is an isomorphism means that we can speak of $R$
as the {\it fixed ring\/} of $T$.  The homomorphism $h$ measures to what
extent the extension $R \to T$ is {\it ramified,} and Galois extensions
are required to be unramified.  The injectivity of $j$ is a form of
Dedekind's lemma, and ensures that the action by $G$ is effective.

\example{Example 2.3.3}
If $K \to L$ is a $G$-Galois extension of number fields, then the
corresponding extension $R = \Cal O_K \to \Cal O_L = T$ of rings of
integers is a $G$-Galois extension of commutative rings if and only if
$K \to L$ is unramified as an extension of number fields \cite{AB59}.
More generally, if $\Sigma$ is a set of prime ideals in $\Cal O_K$, and
$\Sigma'$ the set of primes in $\Cal O_L$ above those in $\Sigma$, then
the extension $\Cal O_{K,\Sigma} \to \Cal O_{L,\Sigma'}$ of rings of
$\Sigma$-integers is $G$-Galois if and only if $\Sigma$ contains all
the primes that ramify in $L/K$ \cite{Gre92, 0.4.1}.  Here $\Cal
O_{K,\Sigma}$ is defined as the ring of elements $x \in K$ that have
non-negative valuation $v_{\goth p}(x) \ge 0$ for all prime ideals
$\goth p \notin \Sigma$.  Thus $\Cal O_K \to \Cal O_L$ becomes a
$G$-Galois extension precisely upon localizing away from (= inverting)
the ramified primes.

To see this, note that if $T = R\{t_1, \dots, t_n\}$ is a free $R$-module
of rank~$n$, then $T \otimes_R T$ is a free $T$-module on the generators
$1 \otimes t_1, \dots, 1 \otimes t_n$, and $h$ is represented as
a $T$-module homomorphism by the square matrix $A = (g(t_i))_{g,i}$
of rank~$n$, with $g \in G$ and $i = 1, \dots, n$.  The discriminant of
$T/R$ is $d = \det(A)^2$, by definition, and the prime ideals in $\Cal
O_K$ that ramify in $L/K$ are precisely the prime ideals dividing the
discriminant.  So $h$ is an isomorphism if and only if $\det(A)$ and
$d$ are units in $R$, or equivalently, if there are no ramified primes.
A local version of the same argument works when $T$ is not free over $R$.
\endexample

Here are some further basic properties of Galois extensions of commutative
rings, which will be relevant to our discussion.

\proclaim{Proposition 2.3.4}
Let $R \to T$ be a $G$-Galois extension.  Then:

(a) $T$ is faithfully flat as an $R$-module, i.e., the functor $(-)
\otimes_R T$ preserves and detects (=reflects) exact sequences.

(b) The trace map $tr \: T \to R$ (taking $t \in T$ to
$\sum_{g\in G} g(t) \in T^G = R$) is a split surjective
$R$-module homomorphism.

(c) $T$ is invertible as an $R[G]$-module, i.e., a finitely generated
projective $R[G]$-module of constant rank~$1$.
\endproclaim

For proofs, see e.g.~\cite{Gre92, 0.1.9}, \cite{Gre92, 0.1.10} and
\cite{Gre92, 0.6.1}.  Beware that part~(b) does not extend well to the
topological setting, as Example~6.4.4 demonstrates.

\head 3. Closed categories of structured module spectra \endhead

\subhead 3.1. Structured spectra \endsubhead

We now adapt these ideas to the context of ``brave new rings,'' i.e.,
of commutative $S$-algebras.  These can in Chapters 2--9 and~12 be
interpreted as the commutative monoids in either one of the popular
symmetric monoidal categories of structured spectra, such as the
$S$-modules of Elmendorf, Kriz, Mandell and May \cite{EKMM97}, the
symmetric spectra in simplicial sets of Hovey, Shipley and Smith
\cite{HSS00}, symmetric spectra in topological spaces or orthogonal
spectra of Mandell, May, Schwede and Shipley \cite{MMSS01} or the
simplicial functors of Segal and Lydakis \cite{Ly98}, according to the
reader's needs or preferences.

However, in Chapters~10 and~11 we make use of the Goerss--Hopkins
obstruction theory for $E_\infty$ mapping spaces \cite{GH04}, which
presumes that one works in a category of spectra that satisfies
Axioms~1.1 and~1.4 in {\it op. cit.}  In particular, this theory is
needed for the proof of parts~(c) and~(d) of our Theorem~1.2, and for
Theorem~1.3.  It is known that $S$-modules, symmetric spectra formed in
topological spaces and orthogonal spectra all satisfy the required
axioms, by \cite{GH04, 1.5}.  To be concrete, and to have a convenient
source for the more technical references, we shall work with the
$S$-modules of Peter May {\it et al.}

Let $S$ be the sphere spectrum, and let $\Cal M_S$ be the category of
$S$-modules.  Among other things, it is a topological category with all
limits and colimits and all topological tensors and cotensors.  A map
$f \: X \to Y$ of $S$-modules is called a {\it weak equivalence\/} if
the induced homomorphism $\pi_*(f) \: \pi_*(X) \to \pi_*(Y)$ of stable
homotopy groups is an isomorphism.  The category $\Cal D_S$ obtained
from $\Cal M_S$ by inverting the weak equivalences is called the {\it
stable homotopy category,} and is equivalent to the homotopy category
of spectra constructed by Boardman \cite{Vo70}.

The smash product $X \wedge Y$ and function object $F(X, Y)$ make $\Cal
M_S$ a closed symmetric monoidal category, with $S$ as the unit object.
For each topological space $T$ the topological tensor $X \wedge T_+$
equals the smash product $X \wedge S[T]$, and the topological cotensor
$Y^T = F(T_+, Y)$ equals the function spectrum $F(S[T], Y)$, where $S[T]
= \Sigma^\infty T_+$ denotes the unreduced suspension $S$-module on $T$.

An (associative) $S$-algebra $A$ is a monoid in $\Cal M_S$, i.e.,
an $S$-module $A$ equipped with a unit map $\eta \: S \to A$ and
a unital and associative multiplication $\mu \: A \wedge A \to A$.
A commutative $S$-algebra $A$ is a commutative monoid in $\Cal M_S$,
i.e., one such that the multiplication $\mu$ is also commutative.
We write $\Cal A_S$ and $\Cal C_S$ for the categories of $S$-algebras
and commutative $S$-algebras, respectively.  More generally, for a
commutative $S$-algebra $A$ we write $\Cal M_A$, $\Cal A_A$ and $\Cal
C_A$ for the categories of $A$-modules, associative $A$-algebras and
commutative $A$-algebras, respectively \cite{EKMM97, VII.1}.

\subhead 3.2. Localized categories \endsubhead

Our first examples of Galois extensions of structured ring spectra will be
maps $A \to B$ of commutative $S$-algebras, with a finite group $G$ acting
on $B$ through $A$-algebra maps, such that there are weak equivalences
$i \: A \simeq B^{hG}$ and $h \: B \wedge_A B \simeq \prod_G B$.
The formal definition appears in Section~4.1 below.  However, there
are interesting examples that only appear as Galois extensions to the
eyes of weaker invariants than the stable homotopy groups $\pi_*(-)$.
More precisely, for a fixed homology theory $E_*(-)$ we shall allow
ourselves to work in the $E$-local stable homotopy category, where have
arranged that each map $f \: X \to Y$ such that $E_*(f) \: E_*(X) \to
E_*(Y)$ is an isomorphism, is in fact a weak equivalence.  In particular,
we will encounter situations where we only have that $E_*(i)$ and
$E_*(h)$ are isomorphisms, in which case we shall interpret $A \to B$
as an $E$-local $G$-Galois extension.

Note the close analogy between the $E$-local theory and the case
(Example~2.3.3) of rings of integers localized away from some set of
primes.  Doug Ravenel's influential treatise on the chromatic
filtration of stable homotopy theory \cite{Ra84, \S5}, brings emphasis
to the tower of cases when $E = E(n)$, the $n$-th Johnson--Wilson
spectrum.  To us, the most interesting case is when $E = K(n)$ is the
$n$-th Morava $K$-theory spectrum.  The $K(n)$-local stable homotopy
category is studied in detail in \cite{HSt99, \S\S7--8}, and captures
the $n$-th layer, or stratum, in the chromatic filtration.

\definition{Definition 3.2.1}
Let $E$ be a fixed $S$-module, with associated homology theory $X
\mapsto E_*(X) = \pi_*(E \wedge X)$.  By definition, an $S$-module $Z$
is said to be {\it $E$-acyclic\/} if $E \wedge Z \simeq *$, and an
$S$-module $Y$ is said to be {\it $E$-local\/} if $F(Z, Y) \simeq *$
for each $E$-acyclic $S$-module $Z$.  Let $\Cal M_{S,E} \subset \Cal
M_S$ be the full subcategory of {\it $E$-local $S$-modules.}  A map $f
\: X \to Y$ of $E$-local $S$-modules is a weak equivalence if and only
if it is an {\it $E_*$-equivalence,} i.e., if $E_*(f)$ is an
isomorphism.

There is a Bousfield localization functor $L_E \: \Cal M_S \to \Cal
M_{S,E} \subset \Cal M_S$ \cite{Bo79}, \cite{EKMM97, VIII.1.6},
and an accompanying natural $E_*$-equivalence $X \to L_E X$ for each
$S$-module $X$.  We may assume that this $E_*$-equivalence is the identity
when $X$ is already $E$-local, so that the localization functor $L_E$
is idempotent.  The homotopy category $\Cal D_{S,E}$ of $\Cal M_{S,E}$
is the {\it $E$-local stable homotopy category.}

More generally, for a commutative $S$-algebra $A$ we let $\Cal M_{A,E}
\subset \Cal M_A$ be the full subcategory of $E$-local $A$-modules, with
homotopy category $\Cal D_{A,E}$.  To be precise, there is an $A$-module
$\F_AE$ of the homotopy type of $A \wedge E$, and a localization functor
$L^A_{\F_AE}\: \Cal M_A \to \Cal M_{A,E}$, with respect to $\F_AE$ in
the category of $A$-modules, which amounts to $E$-localization at the
level of the underlying $S$-modules \cite{EKMM97, VIII.1.7}.  We shall
allow ourselves to simply denote this functor by $L_E$.
\enddefinition

\definition{Notation 3.2.2}
We write
$$
L_n X = L_{E(n)} X
$$
for the Bousfield localization of $X$ with respect to the Johnson--Wilson
spectrum $E(n)$ \cite{JW73}, with $\pi_* E(n) = \Z_{(p)}[v_1, \dots,
v_{n-1}, v_n^{\pm1}]$, for each non-negative integer~$n$, and
$$
L_{K(n)} X
$$
for the Bousfield localization of $X$ with respect to the
Morava $K$-theory spectrum $K(n)$ \cite{JW75}, with $\pi_* K(n) =
\F_p[v_n^{\pm1}]$, for each natural number~$n$.
\enddefinition

We will reserve the symbol $\hat L$ for the Bousfield nilpotent
completion recalled in Definition~8.2.2, and shall therefore not use
this notation for the functor $L_{K(n)}$, unlike e.g.~\cite{HSt99}.

The smash product $X \wedge Y$ of two $E$-local $S$-modules will in
general not be $E$-local, although this is the case when $L_E$ is a
so-called {\it smashing localization,} i.e., one that commutes with
direct limits \cite{Ra84, 1.28}.  The Johnson--Wilson spectra $E =
E(n)$ provide interesting examples of smashing localizations $L_n =
L_{E(n)}$ \cite{Ra92, 7.5.6}, while localization $L_{K(n)}$ with respect
to the Morava $K$-theories $E = K(n)$ is not smashing \cite{HSt99, 8.1}.
Likewise, the unit $S$ for the smash product is rarely $E$-local.  So in
order to work with $S$-algebras and related constructions internally
within $\Cal M_{S,E}$, we first perform each construction as usual in
$\Cal M_S$, and then apply the Bousfield localization functor~$L_E$.

\definition{Definition 3.2.3}
We implicitly give $\Cal M_{S,E}$ all colimits, topological tensors,
smash products and a unit object by applying Bousfield localization
to the constructions in~$\Cal M_S$.  So $\colim_{i \in I} X_i$ means
$L_E(\colim_{i \in I} X_i)$, $X \wedge Y$ means $L_E(X \wedge Y)$, $S$
means $L_E S$ and $S[T]$ means $L_E \Sigma^\infty T_+$.  All limits,
topological cotensors and function objects formed from $E$-local
$S$-modules are already $E$-local, so no Bousfield localization is
required in these cases.  With these conventions, $\Cal M_{S,E}$
is a topological closed symmetric monoidal category with all limits
and colimits.  The same considerations apply for $\Cal M_{A,E}$.
\enddefinition

There is a natural map $L_E X \wedge L_E Y \to L_E(X \wedge Y)$, making
$L_E$ a lax monoidal functor, so that $L_E S$ is always a commutative
$S$-algebra.  When $E$ is smashing, the category $\Cal M_{S,E}$ of
$E$-local $S$-modules is equivalent (at the level of homotopy categories)
to the category $\Cal M_{L_E S}$ of $L_E S$-modules, so the study of
$E$-local $S$-modules is a special case of the study of modules over a
general commutative $S$-algebra $A = L_E S$.  However, when $E$ is not
smashing, as is the case for $E = K(n)$, the two homotopy categories
are not equivalent, and we shall need to consider the more general notion.

When $E = S$, every $S$-module is $E$-local and $\Cal M_{S,E} = \Cal
M_S$, etc., so the $E$-local context specializes to the ``global'',
unlocalized situation.  For brevity, we shall often simply refer to the
$E$-local $S$-modules as $S$-modules, or even as spectra, but except
where we explicitly assume that $E = S$, the discussion is intended to
encompass also the general $E$-local case.

\remark{Remark 3.2.4}
By analogy with algebraic geometry, we may heuristically wish to view
$A$-modules $M$ as suitable sheaves $M^\sim$ over some geometric
``structure space'' $\Spec A$.  This structure space would come with a
Zariski topology, with open subspaces $U_{A,E} \subset \Spec A$
corresponding to the various localization functors $L_E$ on the
category of $A$-modules, in such a way that the restriction of the
sheaf $M^\sim$ over $\Spec A$ to the subspace $U_{A,E}$ would be the
sheaf $(L_E M)^\sim$ corresponding to the $E$-local $A$-module $L_E
M$.  For smashing $E$ this would precisely amount to an $L_E A$-module,
so that $U_{A,E}$ could be identified with the structure space $\Spec
L_E A$.

However, for non-smashing $E$ the condition of being an $E$-local
$A$-module is strictly stronger than being an $L_E A$-module.
Therefore, the geometric structure on $\Spec A$ is not simply that of
an ``$S$-algebra'ed space'' carrying the (commutative) $S$-algebra $L_E
A$ over $U_{A,E}$, by analogy with the ringed spaces of algebraic
geometry.  If we wish to allow non-smashing localizations $E$ to
correspond to Zariski opens, then the geometric structure must also
capture the additional restriction it is for an $L_E A$-module to be an
$E$-local $A$-module.  This exhibits a difference compared to the
situation in commutative algebra, where localization at an ideal
commutes with direct limits, and behaves as a smashing localization,
while completions behave more like non-smashing localizations.  It does
not seem to be so common to do commutative algebra in such implicitly
completed situations, however.

A continuation of this analogy would be to consider other
Grothendieck-type topologies on $\Spec A$, with coverings built from
$E$-local Galois extensions $L_E A \to B$ (Definition~4.1.3) or more
general {\'e}tale extensions (Definition~9.4.1), subject to a combined
faithfulness condition (Definition~4.3.1).  In the unlocalized cases,
such a (big) {\'e}tale site on the opposite category of $\Cal C_S$, and
associated small {\'e}tale sites on the opposite category of each $\Cal
C_A$, have been developed by Bertrand To{\"e}n and Gabriele Vezzosi
\cite{TV05, \S5.2}.  However, the rich source of $K(n)$-local Galois
extensions of $L_{K(n)} S$ discussed in Section~5.4 provides, by
Lemma~9.4.4, an equally rich supply of $K(n)$-local {\'e}tale maps from
$L_{K(n)} S$.  It appears, by extension from the case $n=1$ discussed
in Section~5.5, that these are not globally {\'e}tale maps, in which
case the {\'e}tale topology proposed in \cite{TV05} will be too coarse
to encompass these examples.  The author therefore thinks that a finer
{\'e}tale site, taking non-smashing localizations like $L_{K(n)}$ into
account, would lead to a stronger and more interesting theory.
\endremark

\subhead 3.3. Dualizable spectra \endsubhead

In each closed symmetric monoidal category there is a canonical natural
map
$$
\nu \: F(X, Y) \wedge Z \to F(X, Y \wedge Z) \,.
$$
It is right adjoint to a map $\epsilon\wedge1 \: X \wedge F(X, Y) \wedge
Z \to Y \wedge Z$, where the adjunction counit $\epsilon \: X \wedge F(X,
Y) \to Y$ is left adjoint to the identity map on $F(X, Y)$.

Dold and Puppe \cite{DP80} say that an object $X$ is strongly dualizable
if the canonical map $\nu \: F(X, Y) \wedge Z \to F(X, Y \wedge Z)$ is an
isomorphism for all $Y$ and $Z$.  Lewis, May and Steinberger \cite{LMS86,
III.1.1} say that a spectrum $X$ is finite if it is strongly dualizable
in the stable homotopy category, i.e., if the map $\nu$ is a weak
equivalence.  We shall instead follow Hovey and Strickland \cite{HSt99,
1.5(d)} and briefly call such spectra dualizable.  By \cite{LMS86,
III.1.3(ii)} it suffices to verify this condition in the special case when
$Y = S$ and $Z = X$, so we take this simpler criterion as our definition.

\definition{Definition 3.3.1}
Let $DX = F(X, S)$ be the {\it functional dual\/} of $X$.  We say that $X$
is {\it dualizable\/} if the canonical map $\nu \: DX \wedge X \to F(X,
X)$ is a weak equivalence.  More generally, for an (implicitly $E$-local)
module~$M$ over a commutative $S$-algebra~$A$, let $D_A M = F_A(M, A)$
be the functional dual, and say that $M$ is a dualizable $A$-module if
the canonical map $\nu \: D_A M \wedge_A M \to F_A(M, M)$ is a weak
equivalence.
\enddefinition

\proclaim{Lemma 3.3.2}
(a) If $X$ or $Z$ is dualizable, then the canonical map $\nu \: F(X, Y)
\wedge Z \to F(X, Y \wedge Z)$ is a weak equivalence.

(b) If $X$ is dualizable, then $DX$ is also dualizable and the
canonical map $\rho \: X \to DDX$ is a weak equivalence.

(c) The dualizable spectra generate a thick subcategory, i.e., they are
closed under passage to weakly equivalent objects, retracts, mapping
cones and (de-)suspensions.
\endproclaim

Here $\rho \: X \to DDX = F(F(X, S), S)$ is right adjoint to $F(X, S)
\wedge X \to S$, which is obtained by twisting the adjunction counit
$\epsilon \: X \wedge F(X, S) \to S$.  For proofs, see \cite{LMS86,
III.1.2 and III.1.3}.  We sometimes also use $\nu$ to label the conjugate
map $Y \wedge F(X, Z) \to F(X, Y \wedge Z)$.  The corresponding results
hold for $E$-local $A$-modules, by the same formal proofs.

One justification for the term ``finite'' is the following converse to
Lemma~3.3.2(c), in the unlocalized setting $E = S$.

\proclaim{Proposition 3.3.3}
Let $A$ be commutative $S$-algebra.  A global $A$-module $M$ is
dualizable in $\Cal M_A = \Cal M_{A,S}$ if and only if it is weakly
equivalent to a retract of a finite cell $A$-module.  When $A$ is
connective, this is in turn equivalent to being a retract of a
finite CW $A$-module spectrum.
\endproclaim

The proof \cite{EKMM97, III.7.9} uses in an essential way that stable
homotopy $X \mapsto \pi_*(X) = [A, X]^A_*$ commutes with coproducts,
which amounts to $A$ being small in the homotopy category $\Cal D_A$
of $A$-modules.  This fails in some $E$-local contexts.  For example,
the $K(n)$-local sphere spectrum $L_{K(n)} S$ is not small in the
$K(n)$-local category \cite{HSt99, 8.1}, and consequently $\pi_*(X)$
is not a homology theory on this category.  So in general there will
be more dualizable $E$-local $A$-modules than the {\it semi-finite\/}
ones, i.e., the retracts of the finite cell $L_E A$-modules.  In this
paper we shall prefer to focus on the notion of dualizability, rather
than on being semi-finite, principally because of Proposition~6.2.1
and (counter-)Example~6.2.2 below.

\subhead 3.4. Stably dualizable groups \endsubhead

For our basic theory of $G$-Galois extensions of commutative
$S$-algebras the group action by $G$ appears through the module
action by its suspension spectrum $S[G] = L_E \Sigma^\infty G_+$, and
the finiteness condition on $G$ only enters through the property that
$S[G]$ is a dualizable spectrum.  We then say that $G$ is an $E$-locally
stably dualizable group.  Only when we turn to properties related to
separability will it be relevant that $G$ is discrete, and then usually
finite.  So we shall develop the basic theory in the greater generality
of stably dualizable topological groups $G$.

\definition{Definition 3.4.1}
A topological group $G$ is {\it $E$-locally stably dualizable\/} if
its suspension spectrum $S[G] = L_E \Sigma^\infty G_+$ is dualizable in
$\Cal M_{S,E}$.  Writing $DG_+ = F(G_+, L_E S)$ for its functional dual,
the condition is that the canonical map
$$
\nu \: DG_+ \wedge S[G] \to F(S[G], S[G])
$$
is a weak equivalence in the $E$-local category.
\enddefinition

\example{Examples 3.4.2}
(a)
Each compact Lie group $G$ admits the structure of a finite CW complex,
so $S[G]$ is a finite cell spectrum and $G$ is stably dualizable, for
each~$E$.

(b)
The Eilenberg--Mac\,Lane spaces $G = K(\Z/p, q)$ are loop spaces and
thus admit models as topological groups.  They have infinite mod~$p$
homology for each $q\ge1$, so $S[G]$ is never dualizable in $\Cal M_S$
by Proposition~3.3.3.  However, the Morava $K$-homology $K(n)_* K(\Z/p,
q)$ is finitely generated over $K(n)_*$ by a calculation of Ravenel and
Wilson \cite{RaW80, 9.2}, so $G = K(\Z/p, q)$ is in fact $K(n)$-locally
stably dualizable by \cite{HSt99, 8.6}.  We are curious to see if these
and similar topological Galois groups play any significant role in
the $K(n)$-local Galois theory.
\endexample

\subhead 3.5. The dualizing spectrum \endsubhead

The weak equivalence $S[G] = \bigvee_G S \to \prod_G S = DG_+$
for a finite group $G$ generalizes to an $E$-local self-duality of
the suspension spectrum $S[G]$, when $G$ is an $E$-locally stably
dualizable group.  The self-duality holds up to a twist by a so-called
dualizing spectrum $S^{adG}$.  When $G$ is a compact Lie group this is
the suspension spectrum on the one-point compactification of the adjoint
representation $adG$ of~$G$, thus the notation, and so $S^{adG} = S$
for $G$ finite.  John Klein \cite{Kl01, \S1} introduced dualizing spectra
$S^{adG}$ for arbitrary topological groups, and Tilman Bauer \cite{Bau04,
4.1} established the twisted self-duality of $S[G]$ in the $p$-complete
category, when $G$ is a $p$-compact group in the sense of Bill Dwyer and
Clarence Wilkerson \cite{DW94}.  In \cite{Rog:s} we have extended these
results to all $E$-locally stably dualizable groups, as we now review.

\definition{Definition 3.5.1}
Let $G$ be an $E$-locally stably dualizable group.  The group
multiplication provides the suspension spectrum $S[G] = L_E
\Sigma^\infty G_+$ with mutually commuting left and right $G$-actions.
We define the {\it dualizing spectrum} $S^{adG}$ to be the
$G$-homotopy fixed point spectrum
$$
S^{adG} = S[G]^{hG} = F(EG_+, S[G])^G
$$
of $S[G]$, formed with respect to the right $G$-action \cite{Rog:s,
2.5.1}.  Here $EG = B(*, G, G)$ is the standard free, contractible
right $G$-space.  The remaining left action on $S[G]$ induces a left
$G$-action on $S^{adG}$.
\enddefinition

When $G$ is finite, there is a natural weak equivalence
$$
S^{adG} = S[G]^{hG} \simeq DG_+^{hG} \simeq S \,.
$$
Here the last equivalence involves the collapsing homotopy equivalence $c
\: EG \to *$, which is a $G$-equivariant map, but not a $G$-equivariant
homotopy equivalence.  For general stably dualizable groups $G$, the
dualizing spectrum is indeed dualizable and smash invertible \cite{Rog:s,
3.2.3 and 3.3.4}, so smashing with $S^{adG}$ induces an equivalence of
derived categories.

The left $G$-action on $S[G]$ functorially dualizes to a right
$G$-action on $DG_+$, with associated module action map $\alpha \: DG_+
\wedge S[G] \to DG_+$.  The diagonal map on~$G$ induces a coproduct
$\psi \: S[G] \to S[G] \wedge S[G]$, using \cite{EKMM97, II.1.2}.
These combine to a {\it shear map}
$$
sh \: DG_+ \wedge S[G] @>1\wedge\psi>> DG_+ \wedge S[G] \wedge S[G]
@>\alpha\wedge1>> DG_+ \wedge S[G] \,,
$$
which is equivariant with respect to each of three mutually commuting
$G$-actions \cite{Rog:s, 3.1.2} and is a weak equivalence \cite{Rog:s, 3.1.3}.
Taking homotopy fixed points with respect to the right action of $G$ on
$S[G]$ in the source and the diagonal right action on $DG_+$ and $S[G]$
in the target induces a natural {\it Poincar{\'e} duality equivalence\/}
\cite{Rog:s, 3.1.4}
$$
DG_+ \wedge S^{adG} @>\simeq>> S[G] \,.
\tag 3.5.2
$$
This identification uses the stable dualizability of $G$, and expresses
the twisted self-duality of $S[G]$.  The weak equivalence is equivariant
with respect to both a left and a right $G$-action.  The left $G$-action
is by the inverse of the right action on $DG_+$, the standard left action
on $S^{adG}$ and the standard left action on $S[G]$.  The right $G$-action
is by the inverse of the left action on $DG_+$, the trivial action on
$S^{adG}$ and the standard right action on $S[G]$.

\subhead 3.6. The norm map \endsubhead

Let $X$ be any $E$-local $S$-module with left $G$-action, and equip it
with the trivial right $G$-action.  The smash product $X \wedge S[G]$
then has a diagonal left $G$-action, and a right $G$-action that only
affects $S[G]$.  Consider forming homotopy orbits $(-)_{hG}$ with respect
to the left action and forming homotopy fixed points $(-)^{hG}$ with
respect to the right action, in either order.  There is then a canonical
colimit/limit exchange map
$$
\kappa \: ((X \wedge S[G])^{hG})_{hG} \to ((X \wedge S[G])_{hG})^{hG}
\,.
$$
The source of $\kappa$ receives a weak equivalence from $(X \wedge
S^{adG})_{hG}$ (this uses the stable dualizability of $G$; see the proof
of Lemma~6.4.2), and the target of $\kappa$ maps by a weak equivalence to
$X^{hG}$ (this is easy).  The composite of these three maps is the {\it
(homotopy) norm map\/} \cite{Rog:s, 5.2.2}
$$
N \: (X \wedge S^{adG})_{hG} \to X^{hG} \,.
\tag 3.6.1
$$
If $X = W \wedge G_+ = W \wedge S[G]$ for some spectrum $W$ with left
$G$-action, with $G$ acting in the standard way on $S[G]$, then the
norm map for $X$ is a weak equivalence \cite{Rog:s, 5.2.5}.  That
reference only discusses the case when $G$ acts trivially on~$W$, but
in general there is an equivariant shearing equivalence $\zeta \: w
\wedge g \mapsto g(w) \wedge g$ from $W \wedge S[G]$ with $G$ acting
only on $S[G]$ to $W \wedge S[G]$ with the diagonal $G$-action.

We can define the {\it $G$-Tate construction} $X^{tG}$ to be the cofiber
of the norm map
$$
(X \wedge S^{adG})_{hG} @>N>> X^{hG} @>>> X^{tG} \,.
$$
Then $X^{tG} \simeq *$ if and only if $N$ is a weak equivalence, which in
turn holds if and only if the exchange map $\kappa$ is a weak equivalence.
From this point of view $X^{tG}$ is the obstruction to the commutation
of the $G$-homotopy orbit and the $G$-homotopy fixed point constructions,
when applied to $X \wedge S[G]$.

\head 4. Galois extensions in topology \endhead

\subhead 4.1. Galois extensions of $E$-local commutative $S$-algebras
\endsubhead

Fix an $S$-module $E$, and consider the categories $\Cal M_{S,E}$ and
$\Cal C_{S,E}$ of $E$-local $S$-modules and $E$-local commutative
$S$-algebras, respectively.  These are full subcategories of the
topological (closed) model categories $\Cal M_S$ and $\Cal C_S$,
respectively, as explained in \cite{EKMM97, VII.4}.

The reader may, if preferred, alternatively work with the
``convenient'' $S$-model structures of Jeff Smith and Brooke Shipley
\cite{Sh04}, but this will not be necessary.  There is another
$E$-local model structure on $\Cal M_S$, with $E_*$-equivalences as the
weak equivalences and the $E$-local $S$-modules as the fibrant objects,
see \cite{EKMM97, VIII.1}, but there does not seem to be such an
$E$-local model structure available in the case of $\Cal C_S$.

Let $A \to B$ be a map of $E$-local commutative $S$-algebras, making $B$
a commutative $A$-algebra, and let $G$ be an $E$-locally stably dualizable
group acting continuously on $B$ from the left through commutative
$A$-algebra maps.  For example, $G$ can be a finite discrete group.

Suppose that $A$ is cofibrant as a commutative $S$-algebra, and that $B$
is cofibrant as a commutative $A$-algebra.  The commutative $A$-algebra
$B$ tends not to be cofibrant as an $A$-module, but the smash product
functor $B \wedge_A (-)$ is still homotopically meaningful when applied to
(other) cofibrant commutative $A$-algebras, as explained in \cite{EKMM97,
VII.6}.

Let
$$
i \: A \to B^{hG}
\tag 4.1.1
$$
be the map to the homotopy fixed point $S$-algebra $B^{hG} = F(EG_+,
B)^G$ that is right adjoint to the composite $G$-equivariant map $A
\wedge EG_+ \to A \to B$, collapsing the contractible free $G$-space $EG$
to a point.  Let
$$
h \: B \wedge_A B \to F(G_+, B)
\tag 4.1.2
$$
be the canonical map to the product (cotensor) $S$-algebra $F(G_+, B)$
that is right adjoint to the composite map $B \wedge_A B \wedge G_+
\to B \wedge_A B \to B$, induced by the action $B \wedge G_+ \cong G_+
\wedge B \to B$ of $G$ on $B$, followed by the $A$-algebra multiplication
$B \wedge_A B \to B$ in $B$.

We consider $B \wedge_A B$ and $F(G_+, B)$ as $B$-modules by the
multiplication in the first (left hand) copy of $B$ in $B \wedge_A B$,
and in the target of $F(G_+, B)$.  Then $h$ is a map of $B$-modules.
The group $G$ acts from the left on the second (right hand) copy of $B$
in $B \wedge_A B$, and by right multiplication in the source of $F(G_+,
B)$.  Then $h$ is also a $G$-equivariant map.  These $B$- and
$G$-actions clearly commute, and combine to a left module action by the
group $S$-algebra $B[G]$.

Here is our key definition, which assumes that $E$, $A$, $B$ and $G$ are
as above, and uses the maps $i$ and $h$ just introduced.  We introduce
the related map $j$ in Section~6.1.

\definition{Definition 4.1.3}
We say that $A \to B$ is an {\it $E$-local $G$-Galois extension of
commutative $S$-algebras\/} if the two canonical maps $i \: A \to B^{hG}
= F(EG_+, B)^G$ and $h \: B \wedge_A B \to F(G_+, B)$, formed in the
category of $E$-local $S$-modules, are both weak equivalences.
\enddefinition

The assumption that $A$ and $B$ are $E$-local ensures that $B^{hG}$ and
$F(G_+, B)$ are $E$-local, without any implicit localization.  But $B
\wedge_A B$ formed in $S$-modules needs not be $E$-local, unless $E$ is
smashing.  The condition that $h$ is a weak equivalence in $\Cal M_{S,E}$
amounts to asking that the corresponding map $B \wedge_A B \to F(G_+,
B)$ formed in $\Cal M_S$ is an $E_*$-equivalence, i.e., that $E_*(h)$
is an isomorphism.

\proclaim{Lemma 4.1.4}
Subject to the cofibrancy conditions, the notion of an $E$-local
$G$-Galois extension $A \to B$ is invariant under changes up to weak
equivalence in $A$, $B$ and the stabilized group $S[G] = L_E \Sigma^\infty
G_+$.
\endproclaim

\demo{Proof}
By \cite{EKMM97, VII.6.7} the cofibrancy conditions ensure that the
constructions $A$, $B^{hG}$, $B \wedge_A B$ and $F(G_+, B)$ preserve
weak equivalences in $A$ and $B$, whether implicitly $E$-localized or not.

The natural $E_*$-equivalences $\Sigma^\infty G_+ \to S[G]$ and
$\Sigma^\infty EG_+ \to S[EG]$ induce a (not implicitly localized) map
$$
F_{S[G]}(S[EG], B) \to F_{\Sigma^\infty G_+}(\Sigma^\infty EG_+, B)
\cong F(EG_+, B)^G \,,
$$
which is a weak equivalence when $B$ is $E$-local.  Thus the construction
$B^{hG}$ also preserves weak equivalences in $S[G]$.

Thus the $E$-local Galois conditions, that $G$ is stably dualizable and
the maps $i$ and $h$ are weak equivalences, are invariant under changes
in $A$, $B$ or $G$ that amount to $E$-local weak equivalences of $A$,
$B$ and $S[G]$.
\qed
\enddemo

When $E = S$, so there is no implicit $E$-localization, we may simply say
that $A \to B$ is a $G$-Galois extension, or for emphasis, that $A \to
B$ is a {\it global $G$-Galois extension.}  However, most of the time
we are implicitly working $E$-locally, for a general spectrum $E$, but
omit to mention this at every turn.  Hopefully no confusion will arise.

When $G$ is discrete, we often prefer to write the target $F(G_+, B)$
of $h$ as $\prod_G B$.  When $G$ is finite and discrete, we say that $A
\to B$ is a {\it finite Galois extension.}

\subhead 4.2. The Eilenberg--Mac\,Lane embedding \endsubhead

The Eilenberg--Mac\,Lane functor $H$, which to a commutative ring
$R$ associates a commutative $S$-algebra $HR$ with $\pi_* HR = R$
concentrated in degree $0$, embeds the category of commutative rings
into the category of commutative $S$-algebras.  The two notions of Galois
extension are compatible under this embedding.  For this to make sense,
we must assume that $G$ is finite and that $E = S$.

\proclaim{Proposition 4.2.1}
Let $R \to T$ be a homomorphism of commutative rings, and let $G$
be a finite group acting on $T$ through $R$-algebra homomorphisms.
Then $R \to T$ is a $G$-Galois extension of commutative rings if and
only if the induced map $HR \to HT$ is a global $G$-Galois extension
of commutative $S$-algebras.
\endproclaim

\demo{Proof}
Suppose first that $R \to T$ is $G$-Galois.  Then $T$ is a
finitely generated projective $R$-module by Proposition~2.3.2,
hence flat, so $\Tor^R_s(T, T) = 0$ for $s \ne 0$.  Furthermore,
$T$ is finitely generated projective (of constant rank $1$) as an
$R[G]$-module, by Proposition~2.3.4(c).  There is an isomorphism of
left $R[G]$-modules $R[G] \cong \Hom_R(R[G], R)$, since $G$ is finite,
so $\Ext_{R[G]}^s(R, R[G]) \cong \Ext^s_R(R, R) = 0$ for $s \ne 0$.
Therefore $\Ext_{R[G]}^s(R, T) = 0$ for $s \ne 0$, by the finite
additivity of $\Ext$ in its second argument.

It follows that the homotopy fixed point spectral sequence
$$
E^2_{s,t} = H^{-s}(G; \pi_t HT) = \Ext_{R[G]}^{-s,-t}(R, T)
\Longrightarrow \pi_{s+t}(HT^{hG})
$$
derived from \cite{EKMM97, IV.4.3},
and the K{\"u}nneth spectral sequence
$$
E^2_{s,t} = \Tor^R_{s,t}(T, T)
\Longrightarrow \pi_{s+t}(HT \wedge_{HR} HT)
$$
of \cite{EKMM97, IV.4.2}, both collapse to the origin $s=t=0$.  So
$(HT)^{hG} \simeq H(T^G) = HR$ and $HT \wedge_{HR} HT \simeq H(T \otimes_R
T) \cong H(\prod_G T) \simeq \prod_G HT$ are both weak equivalences.
Thus $HR \to HT$ is a $G$-Galois extension of commutative $S$-algebras.

Conversely, suppose that $HR \to HT$ is $G$-Galois.  Then by the same
spectral sequences $T^G \cong \pi_0(HT^{hG}) \cong \pi_0 HR = R$ and $T
\otimes_R T \cong \pi_0(HT \wedge_{HR} HT) \cong \pi_0(\prod_G HT) \cong
\prod_G T$, so $R \to T$ is a $G$-Galois extension of commutative rings.
\qed
\enddemo

\subhead 4.3. Faithful extensions \endsubhead

Galois extensions of commutative rings are always faithfully flat,
and it will be convenient to consider the corresponding property for
structured ring spectra.  It remains an open problem whether Galois
extensions of commutative $S$-algebras are in fact always faithful,
but we shall verify that this is the case in most of our examples, with
the possible exception of some cases in Section~5.6.

\definition{Definition 4.3.1}
Let $A$ be a commutative $S$-algebra.  An $A$-module $M$ is {\it
faithful\/} if for each $A$-module $N$ with $N \wedge_A M \simeq *$
we have $N \simeq *$.  An $A$-algebra $B$, or $G$-Galois extension $A
\to B$, is said to be faithful if $B$ is faithful as an $A$-module.

A set of $A$-algebras $\{A \to B_i\}_i$ is a {\it faithful cover\/} of $A$
if for each $A$-module $N$ with $N \wedge_A B_i \simeq *$ for every $i$
we have $N \simeq *$.  In particular, a single faithful $A$-algebra $B$
covers $A$ in this sense.
\enddefinition

By the following lemma, this corresponds well to the algebraic notion
of a faithfully flat module \cite{Gre92, 0.1.7}.  Flatness (cofibrancy)
is implicit in our homotopy invariant work, so we only refer to the
faithfulness in our terminology.

\proclaim{Lemma 4.3.2}
Let $M$ be a faithful $A$-module.

(a) A map $f \: X \to Y$ of $A$-modules is a weak equivalence if and only
if $f \wedge 1 \: X \wedge_A M \to Y \wedge_A M$ is a weak equivalence.

(b) A diagram of $A$-modules $X @>f>> Y @>g>> Z$, with a preferred
null-homotopy of $gf$, is a cofiber sequence if and only if $X \wedge_A
M \to Y \wedge_A M \to Z \wedge_A M$, with the associated null-homotopy
of $gf \wedge 1$, is a cofiber sequence.
\endproclaim

\demo{Proof}
(a) Consider the mapping cone $C_f$ of $f$.

(b) Consider the induced map $C_f \to Z$.
\qed
\enddemo

Faithful modules and extensions are preserved under base change,
and are detected by faithful base change.

\proclaim{Lemma 4.3.3}
Let $A \to B$ be a map of commutative $S$-algebras and $M$ a faithful
$A$-module.  Then $B \wedge_A M$ is a faithful $B$-module.
\endproclaim

\demo{Proof}
Let $N$ be a $B$-module such that $N \wedge_B (B \wedge_A M) \simeq *$.
Then $N \wedge_A M \simeq *$, so $N \simeq *$ since $M$ is faithful
over $A$.
\qed
\enddemo

\proclaim{Lemma 4.3.4}
Let $A \to B$ be a faithful map of commutative $S$-algebras and $M$ an
$A$-module such that $B \wedge_A M$ is a faithful $B$-module.  Then $M$
is a faithful $A$-module.
\endproclaim

\demo{Proof}
Let $N$ be an $A$-module such that $N \wedge_A M \simeq *$.  Then $(N
\wedge_A B) \wedge_B (B \wedge_A M) \cong N \wedge_A B \wedge_A M \cong
(N \wedge_A M) \wedge_A B \simeq *$, so $N \wedge_A B \simeq *$ since
$B \wedge_A M$ is faithful over $B$, and thus $N \simeq *$ since $B$
is faithful over $A$.
\qed
\enddemo

\proclaim{Lemma 4.3.5}
For each $G$-Galois extension $R \to T$ of commutative rings, the
induced $G$-Galois extension $HR \to HT$ of commutative $S$-algebras
is faithful.
\endproclaim

\demo{Proof}
Recall that $T$ is faithfully flat over $R$ by Proposition~2.3.4(a).
For each $HR$-module $N$ we have $\pi_*(N \wedge_{HR} HT) \cong \pi_*(N)
\otimes_R T$, by the K{\"u}nneth spectral sequence
$$
E^2_{s,t} = \Tor_{s,t}^R(\pi_*(N), T)
\Longrightarrow \pi_{s+t}(N \wedge_{HR} HT)
$$
and the flatness of $T$.  Therefore $N \wedge_{HR} HT \simeq *$ implies
$\pi_*(N) \otimes_R T = 0$, which in turn implies that $\pi_*(N) = 0$ by
the faithfulness of $T$.  Thus $N \simeq *$ and $HR \to HT$ is faithful.
\qed
\enddemo

\example{Question 4.3.6}
Is every $E$-local $G$-Galois extension $A \to B$ of commutative
$S$-algebras faithful?
\endexample

By Corollary~6.3.4 (or Lemma~6.4.3) the answer is yes when the order
of $G$ is invertible in $\pi_0(A)$, but in some sense this is the less
interesting case.

In the case $E = K(n)$, it is very easy \cite{HSt99, 7.6} to be faithful
over $A = L_{K(n)} S$.

\proclaim{Lemma 4.3.7}
In the $K(n)$-local category, every non-trivial $S$-module is faithful
over $L_{K(n)} S$.
\endproclaim

\demo{Proof}
Let $M$ and $N$ be $K(n)$-local spectra, considered as modules over
$L_{K(n)} S$.  From the K{\"u}nneth formula
$$
K(n)_*(M \wedge_{L_{K(n)} S} N) \cong K(n)_*(M) \otimes_{K(n)_*} K(n)_*(N)
$$
it follows that if $L_{K(n)} (M \wedge_{L_{K(n)} S} N) \simeq *$ then
$K(n)_*(M) = 0$ or $K(n)_*(N) = 0$, since $K(n)_*$ is a graded field.
So if $M$ is non-trivial, we must have $N \simeq *$.  Thus such an $M$
is faithful.
\qed
\enddemo

\head 5. Examples of Galois extensions \endhead

In this chapter we catalog a variety of examples of Galois extensions,
some global and some local, as indicated by the section headings.

\subhead 5.1. Trivial extensions \endsubhead

Let $E$ be any $S$-module and work $E$-locally.  For each cofibrant
commutative $S$-algebra $A$ and stably dualizable group $G$ there is a
{\it trivial $G$-Galois extension\/} from $A$ to $B = F(G_+,A)$, given
by the parametrized {\it diagonal map}
$$
\pi^\# \: A \to F(G_+, A)
$$
that is functionally dual to the collapse map $\pi \: G \to \{e\}$.
Here $G$ acts from the left on $F(G_+, A)$ by right multiplication in
the source.  More precisely, $B$ is a functorial cofibrant replacement
of $F(G_+, A)$ in the category of commutative $A$-algebras, which
inherits the $G$-action by functoriality of the cofibrant replacement.
When $G$ is discrete we can write this extension as $\Delta \: A \to
\prod_G A$.

It is clear that $i \: A \to B^{hG} = F(G_+, A)^{hG}$ is a weak
equivalence, since $(G_+)_{hG} \simeq \{e\}_+$, and that $h \: B
\wedge_A B = F(G_+, A) \wedge_A F(G_+, A) \to F(G_+ \wedge G_+, A)
\cong F(G_+, B)$ is a weak equivalence, since $G$ is stably
dualizable.

The trivial $G$-Galois extension admits an $A$-module retraction
$F(G_+, A) \to A$ functionally dual to the inclusion $\{e\} \to G$, so
$\pi^\# \: A \to F(G_+, A)$ is always faithful.

For any $G$-Galois extension $A \to B$, there is an induced $G$-Galois
extension $B \cong B \wedge_A A \to B \wedge_A B$ (see
Proposition~6.2.1 and Lemma~7.1.3 below), and the map $h \: B \wedge_A
B \to F(G_+, B)$ exhibits an equivalence between this self-induced
extension and the trivial $G$-Galois extension $\pi^\# \: B \to F(G_+,
B)$.

\subhead 5.2. Eilenberg--Mac\,Lane spectra \endsubhead

Let $E = S$.  By Proposition~4.2.1 and Lemma~4.3.5, for each finite
$G$-Galois extension $R \to T$ of commutative rings the induced map of
Eilenberg--Mac\,Lane commutative $S$-algebras $HR \to HT$ is a faithful
$G$-Galois extension.  Proposition~4.2.1 also contains a converse to
this statement.

\subhead 5.3. Real and complex topological $K$-theory \endsubhead

Let $E = S$, and let $KO$ and $KU$ be the real and complex topological
$K$-theory spectra, respectively.  Their connective versions $ko$ and
$ku$ can be realized as the commutative $S$-algebras associated to the
bipermutative topological categories of finite dimensional real and
complex inner product spaces, respectively \cite{May77, VI and VII}.
The periodic commutative $S$-algebras $KO$ and $KU$ are obtained from
these by Bousfield localization, in the $ko$- or $ku$-module categories,
by \cite{EKMM97, VIII.4.3}.

The complexification functor from real to complex inner product spaces
defines maps $c \: ko \to ku$ and $c \: KO \to KU$ of commutative
$S$-algebras, and complex conjugation at the categorical level defines a
$ko$-algebra self map $t \: ku \to ku$ and a $KO$-algebra self map $t \:
KU \to KU$.  Another name for $t$ is the Adams operation $\psi^{-1}$.
Complex conjugation is an involution, so $t^2 = 1$ is the identity in
both cases.  We therefore have an action by $G = \{e, t\} \cong \Z/2$
on $KU$ through $KO$-algebra maps, and can make functorial cofibrant
replacements to keep this property, while making $KO$ cofibrant as a
commutative $S$-algebra and $KU$ cofibrant as a commutative $KO$-algebra.

\proclaim{Proposition 5.3.1}
The complexification map $c \: KO \to KU$ is a faithful $\Z/2$-Galois
extension, i.e., a global quadratic extension.
\endproclaim

See also Example~6.4.4 for more about this extension.

\demo{Proof}
The claim that $i \: KO \to KU^{h\Z/2}$ is a weak equivalence is
well-known to follow from \cite{At66}.  We outline a proof in terms of
the homotopy fixed point spectral sequence
$$
E^2_{s,t} = H^{-s}(\Z/2; \pi_t KU) \Longrightarrow \pi_{s+t}(KU^{h\Z/2})
\,.
$$
Here $\pi_* KU = \Z[u^{\pm1}]$ with $|u| = 2$, $t \in \Z/2$ acts by $t(u)
= -u$ and
$$
E^2_{**} = \Z[a, u^{\pm2}]/(2a)
$$
with $a \in E^2_{-1,2} = H^1(\Z/2; \Z\{u\}) \cong \Z/2$.  A computation
with the Adams $e$-invariant shows that $i$ takes the generator $\eta
\in \pi_1 KO$ to a class represented by $a \in E^{\infty}_{-1,2}$, so
$\eta^3 = 0 \in \pi_3 KO$ implies that $a^3 \in E^2_{-3,6}$ is a
boundary.  The only possibility for this is that $d^3(u^2) = a^3$,
leaving
$$
E^4_{**} = E^\infty_{**} = \Z[a, b, u^{\pm4}]/(2a, a^3, ab, b^2 = 4u^4) \,.
$$
This abutment is isomorphic to $\pi_* KO$, and the graded ring structure
implies that $\pi_*(i)$ is indeed an isomorphism.

To show that $h \: KU \wedge_{KO} KU \to \prod_{\Z/2} KU$ is a weak
equivalence, we use the Bott periodicity cofiber sequence
$$
\Sigma KO @>\eta>> KO @>c>> KU @>\partial>> \Sigma^2 KO
\tag 5.3.2
$$
of $KO$-modules and module maps, up to an implicit weak equivalence
between the homotopy cofiber of $c$ and $\Sigma^2 KO$.  It is the
spectrum level version of the homotopy equivalence $\Omega(U/O) \simeq
\Z \times BU$, and is sometimes stated as an equivalence $KU \simeq KO
\wedge C_\eta$.  Here $\eta$ is given by smashing with the stable Hopf
map $\eta \: S^1 \to S^0$, and $\partial$ is characterized by $\partial
\circ \beta \simeq \Sigma^2 r \: \Sigma^2 KU \to \Sigma^2 KO$, where
$\beta \: \Sigma^2 KU \to KU$ is the Bott equivalence and $r \: KU \to
KO$ is the realification map.  We could write $\partial = \Sigma^2 r
\circ \beta^{-1}$ in $\Cal D_{KO}$.

Inducing~(5.3.2) up along $c \: KO \to KU$, we obtain the upper row in
the following map of horizontal cofiber sequences
$$
\xymatrix{
KU \wedge_{KO} KO \ar[r]^-{1\wedge c} \ar[d]_{\cong}
& KU \wedge_{KO} KU \ar[r]^-{1\wedge\partial} \ar[d]^h
& KU \wedge_{KO} \Sigma^2 KO \ar[d]^{\beta}_{\simeq} \\
KU \ar[r]^-{\Delta}
& \prod_{\Z/2} KU \ar[r]^-{\delta}
& KU
}
\tag 5.3.3
$$
of $KU$-modules and module maps, up to another implicit identification of
the homotopy cofiber of $\Delta$ with $KU$.  Here $h$ is the canonical
map, $\Delta$ is the diagonal inclusion (so the lower row contains the
trivial $\Z/2$-Galois extension of $KU$), $\beta$ is the Bott equivalence
$KU \wedge_{KO} \Sigma^2 KO \cong \Sigma^2 KU \to KU$, and the difference
map $\delta$ is the difference of the two projections from
$\prod_{\Z/2} KU$, indexed by the elements of $\{e, t\} \cong \Z/2$,
written multiplicatively.

The left hand square commutes strictly, since $\Z/2$ acts on $KU$
through $KO$-algebra maps.  To see that the right hand square commutes
up to $KU$-module homotopy, it suffices to prove this after
precomposing with the weak equivalence $1 \wedge \beta \: KU
\wedge_{KO} \Sigma^2 KU \to KU \wedge_{KO} KU$.  To show that the two
resulting $KU$-module maps $KU \wedge_{KO} \Sigma^2 KU \to KU$ are
homotopic, it suffices by adjunction to show that the restricted
$KO$-module maps $\Sigma^2 KU \to KU$ are homotopic.  This is then the
computation
$$
\beta \circ \Sigma^2 c \circ \Sigma^2 r = \delta \circ h \circ (c
\wedge \beta)
$$
in $\Cal D_{KO}$, which follows directly from $\delta \circ h = \mu -
\mu \circ (1 \wedge t) = \mu(1 \wedge (1-t))$ and the well-known relations
$c \circ r = 1 + t$ and $\beta \circ \Sigma^2(1+t) = (1-t) \circ \beta$.

Finally, $c \: KO \to KU$ is faithful.  For if $N$ is a $KO$-module
such that $N \wedge_{KO} KU \simeq *$, then applying $N \wedge_{KO} (-)$
to~(5.3.2) gives a cofiber sequence
$$
\Sigma N @>\eta>> N @>>> N \wedge_{KO} KU \to \Sigma^2 N \,.
$$
The assumption that $N \wedge_{KO} KU \simeq *$ implies that $\eta \:
\Sigma N \to N$ is a weak equivalence.  But $\eta$ is also nilpotent,
since $\eta^4 = 0 \in \pi_4(S)$, so we must have $N \simeq *$.  Therefore
$KU$ is faithful over $KO$.
\qed
\enddemo

The use of nilpotency in this argument may be suggestive of what could
in general be required to answer Question~4.3.6.  We note that the maps
$i \: ko \to ku^{h\Z/2}$ and $h \: ku \wedge_{ko} ku \to \prod_{\Z/2}
ku$ both fail to be weak equivalences.  The homotopy cofiber of $i$ is
$\bigvee_{j<0} \Sigma^{4j} H\Z/2$, and the homotopy cofiber of $h$ is
$H\Z$, as is easily seen by adapting the arguments above.  So $i \: ko
\to ku$ is not Galois.

\subhead 5.4. The Morava change-of-rings theorem \endsubhead

In this section we fix a rational prime $p$ and a natural number~$n$,
and work locally with respect to the $n$-th $p$-primary Morava
$K$-theory $K(n)$.  The work of Devinatz and Hopkins \cite{DH04}
reinterprets the Morava change-of-rings theorem \cite{Mo85, 0.3.3}
as giving a weak equivalence
$$
L_{K(n)} S \simeq E_n^{h\G_n} \,.
$$
We will regard this as a fundamentally important example of a
$K(n)$-local pro-Galois extension $L_{K(n)} S \to E_n$ of commutative
$S$-algebras.  See Definition~8.1.1 for the precise notion of a
pro-Galois extension, which makes most sense after some of the basic
Galois theory has been developed in Chapter~7.

\subsubhead 5.4.1. The Lubin--Tate spectra \endsubsubhead

Recall that $E_n$ is the $n$-th $p$-primary even periodic Lubin--Tate
spectrum, for which
$$
\pi_0(E_n) = \W(\F_{p^n})[[u_1, \dots, u_{n-1}]]
$$
($\W(-)$ denotes the ring of $p$-typical Witt vectors) and $\pi_*(E_n)
= \pi_0(E_n) [u^{\pm1}]$.  Related theories were studied by Morava
\cite{Mo79}, Rudjak \cite{Ru75} and Baker--W{\"u}rgler \cite{BW89},
but in this precise form they seem to have been first considered by
Hopkins and Miller \cite{HG94}, \cite{Re98}.

The height~$n$ Honda formal group law $\Gamma_n$ is defined over $\F_p$
and is characterized by its $p$-series $[p]_n(x) = x^{p^n}$.  Its
Lubin--Tate deformation $\widetilde\Gamma_n$ over $\F_{p^n}$ is the
universal formal group law over a complete local ring with residue
field an extension of~$\F_{p^n}$, whose reduction to the residue field
equals the corresponding extension of $\Gamma_n$.  In this case the
universal complete local ring equals $\pi_0(E_n)$, with maximal ideal
$(p, u_1, \dots, u_{n-1})$ and residue field $\F_{p^n}$.  The
Lubin--Tate spectrum $E_n$ is (at first) the $K(n)$-local complex
oriented commutative ring spectrum that represents the resulting
Landweber exact homology theory $(E_n)_*(X) = \pi_*(E_n)
\otimes_{\pi_*(MU)} MU_*(X)$.

More generally, we can consider $\Gamma_n$ as a formal group law over
the algebraic closure $\bar\F_p$ of $\F_p$.  Its universal deformation
is then defined over the complete local ring
$$
\pi_0(E_n^{nr}) = \W(\bar\F_p)[[u_1, \dots, u_{n-1}]]
$$
and there is a similar $K(n)$-local complex oriented commutative ring
spectrum $E_n^{nr}$ with $\pi_*(E_n^{nr}) = \pi_0(E_n^{nr}) [u^{\pm1}]$.
The superscript ``nr'' is short for ``non ramifi{\'e}e'', indicating
that $\W(\bar\F_p)$ is the $p$-adic completion of the maximal unramified
extension $\colim_f \W(\F_{p^f})$ of $\W(\F_p) = \Z_p$.  (The infinite
product defining $p$-typical Witt vectors only commutes with the colimit
over $f$ after completion.)

\subsubhead 5.4.2. The extended Morava stabilizer group \endsubsubhead

The profinite Morava stabilizer group $\SS_n = \Aut(\Gamma_n /
\F_{p^n})$ of automorphisms defined over $\F_{p^n}$ of the formal group
law $\Gamma_n$ (see \cite{Ra86, \S A2.2, \S6.2}), and the finite Galois
group $\Gal = \Gal(\F_{p^n}/\F_p) \cong \Z/n$ of the extension $\F_p
\subset \F_{p^n}$, both act on the universal deformation
$\widetilde\Gamma_n$, and thus on $\pi_*(E_n)$, by the universal
property.  These actions combine to one by the profinite semi-direct
product $\G_n = \SS_n \rtimes \Gal$.  By the Hopkins--Miller
\cite{Re98} and Goerss--Hopkins theory \cite{GH04, \S7} the ring
spectrum $E_n$ admits the structure of a commutative $S$-algebra, up to
a contractible choice.  Furthermore, the extended Morava stabilizer
group $\G_n$ acts on $E_n$ through commutative $S$-algebra maps, again
up to contractible choice.  However, these actions through commutative
$S$-algebras do not take into account the profinite topology on $\G_n$,
but rather treat $\G_n$ as a discrete group.

It is known by recent work of Daniel G.~Davis \cite{Da:h}, that the
profinite group $\G_n$ acts continuously on~$E_n$ in the category of
$K(n)$-local $S$-modules, but only when $E_n$ is reconsidered as a
pro-object of discrete $\G_n$-module spectra, where the terms have
coefficient groups of the form $\pi_*(E_n)/I_k$ for a suitable
descending sequence of ideals $\{I_k\}_k$ with $\bigcap_k I_k = 0$.
Presently, this kind of limit presentation is not available in the
context of commutative $S$-algebras.  Hopkins has suggested that a
weaker form of structured commutativity, in terms of pro-spectra, may
instead be available.

More generally, the Morava stabilizer group $\SS_n$ and the absolute
Galois group $\Gal(\bar\F_p/\F_p) \cong \hat\Z$ (the Pr{\"u}fer ring)
of $\F_p$ both act on the universal deformation of~$\Gamma_n$ over
the algebraic closure $\bar\F_p$, and thus on $\pi_*(E_n^{nr})$ by the
universal property.  These combine to an action by the profinite group
$$
\G_n^{nr} = \SS_n \rtimes \hat\Z \,.
$$
Note that the (conjugation) action by $\hat\Z$ on $\SS_n$ factors through
the quotient $\hat\Z \to \Z/n = \Gal$, since all the automorphisms of the
height~$n$ Honda formal group law are already defined over $\F_{p^n}$
\cite{Ra86, A2.2.20(a)}.  The Goerss--Hopkins theory cited above again
implies that $E_n^{nr}$ is a commutative $S$-algebra, and the extended
Morava stabilizer group $\G_n^{nr}$ acts on $E_n^{nr}$ through commutative
$S$-algebra maps, up to contractible choices.

\subsubhead 5.4.3. Intermediate $S$-algebras \endsubsubhead

In the Galois theory for fields, the intermediate fields $F \subset E
\subset \bar F$ correspond bijectively (via $E = (\bar F)^K$ and $K =
G_E$) to the closed subgroups $K \subset G_F$ of the absolute Galois
group with the Krull topology, and the finite field extensions $F
\subset E$ correspond to the open subgroups $U \subset G_F$.  Note that
in this topology, the open subgroups are exactly the closed subgroups
of finite index.  Furthermore, $G_F$ acts continuously on $\bar F$ with
the discrete topology, so $\bar F$ is the union over the open subgroups
$U$ of the fixed fields $(\bar F)^U$.

By analogy, it is desirable to construct intermediate $K(n)$-local
commutative $S$-algebras $E_n^{hK}$ for every closed subgroup $K
\subset \G_n$ in the profinite topology.  If $E_n$ were a discrete
$\G_n$-module spectrum, this could be done by the usual definition
$E_n^{hK} = F(EK_+, E_n)^K$, and indeed, for finite (and thus discrete)
subgroups $K \subset \G_n$ the restricted $K$-action is continuous,
$E_n$ is a discrete $K$-module spectrum and $E_n^{hK}$ can well be
defined in this way.  The maximal finite subgroups $M \subset \G_n$
were classified by Hewett \cite{He95, 1.3, 1.4}.  When $M$ is unique up
to conjugacy, $E_n^{hM}$ is known as the {\it $n$-th higher real
$K$-theory\/} spectrum~$EO_n$ of Hopkins and Miller.  Such uniqueness
holds for $p$ odd when $n = (p-1)k$ with $k$ prime to $p$, and for
$p=2$ when $n = 2k$ with $k$ an odd natural number, by {\it loc.~cit.}

However, as recalled in the previous subsection, the spectrum $E_n$ is
not itself a discrete $\G_n$-spectrum, but only an inverse limit of
such, i.e., a pro-discrete $\G_n$-spectrum.  The homotopy invariant way
to form homotopy fixed points of such objects is to take the ordinary
continuous homotopy fixed points for the profinite group acting
discretely at each stage in the limit system, and then to pass to the
homotopy limit, if desired.  Note that the formation of continuous
homotopy fixed points for profinite groups acting on discrete modules
involves a colimit indexed over the finite quotients of the profinite
group, and does therefore not generally commute with limits.  This
procedure describes the approach of \cite{Da:h}, but it only exhibits
the homotopy fixed point spectrum $E_n^{h\G_n}$ as a module spectrum,
and not as an algebra spectrum, precisely because we do not know how to
realize $E_n$ as a pro-object of $\G_n$-discrete associative or
commutative $S$-algebras.

Devinatz and Hopkins circumvent this problem by defining $E_n^{h\G_n}$,
and more generally $E_n^{hU}$ for each open subgroup $U \subset \G_n$,
in a ``synthetic'' way \cite{DH04, Thm.~1}, as the totalization of a
suitably rigidified cosimplicial diagram, to obtain a $K(n)$-local
commutative $S$-algebra of the desired homotopy type.  In particular,
$E_n^{h\G_n} \simeq L_{K(n)} S$.  (See Section~8.2 for further
discussion of the kind of cosimplicial diagram involved, namely the
Amitsur complex.)  For closed subgroups $K \subset \G_n$ they then
define \cite{DH04, Thm.~2}
$$
E_n^{hK} = L_{K(n)} (\colim_i E_n^{hU_iK})
$$
where $\{U_i\}_{i=0}^\infty$ is a fixed descending sequence of open normal
subgroups in $\G_n$ with $\bigcap_{i=0}^\infty U_i = \{e\}$, and the
colimit is the homotopy colimit in commutative $S$-algebras.  For finite
subgroups $K \subset \G_n$ the synthetic construction agrees \cite{DH04,
Thm.~3} with the ``natural'' definition of $E_n^{hK}$ as $F(EK_+, E_n)^K$.

Ethan Devinatz \cite{De05} then proceeds to compare the commutative
$S$-algebras $E_n^{hK}$ and $E_n^{hH}$ for closed subgroups $K$ and $H$
of $\G_n$ with $H$ normal in $K$.  There is a well-defined action by the
quotient group $K/H$ on $E_n^{hH}$ through commutative $S$-algebra maps,
in the $K(n)$-local category \cite{De05, \S3}.

\proclaim{Theorem 5.4.4 (Devinatz--Hopkins)}
(a)
For each pair of closed subgroups $H \subset K \subset \G_n = \Bbb
S_n \rtimes \Gal$ with $H$ normal and of finite index in $K$, the map
$E_n^{hK} \to E_n^{hH}$ is a $K(n)$-local $K/H$-Galois extension.

(b)
In particular, for each finite subgroup $K \subset \G_n$ the map
$E_n^{hK} \to E_n$ is a $K(n)$-local $K$-Galois extension.

(c)
Likewise, for each open normal subgroup $U \subset \G_n$ (necessarily
of finite index) the map
$$
L_{K(n)} S \to E_n^{hU}
$$
is a $K(n)$-local $\G_n/U$-Galois extension.

(d)
A choice of a descending sequence $\{U_i\}_i$ of open normal
subgroups of $\G_n$, with $\bigcap_i U_i = \{e\}$, exhibits
$$
L_{K(n)} S \to E_n
$$
as a $K(n)$-local pro-$\G_n$-Galois extension, in view of the weak
equivalence
$$
L_{K(n)} (\colim_i E_n^{hU_i}) @>\simeq>> E_n \,.
$$
\endproclaim

\demo{Proof}
(a) Let $A = E_n^{hK}$, $B = E_n^{hH}$ and $G = K/H$ (which is finite
and discrete).  By \cite{De05, Prop.~2.3, Thm.~3.1 and Thm.~A.1} the
homotopy fixed point spectral sequence for $\pi_*(B^{hG})$ agrees with a
strongly convergent $K(n)_*$-local Adams spectral sequence converging to
$\pi_*(A)$.  So $i \: A \to B^{hG}$ is a weak equivalence.  By \cite{De05,
Cor.~3.9} the natural map $h \: L_{K(n)} (B \wedge_A B) \to F(G_+, B)$
induces an isomorphism on homotopy groups.

Parts~(b) and (c) are special cases of (a). Part~(d) is contained in
\cite{DH04, Thm.~3(i)}.
\qed
\enddemo

It would be nice to extend the statement of this theorem to the case when
$H$ is normal and closed, but not necessarily of finite index, in $K$.

For $n=2$ and $p=2$, the Morava stabilizer group $\SS_2$ is the group of
units in the maximal order in the quaternion algebra $\Q_2\{1, i, j, k\}$,
and its maximal finite subgroup is the binary tetrahedral group $\hat A_4
= Q_8 \rtimes \Z/3$ of order~24, containing the quaternion group $Q_8 =
\{\pm 1, \pm i, \pm j, \pm k\}$ and the 16 other elements $(\pm 1 \pm
i \pm j \pm k)/2$.  See \cite{CF67, pp.~137--138}, \cite{Ra86, 6.3.27}.
The maximal finite subgroup of $\G_2$ is $G_{48} = \hat A_4 \rtimes \Z/2$,
and $EO_2 = E_2^{hG_{48}}$ is the $K(2)$-localization of the connective
spectrum $eo_2$ with $H^*(eo_2; \F_2) \cong A/\!/A_2$ as a module over
the Steenrod algebra, which is related to the topological modular forms
spectrum $\tmf$ \cite{Hop02, \S3.5}.

\proclaim{Proposition 5.4.5}
At $p=2$, the natural map $EO_2 \to E_2$ is a $K(2)$-local faithful
$G_{48} = \hat A_4 \rtimes \Z/2$-Galois extension.
\endproclaim

\demo{Proof}
This follows from Theorem~5.4.4(b) above and Proposition~5.4.9(b) below,
but we would also like to indicate a direct proof of faithfulness,
using results of Hopkins and Mahowald \cite{HM98}.  There is a finite
CW spectrum $C_\gamma$ obtained as the mapping cone of a map
$$
\gamma \: \Sigma^5 C_\eta \wedge C_\nu \to C_\eta \wedge C_\nu \,,
$$
such that $H^*(C_\gamma; \F_2) \cong DA(1) \cong A(2)/\!/E(2)$ is the
``double'' of $A(1) = \langle Sq^1, Sq^2 \rangle$.  The spectrum
$C_\gamma$ and the self-map $\gamma$ can be obtained by a construction
analogous to that of the spectrum $A_1$ and the map $v_1 \: \Sigma^2 Y
\to Y$ in \cite{DM81, pp.~619--620}, but replacing all the real
projective spaces occurring there by complex projective spaces.
Furthermore, there is a weak equivalence $eo_2 \wedge C_\gamma \simeq
BP\langle2\rangle$ that realizes the isomorphism $A/\!/A(2) \otimes
A(2)/\!/E(2) \cong A/\!/E(2) \cong H^*(BP\langle2\rangle; \F_2)$.
Applying $K(2)$-localization yields
$$
EO_2 \wedge C_\gamma \simeq \widehat{E(2)} \,,
$$
in the notation of~5.4.7, using that $BP\langle2\rangle \to v_2^{-1}
BP\langle2\rangle = E(2)$ is a $K(2)_*$-equivalence.  Since $\eta$,
$\nu$ and $\gamma$ are all nilpotent (for $\eta \in \pi_1(S)$ and $\nu
\in \pi_3(S)$ this is well-known; for $\gamma$ it can be deduced from
the Devinatz--Hopkins--Smith nilpotence theorem \cite{DHS88, Cor.~2}),
it follows as in the proof of Proposition~5.3.1 that $EO_2 \to
\widehat{E(2)}$ is faithful.  And $\widehat{E(2)} \to E_2$ is faithful
by the elementary Proposition~5.4.9(a).
$$
\xymatrix{
& E_2^{\Gal} \ar[r]^{\Z/2} & E_2 \\
EO_2 \ar[r] \ar[ur]^{\hat A_4}
	& {}\widehat{E(2)} \ar[u]^{\F_4^*} \ar[ur]_{\Sigma_3}
}
$$
\qed
\enddemo

\subsubhead 5.4.6. Adjoining roots of unity \endsubsubhead

Including the maximal unramified extensions into this picture, we have
the following diagram of $K(n)$-local extensions.  The groups label Galois
(or pro-Galois) extensions.
$$
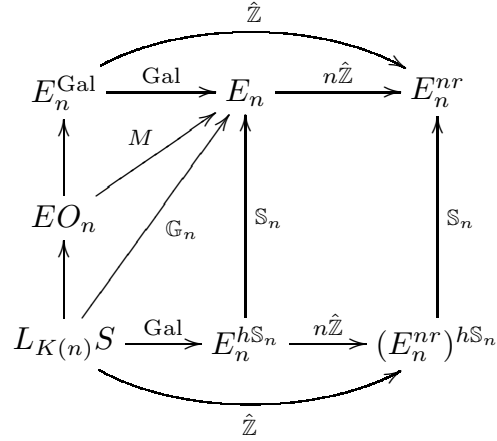
\xymatrix{
E_n^{\Gal} \ar[r]^-{\Gal} \ar@(ur,ul)[rr]^-{\hat\Z}
	& E_n \ar[r]^-{n\hat\Z} & E_n^{nr} \\
EO_n \ar[u] \ar[ur]^-M \\
L_{K(n)} S \ar[r]^-{\Gal} \ar@(dr,dl)[rr]_-{\hat\Z} \ar[u] \ar[uur]_-{\G_n}
	& E_n^{h\SS_n} \ar[r]^-{n\hat\Z} \ar[uu]_{\SS_n}
	& (E_n^{nr})^{h\SS_n} \ar[uu]_{\SS_n}
}
$$
The maximal extension $L_{K(n)} S \to E_n^{nr}$ is $K(n)$-locally
pro-$\G_n^{nr}$-Galois.  

The $\hat\Z$-extension along the bottom is that obtained by adjoining
all roots of unity of order prime to $p$ to the $p$-complete commutative
$S$-algebra $L_{K(n)} S$.  We might write $E_n^{h\SS_n} = L_{K(n)}
S(\mu_{p^n{-}1})$ and $(E_n^{nr})^{h\SS_n} = L_{K(n)} S(\mu_{\infty,p})$,
where $\mu_m$ denotes the group of $m$-th order roots of unity and
$\mu_{\infty,p} = \colim_{p \nmid m} \mu_m$ denotes the group of all
roots of unity of order prime to $p$.  Note that in the latter case,
the infinite colimit of spectra must be implicitly $K(n)$-completed.
Similarly, $E_n = E_n^{\Gal}(\mu_{p^n{-}1})$ and $E_n^{nr} =
E_n^{\Gal}(\mu_{\infty,p})$.

The process of adjoining $m$-th roots of unity makes sense when applied
to a $p$-local commutative $S$-algebra~$A$, for $p \nmid m$, following
Roland Schw{\"a}nzl, Rainer Vogt and Waldhausen \cite{SVW99}, since
$A(\mu_m)$ can be obtained from the group $A$-algebra $A[C_m] = A
\wedge C_{m+}$ of the cyclic group of order~$m$ by localizing with
respect to a $p$-locally defined idempotent.  Likewise, adjoining an
$m$-th root of unity to a $p$-complete commutative $S$-algebra $A$, for
$m = p^f-1$, can be achieved by localizing with respect to a further
idempotent.  The situation is analogous to how $\Q_p \otimes_{\Q}
\Q(\mu_m)$ splits as a product of copies of $\Q_p(\mu_m)$, when $m =
p^f-1$.  For more on the process of adjoining roots of unity to
$S$-algebras, see \cite{La03, 3.4} in the associative case and
\cite{BR:r, 2.2.5 and~2.2.8} in the commutative case.

These observations may justify thinking of the projection $d \: \G_n^{nr}
= \SS_n \rtimes \hat\Z \to \hat\Z$ as the {\it degree map\/} of a
$K(n)$-local class field theory for structured ring spectra \cite{Ne99,
\S IV.4}.

\subsubhead 5.4.7. Faithfulness \endsubsubhead

Let $\widehat{E(n)} = L_{K(n)} E(n)$ be the $K(n)$-localization of the
Johnson--Wilson spectrum $E(n)$ from~3.2.2, called Morava $E$-theory
in \cite{HSt99}.  By \cite{BW89, 4.1} or \cite{HSt99, \S1.1, 5.2} it
has coefficients
$$
\pi_* \widehat{E(n)} = \Z_{(p)}[v_1, \dots, v_{n-1},
v_n^{\pm1}]^\wedge_{I_n} \,,
$$
where $I_n = (p, v_1, \dots, v_{n-1})$.  The spectrum $\widehat{E(n)}$
was proved to be an associative $S$-algebra in \cite{Bak91}, and is in
fact a commutative $S$-algebra by the homotopy fixed point description
in Proposition~5.4.9(a) below.

\proclaim{Theorem 5.4.8 (Hovey--Strickland)}
$L_{K(n)} S$ is contained in the thick subcategory of $K(n)$-local spectra
generated by $\widehat{E(n)}$, so $\widehat{E(n)}$ is a faithful $L_{K(n)}
S$-module in the $K(n)$-local category.
\endproclaim

\demo{Proof}
The first claim is contained in the proof of \cite{HSt99, 8.9}, which
relies heavily on the construction by Jeff Smith of a suitable finite
$p$-local spectrum $X$, as explained in \cite{Ra92, \S8.3}.  The second
claim follows from the first, but also much more easily from Lemma~4.3.7.
\qed
\enddemo

\proclaim{Proposition 5.4.9}
(a)
The $K(n)$-local Galois extension $\widehat{E(n)} \to E_n$ and the
$K(n)$-local pro-Galois extension $L_{K(n)} S \to E_n$ are both faithful.

(b)
For each pair of closed subgroups $H \subset K \subset \G_n$, with $H$
normal and of finite index in $K$, the $K(n)$-local $K/H$-Galois extension
$E_n^{hK} \to E_n^{hH}$ is faithful.
\endproclaim

\demo{Proof}
(a)
There is a finite subgroup $\F_{p^n}^*$ of $\SS_n$ such that for $K =
\F_{p^n}^* \rtimes \Gal$ we have $\widehat{E(n)} \simeq E_n^{hK}$.  In
more detail, $\SS_n$ contains the unit group $\W(\F_{p^n})^*$ \cite{Ra86,
A2.2.17}, whose torsion subgroup reduces isomorphically to $\F_{p^n}^*$.
For an element of finite order $\omega \in \W(\F_{p^n})^*$, with mod~$p$
reduction $\bar\omega \in \F_{p^n}^*$, the linear formal power series
$g(x) = \bar\omega x$ defines an automorphism of $\Gamma_n$, i.e., an
element $g \in \SS_n$, which acts on $\pi_*(E_n)$ by $g(u) = \omega u$
and $g(u u_k) = \omega^{p^k} u u_k$ for $1 \le k < n$ by \cite{DH95, 3.3,
4.4}, leaving $v_n = u^{1-p^n}$ and $v_k = u^{1-p^k} u_k$ invariant.
Thus $\pi_* E_n^{\Gal} = \Z_p[[u_1, \dots, u_{n-1}]] [u^{\pm1}]$ and
$\pi_* E_n^{hK}$ is the $I_n$-adic completion of $\pi_* E(n)$.

Then for any spectrum $X$,
$$
(E_n)^\vee_*(X) \cong \pi_* E_n \otimes_{\pi_*\widehat{E(n)}}
\widehat{E(n)}{}^\vee_*(X)
$$
with $\pi_* E_n$ a free module of rank $|K| = (p^n{-}1)n$ over
$\pi_*\widehat{E(n)}$.  Here we are using the notation $(E_n)^\vee_*(X)
= \pi_* L_{K(n)} (E_n \wedge X)$, and similarly for $\widehat{E(n)}$,
of \cite{HSt99, 8.3}.  It follows easily from this formula that
$\widehat{E(n)} \to E_n$ is faithful in the $K(n)$-local category.

In combination with~5.4.8 this also shows that the composite extension
$L_{K(n)} S \to \widehat{E(n)} \to E_n$ is faithful, but Lemma~4.3.7
provides a much easier argument.

(b)
The second result follows by faithful base change along $\phi \:
L_{K(n)} S \to E_n$.  There is a commutative diagram (for $H$ and $K$
as in the statement)
$$
\xymatrix{
E_n^{hH} \ar[r] & L_{K(n)} (E_n \wedge E_n^{hH}) \\
E_n^{hK} \ar[r]^-{1\wedge\phi} \ar[u]^{\psi}
	& L_{K(n)} (E_n \wedge E_n^{hK}) \ar[u]_{1\wedge\psi} \\
L_{K(n)} S \ar[r]^-{\phi} \ar[u] & E_n \ar[u] \\
}
$$
where the squares are pushouts in the category of $K(n)$-local commutative
$S$-algebras.  By the Morava change-of-rings theorem and \cite{DH04,
Thm.~1(iii)},
$$
\pi_* L_{K(n)} (E_n \wedge E_n^{hH}) \cong \Map(\G_n/H, \pi_* E_n)
$$
and
$$
\pi_* L_{K(n)} (E_n \wedge E_n^{hK}) \cong \Map(\G_n/K, \pi_* E_n) \,.
$$
See also the proof of Theorem~7.2.3 below.  Here $\Map$ denotes the
unbased continuous maps with respect to the profinite topologies on
$\G_n/H$, $\G_n/K$ and $\pi_* E_n$ (in each degree).  Note that $K/H$
is a finite group acting freely on the Hausdorff space $\G_n/H$, with
orbit space $\G_n/K$, so $\pi \: \G_n/H \to \G_n/K$ is a regular
$K/H$-covering space.  We claim that it admits a continuous section
$\sigma \: \G_n/K \to \G_n/H$, so that there is a homeomorphism $K/H
\times \G_n/K \to \G_n/H$, and
$$
\Map(\G_n/H, \pi_* E_n) \cong \prod_{K/H} \Map(\G_n/K, \pi_* E_n) \,.
$$
Thus $\pi_* L_{K(n)} (E_n \wedge E_n^{hH})$ is a free module of rank
$|K/H|$ over $\pi_* L_{K(n)} (E_n \wedge E_n^{hK})$, so that $L_{K(n)}
(E_n \wedge E_n^{hH})$ is faithful over $L_{K(n)} (E_n \wedge E_n^{hK})$.
The map $1\wedge\phi$ is obtained by base change from $\phi$, which is
faithful by~(a), and is therefore faithful by Lemma~4.3.3, so $\psi \:
E_n^{hK} \to E_n^{hH}$ is faithful by Lemma~4.3.4.

It remains to verify the claim.  Let $\{U_i\}_{i=0}^\infty$ be a
descending sequence of open normal subgroups of $\G_n$, with trivial
intersection as above.  Then $U_i H$ is normal of finite index in $U_i
K$, $K/H$ surjects to $U_i K/U_i H$ and there is a regular covering
space $\pi_i \: \G_n/U_i H \to \G_n/U_i K$, for each $i$.  We have the
following commutative diagram for $i < j$:
$$
\xymatrix{
K/H \ar[r] \ar[d] & U_jK/U_jH \ar[r] \ar[d] & U_iK/U_iH \ar[d] \\
\G_n/H \ar[r] \ar[d]_{\pi} & \G_n/U_jH \ar[r] \ar[d]^{\pi_j}
	& \G_n/U_iH \ar[d]^{\pi_i} \\
\G_n/K \ar[r] & \G_n/U_jK \ar[r] & \G_n/U_iK
}
$$
Since $K/H$ is finite, the surjections $U_jK/U_jH \to U_iK/U_iH$
are isomorphisms for all sufficiently large $i$ and $j$, say for $i,j
\ge i_0$, and then $\pi_j$ is the pullback of $\pi_i$ along $\G_n/U_jK
\to \G_n/U_iK$.  Thus any choice of section $\sigma_i$ to $\pi_i$ pulls
back to a section $\sigma_j$ of $\pi_j$, so that the composite maps
$\G_n/K \to \G_n/U_iK \to \G_n/U_iH$ are compatible for all $i\ge i_0$.
Their limit defines the continuous section $\sigma \: \G_n/K \to \G_n/H$.
\qed
\enddemo

\subhead 5.5. The $K(1)$-local case \endsubhead 

When $n=1$, the discussion in Section~5.4 reduces to more classical
statements about variants of topological $K$-theory, which we now make
explicit, together with a comparison to the even more classical
arithmetic theory of abelian extensions of $\Q_p$ and $\Q$.

\subsubhead 5.5.1. $p$-complete topological $K$-theory \endsubsubhead

Mod~$p$ complex topological $K$-theory, with $\pi_*(KU/p) =
\F_p[u^{\pm1}]$, splits as
$$
KU/p \simeq \bigvee_{i=0}^{p-2} \Sigma^{2i} K(1)
$$
where $\pi_* K(1) = \F_p[v_1^{\pm1}]$.  Bousfield $K(1)$-localization
equals Bousfield $KU/p$-localization, which in turn equals Bousfield
$KU$-localization followed by $p$-adic completion: $L_{K(1)} X = L_{KU/p}
X = (L_{KU} X)^\wedge_p$ \cite{Bo79, 2.11}.

The height~$1$ Honda formal group law over $\F_p$ is the multiplicative
one: $\Gamma_1(x, y) = x + y + xy$, its universal deformation
$\widetilde\Gamma_1$ is the multiplicative formal group law over
$\Z_p$, and the Lubin--Tate spectrum $E_1$ equals $p$-completed complex
topological $K$-theory $KU^\wedge_p$ with $\pi_*(KU^\wedge_p) =
\Z_p[u^{\pm1}]$.  The Morava stabilizer group $\G_1 = \SS_1$ is the
group of $p$-adic units $\Z_p^*$, with its profinite topology, and $k
\in \Z_p^*$ acts on the commutative $S$-algebra $KU^\wedge_p$ by the
$p$-adic Adams operation
$$
\psi^k \: KU^\wedge_p \to KU^\wedge_p \,.
$$
On homotopy, $\psi^k(u) = ku$.

\subsubhead 5.5.2. Subalgebras \endsubsubhead

The homotopy fixed point spectrum $E_1^{h\G_n} = (KU^\wedge_p)^{h\Z_p^*}$
is the $p$-complete (non-connective) image-of-$J$ spectrum $L_{K(1)}S =
J^\wedge_p$, defined for $p=2$ by the fiber sequence
$$
J^\wedge_2 \to KO^\wedge_2 @>\psi^3-1>> KO^\wedge_2
$$
and for $p$ odd by the fiber sequence
$$
J^\wedge_p \to KU^\wedge_p @>\psi^r-1>> KU^\wedge_p
$$
for $r$ a topological generator of $\Z_p^*$.  These identifications
of the $p$-completed $KU$-localization of $S$ with $J^\wedge_p$ are
basically due to Mark Mahowald and Haynes Miller \cite{Bo79, 4.2},
respectively.  (Adams--Baird and Ravenel went on to identify the $p$-local
$KU$-localization of $S$, see \cite{Bo79, 4.3}.)

The Morava stabilizer group $\SS_1 = \Z_p^*$ is isomorphic to the Galois
group of the maximal (totally ramified) $p$-cyclotomic extension $\Q_p
\subset \Q_p(\mu_{p^\infty})$, so the classification of intermediate
commutative $S$-algebras $J_p \to C \to KU^\wedge_p$ of the form
$C = (KU^\wedge_p)^{hK}$ for $K$ closed in $\Z_p^*$ is identical
to the classification of intermediate fields $\Q_p \subset E \subset
\Q_p(\mu_{p^\infty})$.  In this way $J^\wedge_p \to KU^\wedge_p$ provides
a $K(1)$-local ``realization'' of the $K(0)$-local extension $\Q_p \to
\Q_p(\mu_{p^\infty})$.  There are similar $K(n)$-local realizations
of the form $L_{K(n)} S \to E_n^{hK}$, when $K$ is the kernel of the
determinant/abelianization homomorphism $\G_n \to \G_n^{ab} \to \Z_p^*$
\cite{Ra86, 6.2.6(b)}.

When $p=2$, $\Z_2^* \cong \Z_2 \times \Z/2$, where $\Z_2 \cong 1 +
4\Z_2$ is open of index~$2$, and $\Z/2 \cong \{\pm1\} \subset \Z_2^*$
is closed.  There are three different subgroups of index~$2$, namely
the topologically generated subgroups $\langle3\rangle$,
$\langle5\rangle$ and $\langle-1,9\rangle$.  The first of these
corresponds to the complex image-of-$J$ spectrum $JU^\wedge_2 =
(KU^\wedge_2)^{h\langle3\rangle}$ given by the fiber sequence
$$
JU^\wedge_2 \to KU^\wedge_2 @>\psi^3-1>> KU^\wedge_2 \,,
$$
and there is a $K(1)$-local (quadratic) $\Z/2$-Galois extension $c \:
J^\wedge_2 \to JU^\wedge_2$, which is compatible with the
complexification map $c \: KO^\wedge_2 \to KU^\wedge_2$.  See
Example~6.2.2 for more on this quadratic extension.  The closed
subgroup $\Z/2$ of $\Z_2^*$ corresponds to $2$-complete real
$K$-theory: $(KU^\wedge_2)^{h\Z/2} \simeq KO^\wedge_2$.

When $p$ is odd, $\Z_p^* \cong \Z_p \times \F_p^*$ is pro-cyclic.  Let
$r \in \Z_p^*$ be a topological generator, chosen to be a natural
number.  Then $\Z_p^*$ has a unique open subgroup $\langle r^n \rangle$
of index~$n$, for each integer $n$ of the form $n = p^e d$ with $e\ge0$
and $d \mid p-1$.  In addition, it has the closed subgroups that appear
as subgroups of $\F_p^*$.  In particular, $\Z_p^*$ has an open subgroup
$\Z_p \cong 1 + p\Z_p$ of index $(p-1)$, and a closed subgroup $\F_p^*
\subset \Z_p^*$.  The latter corresponds to the $p$-complete Adams
summand $L^\wedge_p = (KU^\wedge_p)^{h\F_p^*}$ with $\pi_*(L^\wedge_p)
= \Z_p[v_1^{\pm1}]$.  There are $K(1)$-local $\F_p^*$-Galois extensions
$J^\wedge_p \to (KU^\wedge_p)^{h\Z_p}$ and $L^\wedge_p \to
KU^\wedge_p$.  Let us write $F\psi^{r^n} = (KU^\wedge_p)^{h\langle r^n
\rangle}$ for the homotopy fixed point spectrum of $\psi^{r^n}$, which
is equivalent to the homotopy fiber of $\psi^{r^n}-1$.  Then there is a
$K(1)$-local $\Z/n$-Galois extension
$$
J^\wedge_p = F\psi^r \to F\psi^{r^n}
$$
for each integer $n = p^e d$ with $d \mid p-1$, as above.

\subsubhead 5.5.3. Extensions \endsubsubhead

Incorporating the roots of unity of order prime to~$p$, we have the
following diagram
$$
\xymatrix{
KU^\wedge_p \ar[r]^-{\hat\Z} & KU^\wedge_p(\mu_{\infty,p}) \\
J^\wedge_p \ar[r]^-{\hat\Z} \ar[u]^{\Z_p^*} \ar[ur]^{\G_1^{nr}}
	& J^\wedge_p(\mu_{\infty,p}) \ar[u]_{\Z_p^*} \\
}
$$
with $E_1^{nr} = KU^\wedge_p(\mu_{\infty,p})$.  Here the maximal Galois
group $\G_1^{nr} = \Z_p^* \times \hat\Z$ is abelian, since $\hat\Z$
acts trivially on~$\SS_1 = \Z_p^*$.  It provides a $K(1)$-local
realization of the Galois group of the maximal abelian extension $\Q_p
\to \Q_p(\mu_\infty)$.

It also appears to be possible to fit the various rational primes
together, so as to obtain $KU$-local realizations of the abelian extensions
of the rational field $\Q$ itself.  The Galois group $G = \hat\Z^*$
of the maximal abelian extension $\Q \to \Q(\mu_\infty)$ contains the
Galois group of $\Q_p \to \Q_p(\mu_\infty)$ as the decomposition group
$D_p$ of the prime ideal~$(p)$.  Let $Z_p = \Q(\mu_\infty)^{D_p}$ be
the corresponding decomposition field \cite{Ne99, I.9.2}.
$$
\Q @>G/D_p>> Z_p @>D_p>> \Q(\mu_\infty) \,.
$$
After base change along $\Q \to \Q_p$ there are weak product splittings
\cite{Ne99, II.8.3}
$$
\Q_p \otimes_{\Q} Z_p
	\cong \wprod_{G/D_p} \Q_p
\qquad\text{and}\qquad
\Q_p \otimes_{\Q} \Q(\mu_\infty)
	\cong \wprod_{G/D_p} \Q_p(\mu_\infty)
\,,
$$
i.e., as colimits of the products over the finite quotients of
$$
G/D_p = \hat\Z^*/(\Z_p^* \times \hat\Z) \cong (\prod_{\ell\ne p}
\Z_\ell^*)/\hat\Z \,.
$$
In the latter profinite quotient, the unit of $\hat\Z$ maps diagonally
to the class of $p$ in each $\Z_\ell^*$.  Hence $G = G_{\Q}^{ab}$ is
realized as the Galois group of
$$
\Q_p @>G/D_p>> \wprod_{G/D_p} \Q_p
	@>D_p>> \wprod_{G/D_p} \Q_p(\mu_\infty)
\,,
$$
where the first map is a (pro-)trivial Galois extension.

We can realize the same groups in the $K(1)$-local category, by the two
pro-Galois extensions
$$
J^\wedge_p @>G/D_p>> \wprod_{G/D_p} J^\wedge_p
	@>D_p>> \wprod_{G/D_p} E_1^{nr} \,.
$$
Here the first is the implicitly $K(1)$-localized colimit of
the trivial Galois extensions of $J^\wedge_p$, indexed over the
finite quotients of $G/D_p$.

For brevity, let $B_p = \prod_{G/D_p}' E_1^{nr}$.  Then $J^\wedge_p
\to B_p$ is a $K(1)$-local realization of the maximal abelian extension
of $\Q$.  It seems plausible to find arithmetic pullback squares
$$
\xymatrix{
L_{KU}S \ar[r] \ar[d] & \prod_p J^\wedge_p \ar[d]
& B \ar[r] \ar[d] & \prod_p B_p \ar[d] \\
L_0S \ar[r] & L_0 \prod_p J^\wedge_p
& L_0B \ar[r] & L_0 \prod_p B_p
}
$$
of commutative $S$-algebras, so as to get an integral $KU$-local
realization $L_{KU}S \to B$ of the same Galois group.  It would be
wonderful if analogous (non-abelian) $K(n)$-local constructions for
$n\ge2$ turn out to detect more of the absolute Galois group of $\Q_p$
in $\G_n^{nr}$, or of the absolute Galois group of $\Q$.  The paper
\cite{Mo05} may be relevant.

\subsubhead 5.5.4. $p$-local topological $K$-theory \endsubsubhead

The $p$-local complex $K$-theory spectrum $KU_{(p)}$ is also a commutative
$S$-algebra, and admits an action by the Adams operation $\psi^r$
and its powers through commutative $S$-algebra maps \cite{BR:g, 9.2}.
However, in this case the $E(1)$-local extension
$$
KU_{(p)}^{h\langle r\rangle} \to KU_{(p)}^{h\langle r^n\rangle}
$$
is not a $\Z/n$-Galois extension.  It even fails to be one rationally,
i.e., $K(0)$-locally.  For
$$
\pi_*(KU_{(p)}^{h\langle r^n\rangle}) \otimes \Q \cong E_\Q(\zeta_{r^n})
$$
is an exterior algebra over $\Q$ on one generator, and $E_\Q(\zeta_r)
\to E_\Q(\zeta_{r^n})$ is an isomorphism, so $\{e\}$-Galois, but not
$\Z/n$-Galois.  In spite of the relatively rich source of $K(n)$-local
Galois extensions, there are ramification phenomena that frequently
enter when several chromatic strata are involved.

The idempotent operation $(p-1)^{-1} \sum_{k \in \F_p^*} \psi^k$ on
$KU^\wedge_p$ that defines the $p$-complete Adams summand $L^\wedge_p$
is in fact $p$-locally defined \cite{Ad69, p.~85}, so as to split off
the $p$-local Adams summand $L_{(p)}$ in
$$
KU_{(p)} \simeq \bigvee_{i=0}^{p-2} \Sigma^{2i} L_{(p)} \,.
$$
However, the $p$-adic Adams operations $\psi^k$ of finite order, for
$k$ in the torsion subgroup $\F_p^* \subset \Z_p^*$, are not defined
over $\Z_{(p)}$, since $\psi^k(u) = ku$ on homotopy.  Therefore the
extension $L_{(p)} \to KU_{(p)}$ only becomes Galois after $p$-adic
completion.  This provides an example of an $E(1)$-local {\'e}tale
extension (in the sense of Section~9.4) that does not extend to a
Galois extension.  Again, this is an instance of $K(0)$-local
ramification of the $E(1)$-local prolongation of a, by definition
unramified, $K(1)$-local Galois extension.  These examples are meant as
partial justification for the last paragraph of Remark~3.2.4.

\subhead 5.6. Cochain $S$-algebras \endsubhead

Let $G$ be a topological group and consider a principal $G$-bundle $\pi \:
P \to X$.  Fix a rational prime $p$ and let $A = F(X_+, H\F_p)$ and $B =
F(P_+, H\F_p)$ be the mod~$p$ cochain $H\F_p$-algebras on $X$ and $P$,
respectively.  Note that $\pi_*(A) = H^{-*}(X; \F_p)$ and $\pi_*(B) =
H^{-*}(P; \F_p)$.  We think of $A$ and $B$ as models for the singular
cochain algebras $C^*(X; \F_p)$ and $C^*(P; \F_p)$, in conformance with
\cite{DGI:d, \S3}.  The direct relation between the differential graded
$E_\infty$ structure on $C^*(X; \F_p)$ and the commutative $S$-algebra
structure on $A = F(X_+, H\F_p)$ seems not to have been made explicit,
however.

The projection $\pi$ induces a map of commutative $H\F_p$-algebras $A
\to B$, the right action of $G$ on $P$ induces a left action of $G$ on
$B$ through commutative $A$-algebra maps, and the weak equivalence $P
\times_G EG \to X$ makes its cochain dual $i \: A \to F((P \times_G
EG)_+, H\F_p) \cong B^{hG}$ a weak equivalence.  We now investigate when
$h \: B \wedge_A B \to F(G_+, B)$ is a weak equivalence.

The K{\"u}nneth spectral sequence
$$
E^2_{s,t} = \Tor_{s,t}^{\pi_*(A)}(\pi_*(B), \pi_*(B))
\Longrightarrow \pi_{s+t}(B \wedge_A B)
\tag 5.6.1
$$
can be derived from the skeleton filtration of the (simplicial) two-sided
bar construction
$$
B^{H\F_p}(B, A, B) \: [q] \mapsto B \wedge A^{\wedge q} \wedge B
$$
with all smash products formed over $H\F_p$ \cite{EKMM97, IV.7.7}.
Dually, let
$$
\Omega(P, X, P) \: [q] \mapsto P \times X^q \times P
$$
be the (cosimplicial) two-sided cobar construction, with totalization
equal to the fiber product $P \times_X P$.  There is a natural simplicial
map
$$
\wedge \: B^{H\F_p}(B, A, B) \to F(\Omega(P, X, P), H\F_p)
$$
which is a degreewise weak equivalence by the K{\"u}nneth formula in
mod~$p$ cohomology, under the assumption that $H_*(X; \F_p)$ and $H_*(P;
\F_p)$ are finite in each degree.  So the K{\"u}nneth spectral sequence
equals the one obtained by applying mod~$p$ cohomology to the cobar
construction, i.e., the mod~$p$ Eilenberg--Moore spectral sequence
$$
E^2_{s,t} = \Tor_{s,t}^{H^*(X; \F_p)}(H^*(P; \F_p), H^*(P; \F_p))
\Longrightarrow H^{-(s+t)}(P \times_X P; \F_p)
\tag 5.6.2
$$
\cite{EM66}.  By \cite{Dw74}, \cite{Sh96, 3.1}, the Eilenberg--Moore
spectral sequence converges strongly if, for example, $\pi_0(G)$
is finite, $X$ is path-connected, and $\pi_1(X)$ acts nilpotently on
$H_*(G; \F_p)$.  The K{\"u}nneth spectral sequence is always strongly
convergent, so this comparison implies that the upper horizontal map in
$$
\xymatrix{
B \wedge_A B \ar[r]^-\wedge \ar[d]_h
	& F((P \times_X P)_+, H\F_p) \ar[d]^{\xi^*} \\
F(G_+, B) \ar[r]^-\cong
	& F((P \times G)_+, H\F_p)
}
$$
is a weak equivalence.  The right hand vertical map is induced by
the homeomorphism $\xi \: P \times G \to P \times_X P$, hence is an
isomorphism, as is the lower horizontal map.  Therefore these hypotheses
ensure that the left hand vertical map $h$ is a weak equivalence.

\proclaim{Proposition 5.6.3}
Let $G$ be a stably dualizable group and $P \to X$ a principal $G$-bundle.

(a)
Suppose that $\pi_0(G)$ is finite, $X$ is path-connected, $\pi_1(X)$
acts nilpotently on $H_*(G; \F_p)$, and that $H_*(X; \F_p)$ and $H_*(P;
\F_p)$ are finite in each degree.  Then the map of cochain
$H\F_p$-algebras
$$
F(X_+, H\F_p) \to F(P_+, H\F_p)
$$
is a $G$-Galois extension.

(b)
In particular, when $G$ is a finite discrete group acting nilpotently
on $\F_p[G]$ (this includes all finite $p$-groups), then there is a
$G$-Galois extension
$$
F(BG_+, H\F_p) \to F(EG_+, H\F_p) \simeq H\F_p
$$
that exhibits $H\F_p$ as a Galois extension by each such group.
\endproclaim

A similar argument applies for the map of rational cochain algebras
$$
F(X_+, H\Q) \to F(P_+, H\Q) \,,
$$
when $H^*(X; \Q)$ and $H^*(P; \Q)$ are finite dimensional over $\Q$ in
each degree.

For each natural number $n$ the Morava $K$-theory spectrum $K(n)$
admits uncountably many associative $S$-algebra structures \cite{Rob89,
2.5}, none of which are strictly commutative (cf.~Lemma~5.6.4).  Therefore
$$
F(X_+, K(n)) \to F(P_+, K(n))
$$
is at best a kind of non-commutative $G$-Galois extension.  As a
further complication, the convergence of the $K(n)$-based
Eilenberg--Moore spectral sequence, analogous to~(5.6.2), is not yet
well understood.

\proclaim{Lemma 5.6.4}
$K(n)$ does not admit the structure of a commutative $S$-algebra.
\endproclaim

\demo{Proof}
Suppose that $K(n)$ is a commutative $S$-algebra.  Then so is its
connective cover $k(n)$, and there is a $1$-connected commutative
$S$-algebra map $u \: k(n) \to H\F_p$.  Then $u_* \: H_*(k(n); \F_p)
\to H_*(H\F_p; \F_p)$ is an injective algebra homomorphism, that
commutes with the Dyer--Lashof operations on both sides \cite{BMMS86,
III.2.3}.  The target equals the dual Steenrod algebra $A_* =
E(\chi\tau_k \mid k\ge0) \otimes P(\chi\xi_k \mid k\ge1)$, and the
image of $u_*$ contains $\chi\tau_{n-1}$, but not $\chi\tau_n$.  This
contradicts the operation $Q^{p^k}(\chi\tau_k) = \chi\tau_{k+1}$ in
$A_*$, in the case $k = n{-}1$.
\qed
\enddemo

\head 6. Dualizability and alternate characterizations \endhead

\subhead 6.1. Extended equivalences \endsubhead

Let $A \to B$ be a map of $E$-local commutative $S$-algebras, and let $G$
be a topological group acting from the left on $B$ through $A$-algebra
maps, say by $\alpha \: G_+ \wedge B \to B$.  For example, $A \to B$
could be a $G$-Galois extension.

The twisted group $S$-algebra $B\langle G\rangle$ is defined to be $B \wedge
G_+$ (implicitly $E$-localized, like $B[G]$), with the
multiplication $B\langle G\rangle \wedge B\langle G\rangle \to B\langle
G\rangle$ obtained from the composite map
$$
G_+ \wedge B @>\Delta\wedge1>> G_+ \wedge G_+ \wedge B @>1\wedge\alpha>>
G_+ \wedge B \cong B \wedge G_+
$$
and the multiplications on $B$ and $G$.  As usual, $\Delta$ is the
diagonal map.  The map $A \to B$ and the unit inclusion $\{e\} \to G$
induce a central map $\eta \: A \to B\langle G\rangle$, which makes
$B\langle G\rangle$ an associative $A$-algebra.  Likewise, the
endomorphism algebra $F_A(B, B)$ of $B$ over $A$ is an associative
$A$-algebra with respect to the composition pairing.

Let
$$
j \: B\langle G\rangle \to F_A(B, B)
\tag 6.1.1
$$
be the canonical map of $A$-algebras that is right adjoint to the
composite map
$$
B \wedge G_+ \wedge_A B @>1\wedge\alpha>> B \wedge_A B @>\mu>> B \,,
$$
induced by the ($A$-linear) action of $G$ on $B$ and the multiplication
on $B$.  Note that $B\langle G\rangle$ and $F_A(B, B)$ are left
$B$-modules, with respect to the action on the target in the latter
case, and that $j$ is a map of $B$-modules.  There is also a diagonal
left action by $G$ on $B \wedge G_+$ and on the target in $F_A(B, B)$,
and $j$ is $G$-equivariant with respect to these actions.  These $B$-
and $G$-actions do not commute, but combine to a left module action by
$B\langle G\rangle$.

For a map $f$ of spectra, we will write $f_\#$ and $f^\#$ for various
maps induced by left and right composition with $f$, respectively.

\proclaim{Lemma 6.1.2}
Let $A \to B$ be a map of commutative $S$-algebras, and let $G$ be a
stably dualizable group acting on $B$ through $A$-algebra maps, such that
$h \: B \wedge_A B \to F(G_+, B)$ is a weak equivalence.  For example,
$A \to B$ could be a $G$-Galois extension.  Then:

(a) For each $B$-module $M$ there is a natural weak equivalence
$$
h_M \: M \wedge_A B \to F(G_+, M) \,.
$$

(b) The canonical map
$$
j \: B\langle G\rangle \to F_A(B, B)
$$
is a weak equivalence.

(c) For each $B$-module $M$ there is a natural weak equivalence
$$
j_M \: M \wedge G_+ \to F_A(B, M) \,.
$$
\endproclaim

\demo{Proof}
(a)
By definition, $h_M$ is the composite map
$$
M \wedge_A B \cong M \wedge_B B \wedge_A B
	@>1\wedge h>> M \wedge_B F(G_+, B) @>\nu>> F(G_+, M) \,,
$$
which is a weak equivalence because $h$ is a weak equivalence and $G$
is stably dualizable.

(b)
This is the special case of (c) below when $M = B$.

(c)
By definition, $j_M$ is right adjoint to the composite map $M \wedge G_+
\wedge_A B \to M \wedge_A B \to M$ induced by the group action of $G$
on $B$ and the module action of $B$ on $M$.  We can factor $j_M$ in the
stable homotopy category as the following chain of weak equivalences:
$$
\multline
M \wedge G_+ @>1\wedge\rho>> M \wedge DDG_+ @>\nu>> F(DG_+, M)
\cong F_B(B \wedge DG_+, M) \\
@<\nu^\#<< F_B(F(G_+, B), M) @>h^\#>> F_B(B \wedge_A B, M) \cong
F_A(B, M) \,.
\endmultline
$$ 
Here the map $h^\#$ makes sense because $h$ is a map of $B$-modules, and
similarly for $\nu^\#$.  Algebraically, $m \wedge g$ lifts over $\nu^\#$
to the map $f \mapsto f(g) \cdot m$ in $F_B(F(G_+, B), M)$, which $h^\#$
takes to $j_M(m \wedge g)$.
\qed
\enddemo

\proclaim{Lemma 6.1.3}
Let $A \to B$ be a $G$-Galois extension.  For each $B$-module
$M$ the canonical map
$$
\nu' \: M \wedge_A B^{hG} \to (M \wedge_A B)^{hG}
$$
is a weak equivalence.
\endproclaim

\demo{Proof}
The weak equivalence $M \wedge_A A \cong M \to F(G_+, M)^{hG}$ factors
as the composite
$$
M \wedge_A A @>1\wedge i>\simeq> M \wedge_A B^{hG} @>\nu'>>
(M \wedge_A B)^{hG} @>h_M^{hG}>\simeq> F(G_+, M)^{hG}
$$
where $i$ and $h_M$ are weak equivalences by hypothesis and the
previous lemma, respectively.  The $G$-equivariance of $h_M$, needed to
make sense of $h_M^{hG}$, follows like that of $h$.
\qed
\enddemo

\subhead 6.2. Dualizability \endsubhead

For each $G$-Galois extension $R \to T$ of commutative rings, $T$ is a
finitely generated projective $R$-module.  The following is the
analogous statement for $E$-local commutative $S$-algebras.

\proclaim{Proposition 6.2.1}
Let $A \to B$ be a $G$-Galois extension.  Then $B$ is a dualizable
$A$-module.
\endproclaim

\demo{Proof}
We must show that the canonical map $\nu \: D_A B \wedge_A B \to F_A(B,
B)$ is a weak equivalence.  To keep the different $B$'s
apart, we observe more generally that for each $B$-module $M$ there is
a commutative diagram
$$
\xymatrix{
M \wedge_A F_A(B, A) \ar[d]_{1\wedge i_\#} \ar[r]^{\nu}
	& F_A(B, M \wedge_A A) \ar[d]^{(1\wedge i)_\#} \\
M \wedge_A F_A(B, B^{hG}) \ar[d]_{\cong} \ar[r]^{\nu}
	& F_A(B, M \wedge_A B^{hG}) \ar[d]^{\nu'_\#} \\
M \wedge_A F_A(B, B)^{hG} \ar[r]^{\nu}
	& F_A(B, M \wedge_A B)^{hG} \\
M \wedge_A (B \wedge G_+)^{hG} \ar[u]^{1\wedge j^{hG}} \ar[r]^{\nu}
	& (M \wedge_A B \wedge G_+)^{hG} \ar[u]_{j_{M \wedge_A B}^{hG}} \\
M \wedge_A (B \wedge G_+ \wedge S^{adG})_{hG} \ar[u]^{1\wedge N}
	\ar[r]^{\cong}
	& (M \wedge_A B \wedge G_+ \wedge S^{adG})_{hG} \ar[u]_{N} \\
}
$$
where $\nu'_\#$ is a weak equivalence by Lemma 6.1.3, the maps induced by
$i \: A \to B^{hG}$ are weak equivalences by hypothesis, the maps involving
$j$ are well-defined by the $G$-equivariance of $j$ (and $j_{M \wedge_A
B}$), and are weak equivalences by Lemma~6.1.2, and finally the norm maps $N$
from~(3.6.1) are weak equivalences because the spectra with $G$-action in
question have the form $W \wedge G_+$, with $G$ acting freely on itself
\cite{Rog:s, 5.2.5}.  Thus all maps in this diagram are weak equivalences.

The special case when $M = B$ then verifies that $B$ is dualizable
over $A$.
\qed
\enddemo

In the global case $E = S$ it follows from Propositions~3.3.3 and~6.2.1
that in any $G$-Galois extension $A \to B$, $B$ is a semi-finite
$A$-module, i.e., it is weakly equivalent to a retract of a finite cell
$A$-module.  For example, by Proposition~5.3.1 the complexification map
$KO \to KU$ is a global quadratic extension, and indeed, $KU \simeq KO
\wedge C_{\eta}$ is a finite $2$-cell $KO$-module.  However, in the
localized cases the following counterexample shows that dualizability
is probably the best one can hope for.

\example{Example 6.2.2}
Let $p=2$, recall that $L_{K(1)}S = J^\wedge_2$, and consider
the $K(1)$-local quadratic Galois extension $c \: J^\wedge_2 \to
JU^\wedge_2$ from~5.5.2.  We claim that $JU^\wedge_2$ is not a semi-finite
$J^\wedge_2$-module, even if it is a dualizable $J^\wedge_2$-module,
in the $K(1)$-local category.  There is a diagram of horizontal and
vertical fiber sequences:
$$
\xymatrix{
J^\wedge_2 \ar[r] \ar[d]_c & KO^\wedge_2 \ar[rr]^{\psi^3-1} \ar[d]^c
	&& KO^\wedge_2 \ar[d]^c \\
JU^\wedge_2 \ar[r] \ar[d] & KU^\wedge_2 \ar[rr]^{\psi^3-1} \ar[d]^{\partial}
	&& KU^\wedge_2 \ar[d]^{\partial} \\
\Sigma^2 X_3 \ar[r]
	& \Sigma^2 KO^\wedge_2 \ar[rr]^{\Sigma^2(3^{-1}\psi^3 - 1)}
	&& \Sigma^2 KO^\wedge_2
}
$$
The factor $3^{-1}$ in the lower row comes from the appearance
of the inverse of the Bott equivalence $\beta \: \Sigma^2 KU \to KU$
in the connecting map $\partial$, and the relation $\psi^k \beta = k
\beta \psi^k$.  By definition, following \cite{HMS94, 2.6}, but using
real $K$-theory, $X_3$ is the homotopy fiber of $3^{-1}\psi^3-1 \:
KO^\wedge_2 \to KO^\wedge_2$.

We can compute the zero-th $E_1 = KU^\wedge_2$-cohomology of the spectra
in the upper left hand square, as modules over the group $\SS_1 = \Z_2^*$
of stable Adams operations, with $k \in \Z_2^*$ acting by $\psi^k$.
First, $E_1^0(KU^\wedge_2) \cong \Z_2[[\Z_2^*]]$ (see also Example~8.1.4),
and the remaining modules are the following quotients:
$$
\xymatrix{
\Z_2 & \Z_2[[\Z_2^*/\langle-1\rangle]] \ar[l] \\
\Z_2[[\Z_2^*/\langle3\rangle]] \ar[u]^{c^*}
	& \Z_2[[\Z_2^*]] \ar[l] \ar[u]
}
$$
Here $\langle3\rangle \subset \Z_2^*$ is the subgroup topologically
generated by $3$.  The map $c^*$ takes $E_1^0(JU^\wedge_2) \cong
\Z_2[[\Z_2^*/\langle3\rangle]] \cong \Z_2\{1, \psi^{-1}\}$ to
$E_1^0(J^\wedge_2) \cong \Z_2\{1\}$ by mapping both $1$ and $\psi^{-1}$
to the generator.  Thus $E_1^0(\Sigma^2 X_3) = \ker(c^*) \cong \Z_2\{1
- \psi^{-1}\}$ is such that $\psi^3$ acts as the identity, but
$\psi^{-1}$ acts by reversing the sign.

We claim that there is no semi-finite spectrum with this Morava
module, i.e., this $E_1$-cohomology as an $\SS_1$-module.  For each
finite cell spectrum $X$ the Atiyah--Hirzebruch spectral sequence
$$
E_2^{s,t} = H^s(X; \pi_{-t}(E_1))
\Longrightarrow E_1^{s+t}(X)
$$
is strongly convergent.  After rationalization (inverting~$2$) it
collapses at the $E_2$-term, yielding the Chern character isomorphism
$$
ch \: E_1^0(X)[2^{-1}] \cong \bigoplus_{i \in \Z} H^{2i}(X; \Q_2)
$$
in degree zero.  Here the $i$-th summand appears as the eigenspace of
weight~$i$, where $\psi^k$ acts by multiplication by $k^i$ for each $k \in
\Z_2^*$.  By naturality, there is also such an eigenspace decomposition
of $E_1^0(X)[2^{-1}]$ for each semi-finite $J^\wedge_2$-module~$X$.
(For general spectra $X$, the Atiyah--Hirzebruch spectral sequence needs
not converge.)

Now note that $E_1^0(\Sigma^2 X_3)[2^{-1}] \cong \Q_2$ has $\psi^3$
acting as the identity, and $\psi^{-1}$ acting by sign, which means that
it should lie both in the weight~$0$ eigenspace and in an eigenspace
of odd weight.  This contradicts the possibility that $\Sigma^2 X_3$ is
semi-finite.  It follows that also $JU^\wedge_2$ cannot be $K(1)$-locally
semi-finite.
\endexample

Dualizable modules are preserved under base change, and are detected by
faithful and dualizable base change.

\proclaim{Lemma 6.2.3}
Let $A \to B$ be a map of commutative $S$-algebras and $M$ a dualizable
$A$-module.  Then $B \wedge_A M$ is a dualizable $B$-module.
\endproclaim

\demo{Proof}
We must verify that the canonical map
$$
\nu \: F_B(B \wedge_A M, B) \wedge_B (B \wedge_A M)
\to F_B(B \wedge_A M, B \wedge_A M)
$$
is a weak equivalence.  It factors as the composite
$$
\multline
F_B(B \wedge_A M, B) \wedge_B (B \wedge_A M) \cong F_A(M, B) \wedge_A M \\
@>\nu>> F_A(M, B \wedge_A M) \cong F_B(B \wedge_A M, B \wedge_A M) \,,
\endmultline
$$
where the middle map is a weak equivalence by Lemma~3.3.2(a), since $M$
is a dualizable $A$-module.
\qed
\enddemo

\proclaim{Lemma 6.2.4}
Let $A \to B$ be a faithful map of commutative $S$-algebras, with $B$
dualizable over $A$, and let $M$ be an $A$-module such that $B \wedge_A
M$ is a dualizable $B$-module.  Then $M$ is a dualizable $A$-module.
\endproclaim

\demo{Proof}
We must verify that $\nu \: F_A(M, A) \wedge_A M \to F_A(M, M)$ is
a weak equivalence.  It suffices to show that the map $1 \wedge \nu$
in the commutative square below is a weak equivalence, since $B$ is
assumed to be faithful over $A$.
$$
\xymatrix{
B \wedge_A F_A(M, A) \wedge_A M \ar[r]^-{1\wedge\nu} \ar[d]_{\nu\wedge1}
	& B \wedge_A F_A(M, M) \ar[d]^{\nu} \\
F_A(M, B) \wedge_A M \ar[r]
	& F_A(M, B \wedge_A M)
}
$$
Here the lower horizontal map is isomorphic to
$$
\nu \: F_B(B \wedge_A M, B) \wedge_B (B \wedge_A M) \to
F_B(B \wedge_A M, B \wedge_A M) \,,
$$
which is a weak equivalence because $B \wedge_A M$ is assumed to be
dualizable over~$B$.  The vertical maps are weak equivalences because $B$
is dualizable over $A$, in view of Lemma~3.3.2(a).  Therefore the upper
horizontal map $1\wedge\nu$ is also a weak equivalence.
\qed
\enddemo

\proclaim{Corollary 6.2.5}
If $A$ is a commutative $S$-algebra and $G$ is a stably
dualizable group, so $S[G]$ is dualizable over $S$, then
$A[G]$ is dualizable over $A$.

Conversely, if $A$ is a faithful commutative $S$-algebra, with $A$
dualizable over $S$, and $G$ is a topological group such that $A[G]$
is dualizable over $A$, then $G$ is stably dualizable.
\endproclaim

The following lemma gives the same conclusion as Lemma~6.1.3, but under
different hypotheses, and will be often used.

\proclaim{Lemma 6.2.6}
Let $A \to B$ be a map of commutative $S$-algebras, let $G$ be a
topological group acting on $B$ through $A$-algebra maps, and let $M$
be a dualizable $A$-module.  Then the canonical map
$$
\nu' \: M \wedge_A B^{hG} \to (M \wedge_A B)^{hG}
$$
is a weak equivalence.
\endproclaim

\demo{Proof}
In the commutative diagram
$$
\xymatrix{
M \wedge_A B^{hG} \ar[r]^-{\rho\wedge1} \ar[d]_{\nu'}
	& D_A D_A M \wedge_A B^{hG} \ar[r]^-{\nu} \ar[d]
	& F_A(D_A M, B^{hG}) \ar[d]^{\cong} \\
(M \wedge_A B)^{hG} \ar[r]^-{(\rho\wedge1)^{hG}}
	& (D_A D_A M \wedge_A B)^{hG} \ar[r]^-{\nu^{hG}}
	& F_A(D_A M, B)^{hG}
}
$$
the horizontal maps derived from $\nu$ and $\rho$ are weak equivalences
because $M$ is dualizable over $A$, and the right hand vertical map
is an isomorphism.  Thus the left hand vertical map $\nu'$ is a weak
equivalence.
\qed
\enddemo

\subhead 6.3. Alternate characterizations \endsubhead

The following alternate characterization of Galois extensions corresponds
to the Auslander--Goldman definition.  Compare Proposition~2.3.2.
Implicit cofibrancy and localization at some $S$-module $E$ is to be
understood.

\proclaim{Proposition 6.3.1}
Let $A \to B$ be a map of commutative $S$-algebras, and let $G$ be a
stably dualizable group acting on $B$ through $A$-algebra maps.
Then $A \to B$ is a $G$-Galois extension if and only if
both $i \: A \to B^{hG}$ and $j \: B\langle G\rangle \to F_A(B, B)$
are weak equivalences and $B$ is a dualizable $A$-module.
\endproclaim

\demo{Proof}
Lemma~6.1.2(b) and Proposition~6.2.1 establish one implication.  For the
converse, suppose that $i$ and $j$ are weak equivalences and that $B$ is
dualizable over~$A$.  We must show that $h \: B \wedge_A B \to F(G_+, B)$
is a weak equivalence.  Again, to keep the $B$'s apart we shall observe
that for each $B$-module $M$ the map $h_M$ factors in the stable homotopy
category as the following chain of weak equivalences:
$$
\multline
M \wedge_A B @>1\wedge\rho>> M \wedge_A D_A D_A B
@>\nu>> F_A(D_A B, M) \cong F_B(D_A B \wedge_A B, M) \\
@<\nu^\#<< F_B(F_A(B, B), M)
@>j^\#>> F_B(B\langle G\rangle, M) \cong F(G_+, M) \,.
\endmultline
$$
Algebraically, the forward image of $m \wedge b$ lifts over $\nu^\#$ to
$f \mapsto f(b) \cdot m$, which maps by $j^\#$ to $h_M(m \wedge b) = \{g
\mapsto g(b) \cdot m\}$.  The hypotheses that $B$ is dualizable over $A$
and $j$ is a weak equivalence thus imply that $h_M$ is a weak equivalence.
The special case $M = B$ lets us conclude that $A \to B$ is $G$-Galois.
\qed
\enddemo

In the presence of faithfulness we have a third characterization of
Galois extensions.  See also Propositions~8.2.8 and~12.1.8.

\proclaim{Proposition 6.3.2}
Let $A \to B$ be a map of commutative $S$-algebras, and let $G$
be a stably dualizable group acting on $B$ through $A$-algebra maps.
Then $A \to B$ is a faithful $G$-Galois extension if and only if $h \:
B \wedge_A B \to F(G_+, B)$ is a weak equivalence and $B$ is faithful
and dualizable as an $A$-module.
\endproclaim

\demo{Proof}
Proposition~6.2.1 provides one implication.  For the converse, suppose that
$h$ is a weak equivalence and that $B$ is dualizable and faithful over
$A$.  We must show that $i \: A \to B^{hG}$ is a weak equivalence, and by
faithfulness it suffices to show that $1 \wedge i \: B \cong B \wedge_A
A \to B \wedge_A B^{hG}$ is a weak equivalence.  In the stable homotopy
category we can identify this map with the chain of weak equivalences
$$
B @>\simeq>> F(G_+, B)^{hG} @<h^{hG}<< (B \wedge_A B)^{hG}
@<\nu'<< B \wedge_A B^{hG} \,.
$$
Here $\nu'$ is a weak equivalence by Lemma~6.2.6, because $B$ is
dualizable over $A$.  We are viewing $h$ as a $G$-equivariant map
with respect to the left $G$-actions specified in Section~4.1.
\qed
\enddemo

Here is a characterization of faithfulness in terms of the norm map.

\proclaim{Proposition 6.3.3}
A $G$-Galois extension $A \to B$ is faithful if and only if the norm
map $N \: (B \wedge S^{adG})_{hG} \to B^{hG}$ is a weak equivalence,
or equivalently, if the Tate construction $B^{tG}$ is contractible.
\endproclaim

\demo{Proof}
If the norm map is a weak equivalence, and $Z$ is an $A$-module so
that $Z \wedge_A B \simeq *$, then $Z \simeq Z \wedge_A B^{hG} \simeq Z
\wedge_A (B \wedge S^{adG})_{hG} \cong (Z \wedge_A B \wedge S^{adG})_{hG}
\simeq *$.  Thus $A \to B$ is faithful.

For the converse, consider $B \wedge_A (-)$ applied to the norm map,
appearing as the left hand vertical map in the following commutative
diagram.
$$
\xymatrix{
B \wedge_A (B \wedge S^{adG})_{hG} \ar[r]^{\cong} \ar[d]_{1\wedge N}
	& (B \wedge_A B \wedge S^{adG})_{hG} \ar[r]^-{(h\wedge1)_{hG}} \ar[d]^N
	& (F(G_+, B) \wedge S^{adG})_{hG} \ar[d]^N \\
B \wedge_A B^{hG} \ar[r]^{\nu'}
	& (B \wedge_A B)^{hG} \ar[r]^{h^{hG}}
	& F(G_+, B)^{hG}
}
$$
The map $\nu'$ is a weak equivalence because $B$ is dualizable over $A$,
by Lemma~6.2.6.  The upper and lower right hand horizontal maps are weak
equivalences since $h$ is $G$-equivariant and a weak equivalence.

The right hand vertical map is the norm map for the spectrum with
$G$-action $F(G_+, B)$.  In the source,
$$
(F(G_+, B) \wedge S^{adG})_{hG}
\simeq (B \wedge DG_+ \wedge S^{adG})_{hG}
\simeq (B \wedge S[G])_{hG} \simeq B
$$
by the stable dualizability of $G$ and the Poincar{\'e} duality
equivalence~(3.5.2).  In the target, $F(G_+, B)^{hG} \simeq B$.  A direct
inspection (inducing up from the case $B = S$, where it suffices to
check on $\pi_0$) verifies that these identifications are compatible
under the norm map.  Therefore the right hand vertical map $N$ is a weak
equivalence, and so the norm map for $B$ must be a weak equivalence,
assuming that $B$ is faithful over $A$.

The second equivalence is obvious from the definition of $B^{tG}$ as
the homotopy cofiber of the norm map.
\qed
\enddemo

\proclaim{Corollary 6.3.4}
Any finite $G$-Galois extension $A \to B$ is faithful if the order $|G|$
of $G$ is invertible in $\pi_0(B)$.
\endproclaim

\demo{Proof}
Under these hypotheses $\pi_*(B_{hG}) \cong \pi_*(B)/G$, $\pi_*(B^{hG})
\cong \pi_*(B)^G$ and the composite
$$
\pi_*(B) \to \pi_*(B)/G @>N_*>> \pi_*(B)^G \to \pi_*(B)
$$
is multiplication by $|G|$, so the norm map $N$ must induce an isomorphism
in homotopy.
\qed
\enddemo

The same conclusion, under different hypotheses (allowing ramification)
appears in Lemma~6.4.3.

\subhead 6.4.  The trace map and self-duality \endsubhead

In this section we work principally in the derived category, i.e., in
the stable homotopy category $\Cal D_{A,E}$.

Let $A \to B$ be a map of $E$-local commutative $S$-algebras, and let
$G$ be a stably dualizable group acting on $B$ through $A$-algebra maps.
Suppose that $i \: A \to B^{hG}$ is a weak equivalence.

\definition{Definition 6.4.1}
The {\it trace map} $tr \: B \wedge S^{adG} \to A$ in $\Cal D_{A,E}$
is defined by the natural chain of maps
$$
B \wedge S^{adG} @>in>> (B \wedge S^{adG})_{hG}
@>N>> B^{hG} @<i<\simeq< A \,,
$$
where $in$ denotes the inclusion induced by $G \subset EG$, and the
wrong-way map $i$ is a weak equivalence.
\enddefinition

When $G$ is finite, the dualizing spectrum $S^{adG} = S$ can of course
be ignored.

\proclaim{Lemma 6.4.2}
The trace map $tr \: B \wedge S^{adG} \to A$ equals the composite map
$$
B \wedge S^{adG} = B \wedge S[G]^{hG} @>\nu>\simeq> (B \wedge S[G])^{hG}
@>(\alpha')^{hG}>> B^{hG} @<i<\simeq< A
$$
where $\alpha' \: B \wedge S[G] \to B$ is the right action derived from
$\alpha \: G_+ \wedge B \to B$ by way of the group inverse.
\endproclaim

\demo{Proof}
The canonical map $\nu \: B \wedge S^{adG} \to (B \wedge S[G])^{hG}$
can be identified with the chain of weak equivalences
$$
B \wedge S^{adG} @>\simeq>> F(G_+, B \wedge S^{adG})^{hG}
@<\nu^{hG}<\simeq< (B \wedge DG_+ \wedge S^{adG})^{hG}
@>\simeq>> (B \wedge S[G])^{hG} \,,
$$
using that $G$ is stably dualizable and the (right $G$-equivariant)
Poincar{\'e} duality equivalence~(3.5.2).  In particular, $\nu$ itself
is a weak equivalence.

The claim is then clear from the commutative diagram
$$
\xymatrix{
B \wedge S^{adG} \ar[r]^-{\nu}_-{\simeq} \ar[d]_{in}
& (B \wedge S[G])^{hG} \ar[r]^{=} \ar[d]_{in}
& (B \wedge S[G])^{hG} \ar[d]_{(in)^{hG}} \ar[dr]^{(\alpha')^{hG}} \\
(B \wedge S^{adG})_{hG} \ar[r]^-{\nu_{hG}}_-{\simeq}
& ((B \wedge S[G])^{hG})_{hG} \ar[r]^-{\kappa}
& ((B \wedge S[G])_{hG})^{hG} \ar[r]_-{\simeq}
& B^{hG}
}
$$
where $\kappa$ is the canonical hocolim/holim exchange map and the bottom
row defines the norm map $N$, as in \cite{Rog:s, 5.2.2}.  The right hand
triangle uses that the homotopy orbits $(B \wedge S[G])_{hG}$ are formed
with respect to the diagonal left $G$-action, so the identification with
$B$ extends the right action map $\alpha'$.  Algebraically, $b \wedge g$
in $B \wedge S[G]$ is identified with $g^{-1}b \wedge e$ in the homotopy
orbits, which maps to $\alpha'(b \wedge g) = g^{-1}b$ in $B$.
\qed
\enddemo

\proclaim{Lemma 6.4.3}
When $G$ is finite the composite $B @>tr>> A \to B$ is homotopic to the
sum over all $g \in G$ of the group action maps $g \: B \to B$, and the
composite $A \to B @>tr>> A$ is homotopic to the map multiplying by the
order $|G|$ of $G$.

Thus, if $|G|$ is invertible in $\pi_0(A)$ then $tr$ is a split surjective
map of $A$-modules, up to homotopy, and $B$ is a faithful $A$-module.
In particular, every $G$-Galois extension $A \to B$ with $|G|$ invertible
in $\pi_0(A)$ is faithful.
\endproclaim

\demo{Proof}
When $G$ is finite, the composite $B @>tr>> A \to B$ can be expressed by
continuing the factorization in Lemma~6.4.2 with the map $B^{hG} \to B$
that forgets homotopy invariance, and therefore factors as
$$
B @>1\wedge\Delta>> B \wedge S[G] @>\alpha'>> B \,,
$$
where $\Delta \: S \to S[G] \simeq \prod_G S$ is the diagonal map.
Clearly this is the sum over the elements $g \in G$ of the group action
maps $g \: B \to B$, up to homotopy.

On the other hand, the composite $A \to B @>tr>> A$ is the map of
$G$-homotopy fixed points induced by the same composite displayed above.
Since the action of each group element is homotopic to the identity when
restricted to the homotopy fixed points, their sum equals multiplication
by the group order $|G|$, up to homotopy.
\qed
\enddemo

\example{Example 6.4.4}
In the $\Z/2$-Galois extension $c \: KO \to KU$ the trace map~$tr$ is
homotopic to the realification map $r \: KU \to KO$, as a $KO$-module map,
and therefore also as an $S$-module map.  For $c^\# \: \Cal D_{KO}(KU,
KO) \to \Cal D_{KO}(KO, KO)$ is injective, and both $tr \circ c$ and $r
\circ c$ are homotopic to the multiplication by~$2$ map $KO \to KO$,
by Lemma~6.4.3.

To justify the claim just made, that $c^\#$ is injective, we use the
equivalence $KU \simeq KO \wedge C_\eta$ and adjunction to identify $c^\#$
with $i^\#$ in the exact sequence
$$
\pi_1(KO) @>\eta^\#>> \pi_2(KO) @>j^\#>> [C_\eta, KO] @>i^\#>> \pi_0(KO)
$$
induced by the cofiber sequence $S^0 @>i>> C_\eta @>j>> S^2 @>\eta>> S^1$.
Here $i^\#$ is injective because $\eta^\#$ is well-known to be surjective.

In particular, the trace map $tr = r \: KU \to KO$ is not split
surjective up to homotopy (it is not even surjective on homotopy
groups), so the analog of the algebraic Proposition~2.3.4(b) does not
hold in topology.
\endexample

Recall from Section~3.6 the shearing equivalence $\zeta \: B \wedge
S[G] \to B \wedge S[G]$ that takes the left action on $S[G]$ to the
diagonal left action on $B$ and $S[G]$.

\definition{Definition 6.4.5}
The {\it trace pairing} $B \wedge_A B \wedge S^{adG} \to A$ in $\Cal
D_{A,E}$ is defined as the composite
$$
B \wedge_A B \wedge S^{adG} @>\mu\wedge1>> B \wedge S^{adG} @>tr>> A
\,.
$$
The {\it discriminant map} $\goth d_{B/A} \: B \wedge S^{adG} \to D_AB$
in $\Cal D_{A,E}$ is defined as the composite
$$
\multline
B \wedge S^{adG} = B \wedge S[G]^{hG}
	@>\nu>\simeq> (B \wedge S[G])^{hG}
	@>\zeta^{hG}>\simeq> (B \wedge S[G])^{hG} \\
@>j^{hG}>> F_A(B, B)^{hG} \cong F_A(B, B^{hG})
	@<i_\#<\simeq< F_A(B, A) = D_AB \,.
\endmultline
$$
Here $j$ is $G$-equivariant with respect to the left $G$-action from
Section~6.1.
\enddefinition

We define $\Pic_E = \Pic_E(S)$ in Definition~6.5.1 below to be the
group of weak equivalence classes of $E$-locally smash invertible
spectra.  The dualizing spectrum $S^{adG}$ is one such \cite{Rog:s,
3.3.4}.  By the $\Pic_E$-graded homotopy groups $\pi_*(Y)$ of a
spectrum $Y$ we mean the collection of groups $\pi_X(Y) = [X, Y]$,
where $X$ ranges through $\Pic_E$.  See \cite{HSt99, 14.1}.  This
includes the ordinary stable homotopy groups as the cases $X = S^n$, $n
\in \Z$, as well as the possibly exceptional case $X = S^{adG}$.

\proclaim{Lemma 6.4.6}
The trace pairing $B \wedge_A B \wedge S^{adG} \to A$ is left adjoint
to the discriminant map $\goth d_{B/A} \: B \wedge S^{adG} \to D_AB$.
Thus $\goth d_{B/A}$ is in fact a map in $\Cal D_{B,E}$, and represents
a $\Pic_E$-graded class in $\pi_* D_A(B)$.
\endproclaim

\demo{Proof}
The first claim is a chase of definitions.  The multiplications by $B$
in the two copies of $B$ in the source of the trace pairing get equalized
by $\mu$, so the adjoint (weak) map $\goth d_{B/A}$ commutes with the
obvious $B$-module actions on $B \wedge S^{adG}$ and $D_AB$.
\qed
\enddemo

\proclaim{Proposition 6.4.7}
If $A \to B$ is a $G$-Galois extension, then the discriminant map
$\goth d_{B/A} \: B \wedge S^{adG} \to D_AB$ is a weak equivalence.
In particular, $B$ is self-dual as an $A$-module, up to an invertible
shift by $S^{adG}$.
\endproclaim

\demo{Proof}
When $A \to B$ is $G$-Galois, $j \: B \wedge G_+ \to F_A(B, B)$ is
a weak equivalence by Lemma~6.1.2(b), so the discriminant map is
defined as a composite of weak equivalences.
\qed
\enddemo

In general, we think of the discriminant map $\goth d_{B/A}$ as a measure
of the extent to which $A \to B$ is ramified.  When it is an equivalence,
we think of the trace pairing as a perfect pairing.

\subhead 6.5. Smash invertible modules \endsubhead

The $K(n)$-local Picard group $\Pic_n = \Pic_{K(n)}(S)$ was introduced in
\cite{HMS94}.  Here is a slight generalization.

\definition{Definition 6.5.1}
Let $A$ be a commutative $S$-algebra, and work locally with respect to
the fixed spectrum $E$.  An $A$-module $M$ is {\it smash invertible\/}
if there exists an $A$-module $N$ such that $N \wedge_A M \simeq A$ as
(implicitly $E$-local) $A$-modules.

Let $\Pic_E(A)$ be the class of weak equivalence classes of $E$-locally
smash invertible $A$-modules.  When $\Pic_E(A)$ is a set we call it the
{\it $E$-local Picard group\/} of $A$, with the group structure induced
by the (implicitly $E$-local) smash product of $A$-modules.
\enddefinition

The following proof of the analog of Proposition~2.3.4(c) is close to
one found by Andy Baker and Birgit Richter in the case of a finite
abelian group $G$.

\proclaim{Proposition~6.5.2}
Let $A \to B$ be a faithful abelian $G$-Galois extension, i.e., one
with $G$ an ($E$-locally stably dualizable) abelian group.  Then $B$
is smash invertible as an $A[G]$-module.
\endproclaim

\demo{Proof}
We consider $B$ as an $A[G]$-module by way of the given left $G$-action.
The smash inverse of $B$ over $A[G]$ will be its functional dual
$D_{A[G]}(B) = F_{A[G]}(B, A[G])$ in the category $\Cal M_{A[G],E}$.
There is a natural counit map
$$
\epsilon \: F_{A[G]}(B, A[G]) \wedge_{A[G]} B \to A[G] \,,
$$
that is left adjoint to the identity map on $F_{A[G]}(B, A[G])$ in
the category of $A[G]$-modules.  In symbols, $\epsilon \: f \wedge x
\mapsto f(x)$.  The claim is that $\epsilon$ is a weak equivalence.
By assumption $B$ is faithful over $A$, so it suffices to verify that
$\epsilon$ becomes an equivalence after inducing up along $A \to B$.
We factor the resulting map $1 \wedge \epsilon$ as
$$
\multline
B \wedge_A F_{A[G]}(B, A[G]) \wedge_{A[G]} B @>\nu'>>
F_{A[G]}(B, B[G]) \wedge_{A[G]} B \\
\cong F_{B[G]}(B \wedge_A B, B[G]) \wedge_{B[G]} (B \wedge_A B)
@>\epsilon_1>> B[G]
\,.
\endmultline
$$
Here $\nu'$ is a weak equivalence because $B$ is dualizable over $A$
(cf.~Lemma~6.2.6), the middle isomorphism is a composite of two standard
adjunctions, and $\epsilon_1$ is a counit of the same sort as $\epsilon$,
now in the category of $B[G]$-modules.  We have left to prove that
$\epsilon_1$ is a weak equivalence.

There is a chain of left $B[G]$-module maps
$$
\multline
(B \wedge_A B) \wedge S^{adG} @>h\wedge1>> F(G_+, B) \wedge S^{adG} \\
@<\nu\wedge1<< B \wedge DG_+ \wedge S^{adG}
@>\simeq>> B[G] @>\chi>> B[G] \,,
\endmultline
\tag 6.5.3
$$
each of which is a weak equivalence.  Here $h$ is a weak equivalence
because $A \to B$ is $G$-Galois, $\nu$ is a weak equivalence because $G$
is stably dualizable, and the unnamed weak equivalence is the identity
on $B$ smashed with the Poincar{\'e} duality equivalence from~(3.5.2).
The latter is left $G$-equivariant with respect to the inverse of the
right $G$-action mentioned in Section~3.5, i.e., with respect to the
left action on $DG_+$ given by right multiplication in the source, the
trivial action on $S^{adG}$, and the inverse of the standard right action
on $B[G]$.  The map $\chi$ is induced by the group inverse in~$G$, and
takes the inverse of the standard right action on $B[G]$ to the standard
left action on $B[G]$.

(When $G$ is finite, the chain simplifies to
$$
B \wedge_A B @>h>> F(G_+, B) @<\kappa<< B[G] @>\chi>> B[G] \,,
$$
where $\kappa$ is the usual inclusion and weak equivalence $B[G] \cong
\bigvee_G B \to \prod_G B =F(G_+, B)$.  Again, the right hand $B[G]$
has the standard left $B[G]$-module structure.)

By \cite{Rog:s, 3.3.4, 3.2.3} the dualizing spectrum $S^{adG}$ is smash
invertible (in the $E$-local stable homotopy category), with smash
inverse its functional dual $S^{-adG} = (DG_+)_{hG}$.  It follows that
the counit map $\epsilon_1$ for the $B[G]$-module $B \wedge_A B$ is the
composite of a weak equivalence and the counit map $\epsilon_2$ for $(B
\wedge_A B) \wedge S^{adG}$.  Furthermore, it follows by naturality
with respect to the chain~(6.5.3) of $B[G]$-module weak equivalences
that the counit map $\epsilon_2$ is related by a chain of
weak equivalences to the counit map
$$
\epsilon_3 \: F_{B[G]}(B[G], B[G]) \wedge_{B[G]} B[G] \to B[G] \,,
$$
for $B[G]$ considered as a left $B[G]$-module in the standard way.  The
latter map $\epsilon_3$ is obviously an isomorphism.
\qed
\enddemo

So each (implicitly $E$-local) abelian $G$-Galois extension $A \to B$
exhibits $B$ as a possibly interesting element in the Picard group
$\Pic_E(A[G])$.

The following converse to Proposition~6.5.2 does not require that $G$
is abelian, but for abelian $G$ it follows that the smash invertibility
of $B$ over $A[G]$ is equivalent to $B$ being faithful over $A$.

\proclaim{Lemma 6.5.4}
Let $A \to B$ be a (not necessarily abelian) $G$-Galois extension.
If $B$ is smash invertible as an $A[G]$-module, i.e., if there exists an
$A[G]$-module $C$ and a weak equivalence $B \wedge_{A[G]} C \simeq A[G]$
of $A$-modules, then $B$ is faithful over $A$.
\endproclaim

\demo{Proof}
If $N \wedge_A B \simeq *$ then $N[G] \cong N \wedge_A A[G] \simeq N
\wedge_A B \wedge_{A[G]} C \simeq *$, and $N$ is a retract of $N[G]$,
so $N \simeq *$.
\qed
\enddemo

\head 7. Galois theory I \endhead

We continue to work locally with respect to some $S$-module $E$.

\subhead 7.1. Base change for Galois extensions \endsubhead

Faithful $G$-Galois extensions $A \to C$ are preserved by base change
along arbitrary maps $A \to B$,
$$
\xymatrix{
C \ar[r] & B \wedge_A C \\
A \ar[u] \ar[r] & B \ar[u]
}
$$
and all Galois extensions are preserved by dualizable base change.
Conversely, (faithful) Galois extensions are detected by faithful and
dualizable base change.  We do not know whether these dualizability
hypotheses are necessary.

\proclaim{Lemma 7.1.1}
Let $A \to B$ be a map of commutative $S$-algebras and $A \to C$ a
faithful $G$-Galois extension.  Then $B \to B \wedge_A C$ is a faithful
$G$-Galois extension.
\endproclaim

\demo{Proof}
The action by $G$ on $C$ through $A$-algebra maps extends uniquely to
an action on $B \wedge_A C$ through $B$-algebra maps, taking $g \: C
\to C$ to $1 \wedge g \: B \wedge_A C \to B \wedge_A C$ on the point
set level, for $g \in G$.  The group $G$ remains stably dualizable,
irrespective of whether it is being regarded as acting on $C$ or $B
\wedge_A C$.

We show that $B \to B \wedge_A C$ is a faithful $G$-Galois extension by
appealing to Proposition~6.3.2.  We know that $C$ is a dualizable
$A$-module by Proposition~6.2.1, and it is faithful by hypothesis.
Therefore $B \wedge_A C$ is a dualizable and faithful $B$-module by the
base change lemmas 6.2.3 and~4.3.3.  It remains to verify that the
canonical map $h \: (B \wedge_A C) \wedge_B (B \wedge_A C) \to F(G_+, B
\wedge_A C)$ is a weak equivalence.  It is the lower horizontal map in
the commutative square
$$
\xymatrix{
B \wedge_A C \wedge_A C \ar[r]^-{1\wedge h} \ar[d]_{\cong} &
B \wedge_A F(G_+, C) \ar[d]^{\nu} \\
(B \wedge_A C) \wedge_B (B \wedge_A C) \ar[r]^-{h} &
F(G_+, B \wedge_A C) \,,
}
\tag 7.1.2
$$
where the upper horizontal map $1 \wedge h$ is a weak equivalence
because $A \to C$ is $G$-Galois, and the right hand vertical map $\nu$
is a weak equivalence because $G$ is stably dualizable.  This verifies
the hypotheses of Proposition~6.3.2, so $B \to B \wedge_A C$ is a faithful
$G$-Galois extension.
\qed
\enddemo

\proclaim{Lemma 7.1.3}
Let $A \to B$ be a map of commutative $S$-algebras, with $B$ dualizable
over $A$, and let $A \to C$ be a $G$-Galois extension.  Then $B \to B
\wedge_A C$ is a $G$-Galois extension.
\endproclaim

\demo{Proof}
The group $G$ is stably dualizable, acts on $B \wedge_A C$ through
$B$-algebra maps, and makes the canonical map $h \: (B \wedge_A C)
\wedge_B (B \wedge_A C) \to F(G_+, B \wedge_A C)$ a weak equivalence,
just as in the previous proof.  In order to verify the conditions in
Definition~4.1.3 of a $G$-Galois extension, it remains to show that the
canonical map $i \: B \to (B \wedge_A C)^{hG}$ is a weak equivalence.
But $B \cong B \wedge_A A \simeq B \wedge_A C^{hG}$, so we can identify
$i$ with $\nu' \: B \wedge_A C^{hG} \to (B \wedge_A C)^{hG}$, which is
a weak equivalence by Lemma~6.2.6 because $B$ is dualizable over $A$.
\qed
\enddemo

\proclaim{Lemma 7.1.4}
Let $A \to B$ and $A \to C$ be maps of commutative $S$-algebras, with $B$
a faithful and dualizable $A$-module, and let $G$ be a stably dualizable
group acting on $C$ through $A$-algebra maps.

(a)
If $B \to B \wedge_A C$ is a $G$-Galois extension, then $A \to C$ is a
$G$-Galois extension.

(b)
If $B \to B \wedge_A C$ is a faithful $G$-Galois extension, then $A \to C$ is a
faithful $G$-Galois extension.
\endproclaim

\demo{Proof}
We must verify that the two maps $i \: A \to C^{hG}$ and $h \: C
\wedge_A C \to F(G_+, C)$ are weak equivalences.  For the first map we
factor the weak equivalence $i \: B \to (B \wedge_A C)^{hG}$ for the
$G$-Galois extension $B \cong B \wedge_A A \to B \wedge_A C$ as the
composite
$$
B \wedge_A A @>1\wedge i>> B \wedge_A C^{hG} @>\nu'>>
(B \wedge_A C)^{hG} \,.
$$
Here the right hand map $\nu'$ is a weak equivalence because $B$ is
dualizable over $A$, by Lemma~6.2.6.  Therefore the left hand map $1
\wedge i$ is a weak equivalence, and so $i \: A \to C^{hG}$ is a weak
equivalence because $B$ is faithful over $A$.

For the second map we use the commutative square~(7.1.2) again.  The
right hand vertical map $\nu$ is a weak equivalence because $G$ is
stably dualizable, and the lower horizontal map $h$ is a weak
equivalence because $B \to B \wedge_A C$ is assumed to be $G$-Galois.
So the upper horizontal map $1 \wedge h$ is a weak equivalence, and so
$h \: C \wedge_A C \to F(G_+, C)$ is a weak equivalence because $B$ is
faithful over~$A$.

Finally, if $B \to B \wedge_A C$ is faithful, then we know that $A \to C$ is
faithful by Lemma~4.3.4.
\qed
\enddemo

\subhead 7.2. Fixed $S$-algebras \endsubhead

Let $G$ be a stably dualizable group and let $A \to B$ be a $G$-Galois
extension.  We consider the sub-extensions that occur as the homotopy
fixed points $C = B^{hK}$, for suitable subgroups $K$ of $G$.

\definition{Definition 7.2.1}
Let $K \subset G$ be a topological subgroup.  We say that $K$ is an
{\it allowable subgroup\/} if (a) $K$ is stably dualizable, (b) the
collapse map $c \: G \times_K EK \to G/K$ induces a stable equivalence
$$
S[G \times_K EK] @>\simeq>> S[G/K] \,,
$$
and (c) as a continuous map of spaces, the projection $\pi \: G \to
G/K$ admits a section up to homotopy.

We consider two allowable subgroups $K$ and $K'$ to be {\it
equivalent\/} if $K \subset K'$ and $S[K] \to S[K']$ is a weak
equivalence, or more generally, if $K$ and $K'$ are related by a chain
of such (elementary) equivalences.  We say that $K$ is an {\it
allowable normal subgroup\/} if, furthermore, $K$ is a normal subgroup
of $G$.
\enddefinition

It follows immediately from~(c) above that the orbit space $G/K$ is
stably dualizable, since $S[G/K]$ is a retract up to homotopy of
$S[G]$, and that there is a homotopy equivalence $G \simeq K \times
G/K$ compatible with the obvious projections $\pi$ and $pr_2$ to
$G/K$.  If $K$ is an allowable normal subgroup then $G/K$ is a stably
dualizable group.

\example{Example 7.2.2}
When $G$ is discrete the allowable subgroups of $G$ are just the
subgroups of $G$ in the usual sense, for then $G$ is a disjoint union
of free $K$-orbits, so $c \: G \times_K EK \to G/K$ is already a weak
equivalence, and there is no difficulty in finding a continuous section
to $\pi \: G \to G/K$.
\endexample

For $A \to B$ a $G$-Galois extension and $K \subset G$ an allowable
subgroup, we can form the following maps of commutative $A$-algebras
$$
F(EG_+, B)^G \to F(EG_+, B)^K \to F(EG_+, B) \,.
$$
In view of the natural weak equivalences $A \to F(EG_+, B)^G$ and $F(EG_+,
B) \to B$, we will keep the notation simple by writing the maps above as
$$
A \to B^{hK} \to B \,.
$$
So to be precise, we interpret $B$ as $F(EG_+, B)$, which then admits a
$K$-action through $B^{hK}$-algebra maps.  Likewise, if $K$ is normal in
$G$ then $B^{hK}$ admits a $G/K$-action through $B^{hG}$-algebra maps,
which in turn are $A$-algebra maps.  An implicit cofibrant replacement
is also necessary at this stage.

Here is the forward part of the Galois correspondence for $E$-local
commutative $S$-algebras.

\proclaim{Theorem 7.2.3}
Let $A \to B$ be a faithful $G$-Galois extension and $K \subset G$ any
allowable subgroup.  Then $C = B^{hK} \to B$ is a faithful $K$-Galois
extension.

If furthermore $K \subset G$ is an allowable normal subgroup,
then $A \to C = B^{hK}$ is a faithful $G/K$-Galois extension.
\endproclaim

\demo{Proof}
We shall detect that $C \to B$ (resp.~$A \to C$) is faithfully Galois by
applying Lemma~7.1.4 to the case of faithful and dualizable base change
along $C \to B \wedge_A C$ (resp.~$A \to B$).  Here $B$ is faithful
and dualizable as an $A$-module by hypothesis and Proposition~6.2.1, so $B
\wedge_A C$ is faithful and dualizable as a $C$-module by Lemma~4.3.3
and Lemma~6.2.3.  In the commutative diagram
$$
\xymatrix{
B \ar[r]
	& B \wedge_A B \ar[r]^-{h}_-{\simeq}
	& F(G_+, B)
	& F(G_+, B) \ar[l]^-{=} \\
C \ar[u] \ar[r] 
	& B \wedge_A C \ar[u] \ar[r]^-{h'}_-{\simeq}
	& F(G_+, B)^{hK} \ar[u]
	& F(G/K_+, B) \ar[l]_-{c^\#}^-{\simeq} \ar[u]_{\pi^\#} \\
A \ar[u] \ar[r] 
	& B \ar[u] \ar[r]_-{=}
	& B \ar[u]
	& B \ar[l]^-{=} \ar[u]_{\pi^\#}
}
\tag 7.2.4
$$
the left hand squares are base change pushouts in the category of
commutative $S$-algebras.

The middle horizontal maps are weak equivalences.  For $h$ is a weak
equivalence by the assumption that $A \to B$ is $G$-Galois.  The map $h'
\: B \wedge_A C = B \wedge_A B^{hK} \to F(G_+, B)^{hK}$ factors as a
composite weak equivalence
$$
B \wedge_A B^{hK} @>\nu'>\simeq> (B \wedge_A B)^{hK} @>h^{hK}>\simeq>
F(G_+, B)^{hK}
$$
using that $B$ is dualizable over $A$ (and Lemma~6.2.6) and that $h$
is a weak equivalence.  Here $K$ acts from the left on $B \wedge_A B$
and $F(G_+, B)$ by restriction of the actions by $G$, i.e., on the second
copy of $B$ in $B \wedge_A B$ and by right multiplication in the source
in $F(G_+, B)$, so in particular $h$ is $K$-equivariant.

Likewise, the right hand horizontal maps are weak equivalences.  For
$c^\#$ is the composite map
$$
F(G/K_+, B) @>c^\#>\simeq> F((G_+)_{hK}, B) \cong F(G_+, B)^{hK}
$$
functionally dual to the collapse map $c \: (G_+)_{hK} = (G \times_K
EK)_+ \to G/K_+$, which is a stable equivalence by part~(b) of the
hypothesis that $K$ is allowable.

Therefore, the induced extension $B \wedge_A C \to B \wedge_A B$ is
weakly equivalent to the map $\pi^\# \:  F(G/K_+, B) \to F(G_+, B)$
functionally dual to the projection $\pi \:  G \to G/K$.  By part~(c)
of the hypothesis that $K$ is allowable there is a weak equivalence
$$
F(G_+, B) \simeq F((K \times G/K)_+, B) \cong F(K_+, F(G/K_+, B)) \,,
$$
compatible with the commutative $S$-algebra maps $\pi^\#$ and $pr_2^\#$
from $F(G/K_+, B)$, so that $\pi^\#$ is indeed weakly equivalent to the
trivial $K$-Galois extension (Section~5.1) of $F(G/K_+, B)$.  In
particular, $B \wedge_A C \to B \wedge_A B$ is faithfully $K$-Galois,
and so by the faithful and dualizable detection result Lemma~7.1.4 it
follows that $C \to B$ is faithfully $K$-Galois.

If furthermore $K$ is normal in $G$, then the induced extension $B \to
B \wedge_A C$ is weakly equivalent to the map $\pi^\# \: B \to F(G/K_+,
B)$ functionally dual to the collapse map $\pi \: G/K \to \{e\}$, i.e.,
to the trivial $G/K$-Galois extension of $B$.  So $B \to B \wedge_A C$
is faithfully $G/K$-Galois, and by Lemma~7.1.4 we can conclude that $A
\to C$ is faithfully $G/K$-Galois.
\qed
\enddemo

The following lemma will be applied in Section~9.1, when we discuss
separable extensions.

\proclaim{Lemma 7.2.5}
Let $A \to B$ be a faithful $G$-Galois extension and $K \subset G$ an
allowable subgroup.  Then $C = B^{hK}$ is faithful and dualizable over
$A$, and the canonical map $\kappa \: B^{hK} \wedge_A B^{hK} \to (B
\wedge_A B)^{h(K \times K)}$ is a weak equivalence.
\endproclaim

\demo{Proof}
It is formal that $A \to C$ is faithful when the composite $A \to C \to B$
is faithful.  For if $N \in \Cal M_A$ has $N \wedge_A C \simeq *$ then
$N \wedge_A B \cong N \wedge_A C \wedge_C B \simeq *$, so $N \simeq *$.

The extension $A \to B$ is faithful with $B$ dualizable over $A$ by
Proposition~6.2.1, and $B \wedge_A C \simeq F(G/K_+, B)$ as in~(7.2.4)
is dualizable over $B$, since $S[G/K]$ is assumed to be a retract up to
homotopy of $S[G]$ and therefore is dualizable over~$S$.  Thus $C$ is
dualizable over $A$ by Lemma~6.2.4.

The map $\kappa$ factors as the composite of two weak equivalences
$$
B^{hK} \wedge_A B^{hK} @>\nu'>\simeq> (B \wedge_A B^{hK})^{hK}
	@>(\nu')^{hK}>\simeq> (B \wedge_A B)^{h(K \times K)}
$$
derived from Lemma~6.2.6, where the first uses that $C = B^{hK}$ (on
the right hand side of the smash product) is dualizable over $A$, and
the second uses that $B$ (on the left hand side of the smash product)
is dualizable over $A$.
\qed
\enddemo

\head 8. Pro-Galois extensions and the Amitsur complex \endhead

We continue to let $E$ be a fixed $S$-module and to work entirely in
the $E$-local category.

\subhead 8.1. Pro-Galois extensions \endsubhead

\definition{Definition 8.1.1}
Let $A$ be an $E$-local cofibrant commutative $S$-algebra, and consider
a directed system of $E$-local finite $G_\alpha$-Galois extensions $A
\to B_\alpha$, such that $B_\alpha \to B_\beta$ is a cofibration of
commutative $A$-algebras for each $\alpha \le \beta$.  Suppose further
that each $A \to B_\alpha$ is an $E$-local sub-Galois extension of $A
\to B_\beta$, so/such that there is a preferred surjection $G_\beta \to
G_\alpha$ with kernel $K_{\alpha\beta}$, and a natural weak equivalence
$B_\alpha \simeq B_\beta^{hK_{\alpha\beta}}$.  Let $B = \colim_\alpha
B_\alpha$, where the colimit is formed in $\Cal C_{A,E}$, and let $G =
\lim_\alpha G_\alpha$, with the (profinite) limit topology.  Then, by
definition, $A \to B$ is an {\it $E$-local pro-$G$-Galois extension.}
\enddefinition

More generally, one might consider a directed system of ($E$-local)
Galois extensions with stably dualizable (rather than finite) Galois
groups $G_\alpha$, arranging that each normal subgroup
$K_{\alpha\beta}$ is stably dualizable.  We prefer to wait for some
relevant examples before discussing the analog of the Krull topology on
the resulting limit group $G$, but compatibility with the ``natural
topology'' on $E$-local $\Hom$-sets (see \cite{HPS97, \S4.4} and
\cite{HSt99, \S11}) is certainly desirable.

For each $\alpha$ the weak equivalence $h_\alpha \: B_\alpha \wedge_A
B_\alpha \to F(G_{\alpha+}, B_\alpha)$ extends by Lemma~6.1.2(a) to a weak
equivalence $h_{\alpha,B} \: B \wedge_A B_\alpha \to F(G_{\alpha+}, B)$.
The colimit of these over $\alpha$ is a weak equivalence
$$
h \: B \wedge_A B \to F(\!(G_+, B)\!) \,,
\tag 8.1.2
$$
where by definition $F(\!(G_+, B)\!) = \colim_\alpha F(G_{\alpha+}, B)$
is the ``continuous'' mapping spectrum with respect to the Krull topology,
and $\colim_\alpha B \wedge_A B_\alpha = B \wedge_A \colim_\alpha B_\alpha
= B \wedge_A B$, since pushout with $B$ commutes with colimits in the
category of commutative $A$-algebras.

Likewise, for each $\alpha$ the weak equivalence $j_\alpha \: B_\alpha
\langle G_\alpha \rangle \to F_A(B_\alpha, B_\alpha)$ extends by
Lemma~6.1.2(c) to a weak equivalence $j_{\alpha,B} \: B \langle G_\alpha
\rangle \to F_A(B_\alpha, B)$.  The limit of these over $\alpha$ is a
weak equivalence
$$
j \: B \llangle G \rrangle \to F_A(B, B) \,,
\tag 8.1.3
$$
where by definition $B \llangle G \rrangle = \lim_\alpha B \langle
G_\alpha\rangle$ is the ``completed'' twisted group $A$-algebra, and
$\lim_\alpha F_A(B_\alpha, B) \cong F_A(\colim_\alpha B_\alpha, B) =
F_A(B, B)$.

\example{Example 8.1.4}
In the case of the $K(n)$-local pro-$\G_n$-Galois extension $L_{K(n)}
S \to E_n$, these weak equivalences induce the isomorphism
$$
\Phi \: E_n{}^\vee_*(E_n) \cong \Map(\G_n, \pi_*(E_n))
$$
that is implicit in \cite{Mo85} and explicit in \cite{St00, Thm.~12}
and \cite{Hov04, 4.11}, and the isomorphism
$$
\Psi \: E_{n*} \llangle \G_n \rrangle \cong E_n^*(E_n)
$$
from \cite{St00, p.~1029} and \cite{Hov04, 5.1}.  The appearance of
the continuous mapping space and the completed twisted group ring
corresponds to the spectrum level colimits and limits above, combined
with the $I_n$-adic completion at the level of homotopy groups induced
by the implicit $K(n)$-localization \cite{HSt99, 7.10(e)}.
\endexample

The pro-Galois formalism thus accounts for the first steps in a proof
of Gross--Hopkins duality \cite{HG94}, following \cite{St00}.  The next
step would be to study the $K(n)$-local functional dual of $E_n$ as the
continuous homotopy fixed point spectrum
$$
L_{K(n)} DE_n = F(E_n, L_{K(n)} S) \simeq F(E_n, E_n)^{h\G_n}
	\simeq (E_n \llangle \G_n \rrangle)^{h\G_n} \,,
$$
but here technical issues related to the continuous cohomology of
profinite groups arise, which are equivalent to those handled by
Strickland.

\subhead 8.2. The Amitsur complex \endsubhead

As usual, let $A$ be a cofibrant commutative $S$-algebra and $B$ a
cofibrant commutative $A$-algebra.

\definition{Definition 8.2.1}
The (additive) {\it Amitsur complex\/} \cite{Am59, \S5}, \cite{KO74,
\S II.2} is the cosimplicial commutative $A$-algebra
$$
C^\bullet(B/A) \: [q] \mapsto B \otimes_A [q]
	= B \wedge_A \dots \wedge_A B
$$
($(q+1)$ copies of $B$), coaugmented by $A \to B = C^0(B/A)$.  Here $B
\otimes_A [q]$ refers to the tensored structure in $\Cal C_{A,E}$, and
the cosimplicial structure is derived from the functoriality of this
construction.  In particular, the $i$-th coface map is induced by
smashing with $A \to B$ after the $i$ first copies of $B$, and the
$j$-th codegeneracy map is induced by smashing with $B \wedge_A B \to
B$ after the $j$ first copies of $B$.

Let the {\it completion of $A$ along $B$} be the totalization
$A^\wedge_B = \Tot C^\bullet(B/A)$ of this cosimplicial resolution.
The coaugmentation induces a natural completion map $\eta \: A \to
A^\wedge_B$ of commutative $A$-algebras.
\enddefinition

Gunnar Carlsson has considered this form of completion in his work on
the descent problem for the algebraic $K$-theory of fields \cite{Ca:d,
\S3}, and Bendersky--Thompson have considered an unstable analog in
\cite{BT00}.  It compares perfectly with Bousfield's $B$-nilpotent
completion \cite{Bo79, \S5}, as extended from spectra to the context of
$A$-modules.

\definition{Definition 8.2.2}
The canonical $B$-based Adams resolution of $A$ in $A$-modules is the
diagram below, inductively defined from $D_0 = A$ by letting $D_{s+1}$
be the homotopy fiber of the natural map $D_s = A \wedge_A D_s \to B
\wedge_A D_s$ for all $s\ge0$.
$$
\xymatrix{
A \ar[d] & D_1 \ar[d] \ar[l] & D_2 \ar[d] \ar[l] & \dots \ar[l] \\
B \ar@{-->}[ur] & B \wedge_A D_1 \ar@{-->}[ur] & B \wedge_A D_2 \ar@{-->}[ur]
}
$$
Continuing, $K_s$ is defined to be the homotopy cofiber of the
composite map $D_s \to A$, and the {\it $B$-nilpotent completion\/} of
$A$ in $A$-modules is the homotopy limit $\hat L^A_B A = \holim_s
K_s$.
\enddefinition

\proclaim{Lemma 8.2.3}
The completion $A^\wedge_B$ of $A$ along $B$ is weakly equivalent to
the Bousfield $B$-nilpotent completion $\hat L^A_B A$ of $A$ formed in
$A$-modules.
\endproclaim

\demo{Proof}
One proof uses Bousfield's paper \cite{Bo03} on cosimplicial
resolutions.  The functor $\Gamma(M) = B \wedge_A M$ defines a triple,
or monad, on $\Cal M_A$, and $A \to C^\bullet(B/A)$ is the
corresponding triple resolution of $A$ \cite{Bo03, \S7}.  The
$B$-module spectra define a class $\Cal G$ of injective models in $\Cal
D_A$, whose $\Cal G$-completion is the $B$-nilpotent completion in
$A$-modules, by \cite{Bo79, 5.8} and \cite{Bo03, 5.7}.  It agrees with
the totalization of the triple resolution by \cite{Bo03, 6.5}, which by
definition is the completion of $A$ along $B$, in the sense above.

A more computational proof follows the unstable case of \cite{BK73,
\S3--5}, especially~5.3.  There is an ``iterated boundary isomorphism''
from the $E_1$-term of the Bousfield--Kan spectral sequence associated
to the $\Tot$-tower of the cosimplicial spectrum $C^\bullet(B/A)$, to
the $E_1$-term of the Adams spectral sequence associated to the tower
of derived spectra $\{D_s\}_s$.  The isomorphisms persist, with a shift
in indexing, upon passage to the tower of cofibers $\{K_s\}_s$.  Since
$\Tot_0 C^\bullet(B/A) = B \simeq K_1$, it follows that the homotopy
limits $A^\wedge_B$ and $\hat L^A_B A$ are also weakly equivalent.
\qed
\enddemo

More generally, for each functor $F$ from commutative $A$-algebras to
a category of spaces or spectra, like the units functor $U = GL_1$,
the Amitsur complex $C^\bullet(B/A; F)$ is the cosimplicial object $[q]
\mapsto F(B \otimes_A [q])$.  It is natural to consider the colimit of its
totalization, as $B$ ranges over a class of $A$-algebras.  When $F$ is the
identity functor, this is the completion defined above.  When $A \to B$
is Galois, or ranges through all Galois extensions, we obtain forms of
Amitsur cohomology \cite{Am59} and Galois cohomology \cite{CHR65, \S5}.
Note that if $\Spec B$ is thought of as a covering of $\Spec A$, then
$\Spec (B \wedge_A B)$ consists of the covering of $\Spec A$ by double
intersections, or fiber products, from the first covering, and likewise
for $\Spec C^q(B/A)$ and $(q+1)$-fold intersections.  We are therefore
recovering a form of {\v C}ech cohomology.  In general, the appropriate
context for what classes of extensions $A \to B$ to consider is that
of a {\it Grothendieck model topology\/} on the category of commutative
$A$-algebras, or a {\it model site.}  We simply refer to \cite{TV05}
for a detailed exposition on this matter.

The following is a form of faithfully projective descent.

\proclaim{Lemma 8.2.4}
If $B$ is faithful and dualizable over $A$, then $\eta \: A \to
A^\wedge_B$ is a weak equivalence, i.e., $A$ is complete along $B$.
\endproclaim

\demo{Proof}
It suffices to prove that $1 \wedge \eta \: B \wedge_A A \to B \wedge_A
A^\wedge_B$ is a weak equivalence.  Here $B \wedge_A A^\wedge_B \simeq
F_A(D_A B, \Tot C^\bullet(B/A)) \cong \Tot F_A(D_A B, C^\bullet(B/A))
\simeq \Tot B \wedge_A C^\bullet(B/A)$,
and
$$
B \wedge_A C^\bullet(B/A) \: [q] \to
	B \wedge_A (B \otimes_A [q]) \cong B \otimes_A [q]_+
$$
admits a cosimplicial contraction to $B$, so $1 \wedge \eta$ is indeed
a weak equivalence.
\qed
\enddemo

Let $G$ be a topological group acting from the left on an $S$-module $M$,
and let
$$
EG_\bullet = B(G, G, *) \: [q] \mapsto \Map([q], G) \cong G^{q+1}
$$
be the usual free contractible simplicial left $G$-space.

\definition{Definition 8.2.5}
The (group) {\it cobar complex\/} for $G$ acting on $M$ is the
cosimplicial $S$-module
$$
C^\bullet(G; M) = F(EG_{\bullet+}, M)^G
\: [q] \mapsto F(G^{q+1}_+, M)^G \cong F(G^q_+, M)
\,.
$$
Its totalization is the homotopy fixed point spectrum
$M^{hG} = \Tot C^\bullet(G; M)$.
\enddefinition

Here the standard identification $F(G^{q+1}_+, M)^G \cong F(G^q_+, M)$
takes the left $G$-map $f \: G^{q+1}_+ \to M$ to the map $\phi \: G^q_+
\to M$ that satisfies
$$
\align
f(g_0, \dots, g_q) &= g_0 \cdot \phi([g_0^{-1}g_1 | \dots | g_{q-1}^{-1}g_q]) \\
\phi([h_1 | \dots | h_q]) &= f(e, h_1, \dots, h_1 \dots h_q)
\endalign
$$
(adapted as needed to make sense when the target is a spectrum).

In the presence of a left $G$-action on $B$ through commutative
$A$-algebra maps, these two cosimplicial constructions can be compared.

\definition{Definition 8.2.6}
There is a natural map of cosimplicial commutative $A$-algebras $h^\bullet
\: C^\bullet(B/A) \to C^\bullet(G; B)$ given in codegree~$q$ by the map
$$
h^q \: B \wedge_A \dots \wedge_A B \to F(G^{q+1}_+, B)^G
\cong F(G^q_+, B)
$$
given symbolically by
$$
\align
b_0 \wedge \dots \wedge b_q &\mapsto \bigl( f \: (g_0, \dots, g_q)
\mapsto g_0(b_0) \cdot \dots \cdot g_q(b_q) \bigr) \\
&\cong \bigl( \phi \: [h_1 | \dots | h_q] \mapsto
b_0 \cdot h_1(b_1) \cdot \dots \cdot (h_1\dots h_q)(b_q)
\bigr) \,.
\endalign
$$
On totalizations, $h^\bullet$ induces a natural map of commutative
$A$-algebras $h' \: A^\wedge_B \to B^{hG}$.
\enddefinition

In codegree~$1$, we can recognize $h^1 \: B \wedge_A B \to F(G_+, B)$
as the canonical map $h$ from~(4.1.2).  It is not hard to give a formal
definition of $h^q$ as the right adjoint of a $G$-equivariant map $B
\otimes_A [q] \wedge \Map([q], G)_+ \to B$.

\proclaim{Lemma 8.2.7}
Let $G$ be a stably dualizable group acting on $B$ through $A$-algebra
maps, and suppose that $h \: B \wedge_A B \to F(G_+, B)$ is a weak
equivalence.  Then $h^\bullet$ is a codegreewise weak equivalence that
induces a weak equivalence $h' \: A^\wedge_B \to B^{hG}$.
\endproclaim

\demo{Proof}
In each codegree~$q$, the map $h^q$ factors as a composite of weak
equivalences of the form
$$
\multline
B^{\wedge_A i} \wedge_A F(G^j_+, B)
@>\simeq>>
B^{\wedge_A (i-1)} \wedge_A F(G^j_+, B \wedge_A B) \\
@>\simeq>>
B^{\wedge_A (i-1)} \wedge_A F(G^j_+, F(G_+, B))
\cong
B^{\wedge_A (i-1)} \wedge_A F(G^{(j+1)}_+, B)
\endmultline
$$
with $j = 0, \dots, q-1$ and $i+j=q$.  Here the first map is a weak
equivalence because~$G$, and thus $G^j$, is stably dualizable, and the
second map is a weak equivalence because $h \: B \wedge_A B \to F(G_+,
B)$ is assumed to be one.  The claim follows by induction.
\qed
\enddemo

The following is close to Proposition~6.3.2.  See also Proposition~12.1.8.

\proclaim{Proposition 8.2.8}
Let $G$ be a stably dualizable group acting on $B$ through commutative
$A$-algebra maps, and suppose that $h \: B \wedge_A B \to F(G_+, B)$
is a weak equivalence.  Then $A \to B$ is $G$-Galois if and only if $A$
is complete along $B$
\endproclaim

\demo{Proof}
We have $i = h' \circ \eta$, with $h'$ a weak equivalence, so $i \: A \to
B^{hG}$ is a weak equivalence if and only if $\eta \: A \to A^\wedge_B$
is a weak equivalence.
\qed
\enddemo

\head 9. Separable and {\'e}tale extensions \endhead

We now address structured ring spectrum analogs of the unique lifting
properties in covering spaces, continuing to work implicitly in some
$E$-local category.  Throughout, we let $A$ be a cofibrant commutative
$S$-algebra and $B$ a cofibrant associative or cofibrant commutative
$A$-algebra.  (There appear to be interesting intermediate theories of
$E_n$ $A$-ring spectra for $1 \le n \le \infty$, in the operadic sense,
but we shall focus on the extreme cases of $E_1 = A_\infty$ $A$-ring
spectra, i.e., associative $A$-algebras, and $E_\infty$ $A$-ring
spectra, i.e., commutative $A$-algebras.)

Our main observations are that $G$-Galois extensions $A \to B$ with $G$
discrete are necessarily separable and dualizable, hence symmetrically
{\'e}tale (= thh-{\'e}tale) and {\'e}tale (= taq-{\'e}tale).  In most
cases of current interest, including $E = S$ and $E = K(n)$ for $0 \le
n \le \infty$, a discrete group~$G$ is stably dualizable if and only if
it is finite.

\subhead 9.1. Separable extensions \endsubhead

The algebraic definition \cite{KO74, p.~74} of a separable extension of
commutative rings can be adapted to stable homotopy theory as follows.

\definition{Definition 9.1.1}
We say that $A \to B$ is {\it separable\/} if the $A$-algebra
multiplication map $\mu \: B \wedge_A B^{op} \to B$, considered as
a map in the stable homotopy category $\Cal D_{B \wedge_A B^{op}}$
of $B$-bimodules relative to $A$, admits a section $\sigma \: B \to
B \wedge_A B^{op}$.  Equivalently, there is a map $\sigma \: B' \to B
\wedge_A B^{op}$ of $B$-bimodules relative to $A$, such that the composite
$\mu \sigma \: B' \to B$ is a weak equivalence.
\enddefinition

Here $B^{op}$ is $B$ with the opposite $A$-algebra multiplication
$\mu \gamma \: B \wedge_A B \cong B \wedge_A B \to B$.  It equals $B$
precisely when $B$ is commutative.  Since $B$ will rarely be cofibrant
as a $B$-bimodule relative to $A$, it is only reasonable to ask for
the existence of a bimodule section $\sigma$ in the stable homotopy
category.  The condition for $A \to B$ to be separable only involves the
bimodule structure on $B$, so it is quite accessible to verification
by calculation.  For example, it is equivalent to the condition that
the algebra multiplication $\mu$ induces a surjection
$$
\mu_\# \: \THH_A^0(B, B \wedge_A B^{op}) \to \THH_A^0(B, B)
$$
of zero-th topological Hochschild cohomology groups.  See \cite{La01,
9.3} for a spectral sequence computing the latter in many cases.

\proclaim{Lemma 9.1.2}
Let $A \to B$ be a $G$-Galois extension, with $G$ a discrete group.
Then $A \to B$ is separable.
\endproclaim

\demo{Proof}
Let $d \: G_+ \to \{e\}_+$ be the continuous (Kronecker delta) map given
by $d(e) = e$ (the unit element in $G$) and $d(g) = *$ (the base point)
for $g \ne e$.  Its functional dual
$$
in_e = d^\# \: B \cong F(\{e\}_+, B) \to F(G_+, B)
$$
and the canonical weak equivalence $h$ define the required weak
$B$-bimodule section $\sigma = h^{-1} \circ in_e$ to $\mu$, as a morphism
in the stable homotopy category.
$$
\xymatrix{
B \ar[dr]_{in_e} \ar@{-->}[r]^-{\sigma}
& B \wedge_A B \ar[d]^h_{\simeq} \ar[r]^-{\mu}
& B \\
& \prod_G B \ar[ur]_{pr_e}
}
\tag 9.1.3
$$
\qed
\enddemo

\proclaim{Proposition 9.1.4}
Let $A \to B$ be a faithful $G$-Galois extension, with $G$ a discrete
group and $K \subset G$ any subgroup.  Then $A \to C = B^{hK}$ is
separable.
\endproclaim

\demo{Proof}
By Example~7.2.2, any subgroup $K$ of $G$ is allowable.  We are
therefore in the situation of Lemma~7.2.5.

The map $h \: B \wedge_A B \to \prod_G B$ is $(K \times K)$-equivariant
with respect to the action $(k_1, k_2) \cdot (b_1 \wedge b_2) = k_1(b_1)
\wedge k_2(b_2)$ in the source, and the action that takes a sequence
$\{g \mapsto \phi(g)\}$ to the sequence $\{g \mapsto k_1(\phi(k_1^{-1}
g k_2))\}$ in the target.  There are maps
$$
\prod_K B @>in_K>> \prod_G B @>pr_K>> \prod_K B
$$
functionally dual to a characteristic map $d_K \: G_+ \to K_+$
(taking $G \setminus K$ to the base point) and the inclusion $K_+
\subset G_+$, whose composite is the identity.  We give $\prod_K B$
the $(K \times K)$-action that takes $\{k \mapsto \phi(k)\}$ to $\{k
\mapsto k_1(\phi(k_1^{-1} k k_2))\}$, so that $in_K$ and $pr_K$ are $(K
\times K)$-equivariant.  The weak equivalence $B \to (\prod_K B)^{hK}$
induces a natural weak equivalence $B^{hK} \to (\prod_K B)^{h(K \times
K)}$ that makes the following diagram commute:
$$
\xymatrix{
B^{hK} \ar[d]_{=} \ar@{-->}[r]
	& B^{hK} \wedge_A B^{hK} \ar[d]^{\kappa}_{\simeq} \ar[r]^-{\mu}
	& B^{hK} \ar[d]^{=} \\
B^{hK} \ar[d]_{\simeq} \ar@{-->}[r]
	& (B \wedge_A B)^{h(K \times K)} \ar[d]^{h_\#}_{\simeq} \ar[r]
	& B^{hK} \ar[d]^{\simeq} \\
(\prod_K B)^{h(K \times K)} \ar[r]^{in_{K\#}}
	& (\prod_G B)^{h(K \times K)} \ar[r]^{pr_{K\#}}
	& (\prod_K B)^{h(K \times K)}
}
$$
The vertical map $\kappa$ is a weak equivalence by Lemma~7.2.5, and the
maps $h_\#$ and $pr_{K\#} \circ in_{K\#}$ are obtained from weak
equivalences by passage to $(K \times K)$-homotopy fixed points, so a
little diagram chase shows that $\mu \: B^{hK} \wedge_A B^{hK} \to
B^{hK}$ does indeed admit a weak bimodule section.
\qed
\enddemo

\remark{Remark 9.1.5}
It is easy to see that separable extensions are preserved by base change.
To detect separable extensions by faithful base change will require some
additional hypotheses, as in \cite{KO74, III.2.2}.
\endremark

\subhead 9.2. Symmetrically {\'e}tale extensions \endsubhead

The {\it topological Hochschild homology} $\THH^A(B)$ of $B$ relative to
$A$ is the geometric realization of a simplicial $A$-module
$$
\xymatrix{
B \ar[r]
	& B \wedge_A B \ar[l]<0.5ex> \ar[l]<-0.5ex>
	\ar[r]<0.5ex> \ar[r]<-0.5ex>
	& B \wedge_A B \wedge_A B
	\ar[l]<1.0ex> \ar[l] \ar[l]<-1.0ex>
	\ar[r]<1.0ex> \ar[r] \ar[r]<-1.0ex>
	& \dots
	\ar[l]<1.5ex> \ar[l]<0.5ex> \ar[l]<-0.5ex> \ar[l]<-1.5ex>
}
$$
with the smash product of $(q+1)$ copies of $B$ in degree~$q$.  See
\cite{EKMM97, IX.2}.  Alternatively, $\THH^A(B)$ can be computed in the
stable homotopy category as
$$
\Tor^{B \wedge_A B^{op}}(B, B) = B \wedge^L_{B \wedge_A B^{op}} B \,.
$$
In the case $A = S$, we will often write $\THH(B)$ for $\THH^S(B)$,
which agrees with the topological Hochschild homology introduced by
Marcel B{\"o}kstedt \cite{BHM93}.  The inclusion of $0$-simplices
defines a natural map $\zeta \: B \to \THH^A(B)$.  When $B$ is
commutative, $\THH^A(B)$ can be expressed in terms of the topologically
tensored structure on $\Cal C_A$ as $B \otimes_A S^1$.

It is also possible to define $\THH^A(B)$ for non-commutative $A$, by
analogy with the definition of Hochschild homology over a
non-commutative ground ring \cite{Lo98, 1.2.11}, but we have found no
occasion to make use of this more general definition.

\definition{Definition 9.2.1}
We say that $A \to B$ is {\it formally symmetrically {\'e}tale} (=
{\it formally thh-{\'e}tale}) if the map $\zeta \: B \to \THH^A(B)$ is a
weak equivalence.  If furthermore $B$ is dualizable as an $A$-module,
then we say that $A \to B$ is {\it symmetrically {\'e}tale} (= {\it
thh-{\'e}tale}).
\enddefinition

\remark{Remark 9.2.2}
This definition of an (symmetrically) {\'e}tale map does not quite
conform to the algebraic case, in that it may be too restrictive to ask
that $B$ is dualizable as an $A$-module.  Instead, it is likely to be
more appropriate to only impose the dualizability condition locally
with respect to some Zariski open cover of $\Spec A$.  This may be
taken to mean that for some set of (smashing, Bousfield) localization
functors $\{L_{E_i}\}_i$, such that the collection $\{A \to L_{E_i}A\}_i$
is a faithful cover in the sense of Definition~4.3.1, each localization
$L_{E_i}B$ is dualizable as an $L_{E_i}A$-module.  The author is undecided
about exactly which localization functors to allow.  However, for Galois
extensions the stronger (global) dualizability hypothesis will always
be satisfied, and this may permit us to leave the issue open.
\endremark

\example{Example 9.2.3}
Note that Definition~9.2.1 implicitly takes place in an $E$-local
category.  By McClure--Staffeldt \cite{MS93, 5.1} at odd primes $p$,
and Angeltveit--Rognes \cite{AnR:h, 8.10} at $p=2$, the inclusion
$\zeta \: \ell \to \THH(\ell)$ is a $K(1)$-local equivalence, where
$\ell = BP\langle 1\rangle$ is the $p$-local connective Adams summand
of topological $K$-theory, so $S \to \ell$ is $K(1)$-locally formally
symmetrically {\'e}tale.  It also follows that the localization of this
map, $J^\wedge_p = L_{K(1)}S \to L_{K(1)} \ell = L^\wedge_p$ is
$K(1)$-locally formally symmetrically {\'e}tale.  Here $L^\wedge_p$ is
the $p$-complete periodic Adams summand, as in~5.5.2.

These maps are not $K(1)$-locally symmetrically {\'e}tale, because
$L^\wedge_p$ is not dualizable as a $J^\wedge_p$-module.  More
globally, $S \to L^\wedge_p$ fails to be $E(1)$-locally formally
symmetrically {\'e}tale.  For by \cite{MS93, 8.1}, $\THH(L^\wedge_p)
\simeq L^\wedge_p \vee L_0(\Sigma L^\wedge_p)$, so $\zeta$ has a
rationally non-trivial cofiber.

Similarly, $\zeta \: ku \to \THH(ku)$ is a $K(1)$-homology equivalence
by Christian Ausoni's calculation \cite{Au:t, 6.5} for $p$ odd, and
\cite{AnR, 8.10} again for $p=2$, so the map $S \to ku$ to connective
topological $K$-theory, and its $K(1)$-localization $J^\wedge_p \to
KU^\wedge_p$, are $K(1)$-locally formally symmetrically {\'e}tale.  The
map $L^\wedge_p \to KU^\wedge_p$ is $K(1)$-locally $\F_p^*$-Galois, as
noted in~5.5.2, so by Lemma~9.2.6 below $L^\wedge_p \to KU^\wedge_p$
is $K(1)$-locally symmetrically {\'e}tale.  In other words, $\zeta \:
ku \to \THH^{\ell}(ku)$ and $\zeta \: KU^\wedge_p \to
\THH^{L^\wedge_p}(KU^\wedge_p)$ are $K(1)$-local equivalences.
\endexample

The terminology ``thh-{\'e}tale'' is that of Randy McCarthy and Vahagn
Minasian \cite{MM03, 3.2}, except that for brevity they suppress the
distinction between the formal and non-formal cases.  The author's
lengthier term ``symmetrically {\'e}tale'' was motivated by the following
definitions and result.

\definition{Definition 9.2.4}
Let $M$ be a $B$-bimodule relative to $A$, i.e., a $B \wedge_A
B^{op}$-module.  The space of {\it associative $A$-algebra derivations
of $B$ with values in~$M$} is defined to be the derived mapping space
$$
\ADer_A(B, M) := (\Cal A_A/B)(B, B \vee M)
$$
in the topological model category of associative $A$-algebras over $B$,
where $pr_1 \: B \vee M \to B$ is the square-zero $A$-algebra extension
of $B$ with fiber $M$.  We say that a $B$-bimodule relative to $A$ is
{\it symmetric\/} if it has the form $\mu^! N$ for some $B$-module $N$,
i.e., if the bimodule action is obtained by composing with the $A$-algebra
multiplication map $\mu \: B \wedge_A B^{op} \to B$.
\enddefinition

\proclaim{Proposition 9.2.5}
$A \to B$ is formally symmetrically {\'e}tale if and only if the space
of associative derivations $\ADer_A(B, M)$ is contractible for each
symmetric $B$-bimodule~$M$.
\endproclaim

\demo{Proof}
Let $\Omega_{B/A}$ be a cofibrant replacement of the homotopy fiber
of $\mu \: B \wedge_A B^{op} \to B$ in the category of $B$-bimodules
relative to $A$.  There is a cofiber sequence
$$
B \wedge_{B \wedge_A B^{op}} \Omega_{B/A} \to B @>\zeta>> \THH^A(B)
$$
and for each $B$-module $N$, with associated symmetric $B$-bimodule $M =
\mu^! N$, there is an adjunction equivalence
$$
\Cal M_{B \wedge_A B^{op}}(\Omega_{B/A}, M) \simeq
\Cal M_B(B \wedge_{B \wedge_A B^{op}} \Omega_{B/A}, N) \,.
$$
Furthermore, there is an equivalence (for each $B \wedge_A B^{op}$-module
$M$)
$$
\ADer_A(B, M) = (\Cal A_A/B)(B, B \vee M)
\simeq \Cal M_{B \wedge_A B^{op}}(\Omega_{B/A}, M)
$$
obtained by Lazarev \cite{La01, 2.2}.  So $\zeta$ is an equivalence if
and only if $B \wedge_{B \wedge_A B^{op}} \Omega_{B/A} \simeq *$, which is
equivalent to $\ADer_A(B, M) \simeq \Cal M_B(B \wedge_{B \wedge_A B^{op}}
\Omega_{B/A}, N)$ being contractible for each symmetric $B$-bimodule $M =
\mu^! N$.

In the $E$-local context, this argument shows that $E_*(\zeta)$ is an
isomorphism if and only if $\ADer_A(B, M) \simeq *$ for each $E$-local
symmetric $B$-module $M$.  For $\Cal A_{A,E}/B$ is a full subcategory
of $\Cal A_A/B$, and likewise for the homotopy categories.
\qed
\enddemo

\proclaim{Lemma 9.2.6}
Each separable extension $A \to B$ of commutative $S$-algebras is formally
symmetrically {\'e}tale.  In particular, each $G$-Galois extension $A
\to B$ with $G$ discrete is symmetrically {\'e}tale.
\endproclaim

\demo{Proof}
By assumption there is a bimodule section $\sigma$ so that the composite
$B @>\sigma>> B \wedge_A B^{op} @>\mu>> B$ is homotopic to the identity.
Smashing with $B$ over $B \wedge_A B^{op}$ tells us that the composite
$$
\THH^A(B) @>\sigma\wedge1>> B @>\zeta>> \THH^A(B)
$$
is an equivalence.  Furthermore, there is a retraction $\rho \: \THH^A(B)
\to B$ given in simplicial degree $q$ by the iterated multiplication map
$\mu^{(q)} \: B \wedge_A \dots \wedge_A B \to B$, since we are assuming
that $B$ is commutative.  Therefore $\zeta$ admits a right and a left
inverse, up to homotopy, and is therefore a weak equivalence.

When $A \to B$ is $G$-Galois with $G$ discrete, we showed in
Lemma~9.1.2 that $A \to B$ is separable and in Proposition~6.2.1 that $B$
is a dualizable $A$-module.  The above argument then implies that $A
\to B$ is symmetrically {\'e}tale.
\qed
\enddemo

\subhead 9.3. Smashing maps \endsubhead

Maps $A \to B$ having the corresponding property to the conclusion
of Proposition~9.2.5 for associative derivations into arbitrary (not
necessarily symmetric) $B$-bimodules relative to $A$, also have a familiar
characterization.  This material is not needed for our Galois theory,
but nicely illustrates the relation of smashing localizations (and
Zariski open sub-objects) to {\'e}tale and symmetrically {\'e}tale maps.

\definition{Definition 9.3.1}
We say that $A \to B$ is {\it smashing\/} if the algebra multiplication
map $\mu \: B \wedge_A B^{op} \to B$ is a weak equivalence.
\enddefinition

In view of the following proposition, smashing maps could also be called
formally associatively {\'e}tale extensions.

\proclaim{Proposition 9.3.2}
$A \to B$ is smashing if and only if $\ADer_A(B, M)$ is
contractible for each $B$-bimodule $M$ relative to $A$.
\endproclaim

\demo{Proof}
This is immediate from the equivalence
$$
\ADer_A(B, M) \simeq \Cal M_{B \wedge_A B^{op}}(\Omega_{B/A}, M)
$$
from \cite{La01}, since $A \to B$ is smashing if and only if $\Omega_{B/A}
\simeq *$.
\qed
\enddemo

The terminology is explained by the following result, one part
of which the author learned from Mark Hovey.

\proclaim{Proposition 9.3.3}
$A \to B$ is smashing if and only if $LM = B \wedge_A M$ defines a
smashing Bousfield localization functor on $\Cal M_A$, in which case $B
= LA$.  In particular, $B$ will be a commutative $A$-algebra.
\endproclaim

\demo{Proof}
Let $B^A_*(-)$ be the homotopy functor on $\Cal M_A$ defined by
$B^A_*(M) = \pi_*(B \wedge_A M)$.  The natural map $M \to B \wedge_A M$
is a $B^A_*$-equivalence, since $A \to B$ is smashing, and $B \wedge_A
M$ is $B^A_*$-local by the prototypical ring spectrum argument of Adams
\cite{Ad71}:  if $B \wedge_A Z \simeq *$ then any map $f \: Z \to B
\wedge_A M$ factors as
$$
Z \to B \wedge_A Z @>1\wedge f>> B \wedge_A B \wedge_A M @>\mu\wedge1>>
B \wedge_A M
$$
and is therefore null-homotopic.  So $LM = B \wedge_A M$ defines a
(Bousfield) localization functor $L$ on $\Cal M_A$.

Conversely, a smashing localization functor $L$ on $\Cal M_A$ produces
an associative $A$-algebra $B = LA$, by \cite{EKMM97, VIII.2.1}, such
that $LM \simeq B \wedge_A M$ (since $L$ is assumed to be smashing).
The idempotency of $L$ then ensures that the multiplication map $B
\wedge_A B^{op} \to B$ is a weak equivalence.
\qed
\enddemo

\proclaim{Lemma 9.3.4}
Each smashing map $A \to LA$ is separable, hence formally symmetrically
{\'e}tale.
\endproclaim

\demo{Proof}
If $A \to B = LA$ is smashing, then $\mu \: B \wedge_A B^{op} \to B$
is an equivalence.  It therefore admits a bimodule section $\sigma$
up to homotopy, so $A \to B$ is separable.
\qed
\enddemo

In general, $LA$ is not dualizable as an $A$-module, as easy algebraic
examples illustrate ($\Z \subset \Z_{(p)}$).  Instead, the local
dualizability of Remark~9.2.2 is more appropriate.

\subhead 9.4. {\'E}tale extensions \endsubhead

We keep on working implicitly in an $E$-local category, now with $B$
a cofibrant commutative $A$-algebra.

For a map $A \to B$ of commutative $S$-algebras, the {\it topological
Andr{\'e}--Quillen homology} $\TAQ(B/A)$ is defined in \cite{Bas99, 4.1} as
$$
\TAQ(B/A) := (LQ_B)(RI_B)(B \wedge^L_A B) \,,
$$
i.e., as the $B$-module of (left derived) indecomposables in the
non-unital $B$-algebra given by the (right derived) augmentation ideal in
the augmented $B$-algebra defined by the (left derived) smash product $B
\wedge^L_A B$, augmented over $B$ by the $A$-algebra multiplication $\mu$.

\definition{Definition 9.4.1}
Let $A \to B$ be a map of commutative $S$-algebras.  We say that $A
\to B$ is {\it formally {\'e}tale} (= {\it formally taq-{\'e}tale})
if $\TAQ(B/A)$ is weakly equivalent to $*$.  If furthermore $B$ is
dualizable as an $A$-module, then we say that $A \to B$ is {\it {\'e}tale}
(= {\it taq-{\'e}tale}).
\enddefinition

Like in Remark~9.2.2, the condition that $B$ is dualizable over $A$ is
likely to be stronger than necessary for $B$ to qualify as {\'e}tale
over $A$, and should eventually be replaced with a local condition over
each subobject in an open cover of $A$.  The apologetic discussion from
the associative/symmetric case applies in the same way here.

The terminology is justified by the following definition and result from
\cite{Bas99}.  The vanishing of $\TAQ(B/A)$ gives a unique infinitesimal
lifting property, up to contractible choice, for geometric maps into
the affine covering represented (in the opposite category) by a formally
{\'e}tale map $A \to B$.
$$
\xymatrix{
B \ar[r]^{=} \ar@{-->}[dr] & B \\
A \ar[u] \ar[r] & B \vee M \ar[u]
}
$$
Compare \cite{Mil80, I.3.22}.

\definition{Definition 9.4.2}
Let $A \to B$ be a map of commutative $S$-algebras and let $M$ be
a $B$-module.  The space of {\it commutative $A$-algebra derivations
of $B$ with values in~$M$} is defined to be the derived mapping space
$$
\CDer_A(B, M) := (\Cal C_A/B)(B, B \vee M)
$$
in the topological model category of commutative $A$-algebras over $B$,
where $pr_1 \: B \vee M \to B$ is the square-zero extension of $B$
with fiber $M$.
\enddefinition

\proclaim{Proposition 9.4.3}
A map $A \to B$ of commutative $S$-algebras is formally {\'e}tale if
and only if $\CDer_A(B, M)$ is contractible for each $B$-module $M$.
\endproclaim

\demo{Proof}
There is an equivalence
$$
\CDer_A(B, M) = (\Cal C_A/B)(B, B \vee M) \simeq
\Cal M_B(\TAQ(B/A), M)
$$
for each $B$-module $M$, by \cite{Bas99, 3.2}.  By considering the
universal example $M = \TAQ(B/A)$, we conclude that $\TAQ(B/A) \simeq
*$ if and only if $\CDer_A(B, M) \simeq *$ for each $B$-module $M$.
In the implicitly local context only $E$-local $M$ occur, so we can
conclude that $\TAQ(B/A)$ is $E$-acyclic, i.e., $E$-locally weakly
equivalent to $*$.
\qed
\enddemo

For a finite commutative $R$-algebra $T$, the two conditions $T \cong
HH^R_*(T)$ and $D_*(T/R) = AQ_*(T/R) = 0$ are logically equivalent
\cite{Gro67, 18.3.1(ii)}, where $HH^R_*$ denotes Hochschild homology
and $D_* = AQ_*$ denotes Andr{\'e}--Quillen homology.  In the context
of commutative $S$-algebras this is only true subject to a connectivity
hypothesis \cite{Min03, 2.8}, due to a convergence issue in the analog
of the Quillen spectral sequence from Andr{\'e}--Quillen homology to
Hochschild homology.  However, one implication (from symmetrically
{\'e}tale to {\'e}tale) does not depend on the connectivity hypothesis
stated there.  In other words, if $\zeta \: B \to \THH^A(B)$ is a weak
equivalence, then $\TAQ(B/A) \simeq *$.  We discuss a proof below,
based on \cite{BMa05}.

There is a counterexample to the opposite implication, due to Mike
Mandell, which is discussed in \cite{MM03, 3.5}.  For $n\ge2$ let $X
= K(\Z/p, n)$ be an Eilenberg--Mac\,Lane space and let $B = F(X_+,
H\F_p)$ be its mod~$p$ cochain $H\F_p$-algebra, with $\pi_*(B) =
H^{-*}(K(\Z/p, n); \F_p)$.  Then $H\F_p \to B$ is formally {\'e}tale, but
not symmetrically (=thh-){\'e}tale.  So, any converse statement deducing
that an {\'e}tale map is symmetrically {\'e}tale must contain additional
hypotheses to exclude this example.

\proclaim{Lemma 9.4.4}
Each (formally) symmetrically {\'e}tale extension $A \to B$ of commutative
$S$-algebras is (formally) {\'e}tale.  In particular, each $G$-Galois
extension $A \to B$ with $G$ discrete is {\'e}tale, and each smashing
localization $A \to LA = B$ is formally {\'e}tale.
\endproclaim

\demo{Proof}
Recall that $\THH^A(B) \simeq B \otimes_A S^1$ as commutative $A$-algebras.
Here $\otimes_A$ denotes the tensored structure on $\Cal C_A$
over unbased topological spaces.  To describe the commutative
$B$-algebra structure on $\THH^A(B)$ in similar terms, and to relate it to
the $B$-module $\TAQ(B/A)$, we will need a tensored structure over based
topological spaces.  This makes sense when we replace $\Cal C_A$ by the
pointed category $\Cal C_B/B$ of commutative $B$-algebras augmented over
$B$.  There is then a (reduced) tensor structure $(-) \widetilde\otimes_B
X$ on $\Cal C_B/B$ over based topological spaces $X$, with
$$
(\Cal C_B/B)(C \widetilde\otimes_B X, C')
	\cong \Map_*(X, (\Cal C_B/B)(C, C')) \,,
$$
where $\Map_*$ denotes the base-point preserving mapping space.
It follows that $(C \widetilde\otimes_B X) \widetilde\otimes_B Y \cong
C \widetilde\otimes_B (X \wedge Y)$.  The unbased and based tensored
structures are related by $C \widetilde\otimes_B X \cong B \wedge_C
(C \otimes_B X)$ and $C \otimes_B T \cong C \widetilde\otimes_B (T_+)$,
for unbased spaces~$T$.

There is a pointed model structure on $\Cal C_B/B$, and the associated
Quillen suspension functor $\E$ is given on cofibrant objects by the
reduced tensor $\E(C) = C \widetilde\otimes_B S^1$ with the based circle.
For each $n\ge0$ we can form the $n$-fold iterated suspension
$$
\E^n(C) = C \widetilde\otimes_B S^n
$$
in $\Cal C_B/B$, so that $\E(\E^n(C)) \cong \E^{n+1}(C)$, and these
objects assemble to a sequential suspension spectrum $\E^\infty(C)$,
in this category.  By \cite{BMa05, Thm.~3}, the homotopy category of
such spectra, up to stable equivalence, is equivalent to the homotopy
category $\Cal D_B$ of $B$-modules, up to weak equivalence.

Base change along $A \to B$ takes $B$ to $B \wedge_A B$, which
is a cofibrant commutative $B$-algebra, augmented over $B$ by the
multiplication map $\mu \: B \wedge_A B \to B$.  Hereafter, write
$C = B \wedge_A B$ for brevity.  By \cite{BMa05, Thm.~4}, the cited
equivalence takes $\E^\infty(C)$ to the topological Andr{\'e}--Quillen
homology spectrum $\TAQ(B/A)$.  So $\E^\infty(C)$ is stably trivial
if and only if $\TAQ(B/A) \simeq *$, i.e., if and only if $A \to B$
is formally {\'e}tale.

On the other hand,
$$
\E(C) = C \widetilde\otimes_B S^1 \simeq B
\wedge_C \THH^B(C) \cong \THH^A(B) \,,
$$
now as commutative $B$-algebras.  So $\E(C)$ is weakly trivial, i.e.,
weakly equivalent to the base point $B$ in $\Cal C_B/B$, if and only if
$\zeta \: B \to \THH^A(B)$ is a weak equivalence.

The proof of the lemma is now straightforward.  If $A \to B$ is formally
symmetrically {\'e}tale, then $\E(C)$ is weakly trivial, and therefore
so is each of its suspensions $\E^n(C) = \E^{n-1}(\E(C))$ for $n\ge1$.
Thus the suspension spectrum $\E^\infty(C)$ is stably trivial (in a
very strong sense), and so $\TAQ(B/A)$ is weakly equivalent to the
trivial $B$-module.
\qed
\enddemo

In the notation of the above proof: $C = B \wedge_A B$ is weakly trivial
in $\Cal C_B/B$ if and only if $A \to B$ is smashing, $\E(C) = \THH^A(B)$
is weakly trivial if and only if $A \to B$ is formally symmetrically
{\'e}tale, and $\E^\infty(C)$ is stably trivial if and only if $A \to B$
is formally {\'e}tale.

\subhead 9.5. Henselian maps \endsubhead

By definition, an {\'e}tale map $A \to B$ has the unique lifting
property up to contractible choice for each square-zero extension of
commutative $A$-algebras $B \vee M \to B$, and satisfies a finiteness
condition.  In this chapter we conversely ask which extensions $D \to C$
of commutative $A$-algebras are such that each {\'e}tale map $A \to B$,
with $B$ mapping to $C$, has this homotopy unique lifting property with
respect to $D \to C$.
$$
\xymatrix{
B \ar[r] \ar@{-->}[dr] & C \\
A \ar[u] \ar[r] & D \ar[u]
}
$$
We shall refer to such $D \to C$ as Henselian maps.  Section~9.6 will
exhibit some interesting examples of Henselian maps.

In the opposite category to that of commutative $A$-algebras, of affine
algebro-geometric objects in a homotopy-theoretic sense \cite{TV05,
\S5.1}, we can view the square-zero extensions as infinitesimal
thickenings of a special kind, forming a generating class of acyclic
cofibrations.  The {\'e}tale extensions then correspond to smooth
and unramified covering maps, and constitute a class of fibrations
characterized by their right lifting property with respect to these
generating acyclic cofibrations, together with a finiteness hypothesis.
The Henselian maps, in turn characterized by their left lifting property
with respect to these fibrations, then form a class of thickenings that
contains all composites of the generating acyclic cofibrations of the
theory, i.e., all infinitesimal thickenings, but which also encompasses
many other maps.  By comparison, in the algebraic context Hensel's lemma
applies to a complete local ring mapping to its residue field, but also
to many other cases.

For a fixed commutative $S$-algebra $A$, this discussion could take
place as above in the context of commutative $A$-algebras, with maps
from (taq-){\'e}tale extensions $A \to B$, but also in the alternate
context of associative $A$-algebras, with maps from symmetrically (=
thh-){\'e}tale extensions.  To be concrete we shall focus on the
commutative case, although all of the formal arguments carry over to
the associative category and extensions by symmetric bimodules.

Throughout this section we continue to work $E$-locally, and let $A$ be
a cofibrant commutative $S$-algebra, $B \to C$ a map of commutative
$A$-algebras and $M$ any $C$-module.  We sometimes consider $M$ as a
$B$-module by pull-back along $B \to C$.  We always make the cofibrant
and fibrant replacements required for homotopy invariance, implicitly.

\proclaim{Lemma 9.5.1}
The square-zero extension $B \vee M \to B$ is the pull-back in $\Cal C_A$
of the square-zero extension $C \vee M \to C$ along $B \to C$,
$$
\xymatrix{
B \ar[r]^{=} \ar@{-->}[dr] \ar@{-->}[drr] & B \ar[r] & C \\
A \ar[r] \ar[u] & B \vee M \ar[r] \ar[u] & C \vee M \ar[u]_{pr_1}
}
$$
so there is a weak equivalence
$$
(\Cal C_A/B)(B, B \vee M) \simeq (\Cal C_A/C)(B, C \vee M) \,.
$$
In particular, both of these spaces are contractible whenever $A \to B$
is formally {\'e}tale.
\endproclaim

\demo{Proof}
The pullback along $B \to C$ of a fibrant replacement for $C \vee M \to
C$ is a fibrant replacement for $B \vee M \to B$, and forming mapping
spaces from a cofibrant replacement for $B$ in $\Cal C_A$ has a left
adjoint given by the tensored structure, hence commutes with pullbacks
and other limits.  So the homotopy fiber at the identity of $B$ of
$\Cal C_A(B, B \vee M) \to \Cal C_A(B, B)$ is weakly equivalent to the
homotopy fiber at $B \to C$ of $\Cal C_A(B, C \vee M) \to \Cal C_A(B,
C)$.
\qed
\enddemo

\proclaim{Lemma 9.5.2}
The commutative diagram
$$
\xymatrix{
B \ar[rr] \ar@{-->}[drr] \ar@{-->}[dr] && C \\
A \ar[u] \ar[r] & C \ar[r]_-{in_1} \ar[ur]^{=} & C \vee M \ar[u]_{pr_1}
}
$$
yields a homotopy fiber sequence
$$
(\Cal C_A/C \vee M)(B, C) \to (\Cal C_A/C)(B, C) \to
(\Cal C_A/C)(B, C \vee M)
$$
for which the middle space is contractible.  In particular, all three
spaces are contractible whenever $A \to B$ is formally {\'e}tale.
\endproclaim

\demo{Proof}
After replacing first $pr_1$ and then $in_1$ by fibrations, the mapping
spaces in $\Cal C_A$ from a cofibrant replacement for $B$ to these
fibrations sit in two fibrations $p$ and~$i$, whose composite $p \circ
i$ is also a fibration.  The fibers of the $i$, $p \circ i$ and $p$
above $B \to C$ then form the desired fiber sequence.
\qed
\enddemo

The following definition is the commutative analog of that in
\cite{La01, 3.3}.

\definition{Definition 9.5.3}
A map $\pi \: D \to C$ of commutative $A$-algebras is a {\it singular
extension\/} if there is an $A$-linear derivation of $C$ with values in
$M$, i.e., a commutative $A$-algebra map $d \: C \to C \vee M$ over $C$,
and a homotopy pull-back square
$$
\xymatrix{
C \ar[r]^-{d} & C \vee M \\
D \ar[r] \ar[u]^{\pi} & C \ar[u]_{in_1}
}
$$
of commutative $A$-algebras.
\enddefinition

For example, the square-zero extension $C \vee \Sigma^{-1} M \to C$ is
the singular extension pulled back from the trivial derivation $d =
in_1 \:  C \to C \vee M$.  So the class of singular extensions contains
the class of square-zero extensions.

\proclaim{Lemma 9.5.4}
For each singular extension $\pi \: D \to C$ the commutative
diagram
$$
\xymatrix{
B \ar[r] \ar@{-->}[dr] \ar@{-->}[drr]
	& C \ar[r]^-{d} & C \vee M \\
A \ar[u] \ar[r] & D \ar[r] \ar[u]^{\pi} & C \ar[u]_{in_1}
}
$$
induces a weak equivalence
$$
(\Cal C_A/C)(B, D) \simeq (\Cal C_A/C \vee M)(B, C) \,.
$$
In particular, both of these spaces are contractible whenever
$A \to B$ is formally {\'e}tale.
\endproclaim

\demo{Proof}
The first part of the proof is like that of Lemma~9.5.1.  The second
claim follows by using the definition of formally {\'e}tale maps to
deduce that both mapping spaces are contractible for all formally
{\'e}tale maps $A \to B$, in the special case when $\pi \: D = C \vee
\Sigma^{-1} M \to C$ is the square-zero extension pulled back from the
trivial derivation $d = in_1 \: C \to C \vee M$.  The right hand
mapping space does not depend on the particular singular extension, so
it follows from the first claim applied to a general singular extension
$\pi \: D \to C$ that also the left hand mapping space is contractible
for arbitrary singular extensions $D \to C$ and formally {\'e}tale maps
$A \to B$.
\qed
\enddemo

In view of \cite{Gro67, 18.5.5} or \cite{Mil80, I.4.2(d)} we can make
the following definition.

\definition{Definition 9.5.5}
Let $D \to C$ be a map of commutative $S$-algebras.  We say that $D \to
C$ is {\it Henselian\/} if for each {\'e}tale map $A \to B$, with $B$
and $D$ commutative $A$-algebras over $C$,
$$
\xymatrix{
B \ar[r] \ar@{-->}[dr] & C \\
A \ar[u] \ar[r] & D \ar[u]
}
$$
the derived mapping space
$$
(\Cal C_A/C)(B, D) \simeq *
$$
is contractible, i.e., if $A \to B$ has the unique lifting property up
to contractible choice with respect to $D \to C$.

If $D$ is a commutative $S$-algebra and $C$ is an associative $D$-algebra,
we say that $D \to C$ is {\it symmetrically (= thh-)Henselian\/} if
for each symmetrically (= thh-){\'e}tale map $A \to B$, in a diagram
as above, the associative $A$-algebra mapping space $(\Cal A_A/C)(B,
D)$ is contractible.
\enddefinition

By the following lemma it suffices (in the commutative case) to verify
the homotopy unique lifting property for the {\'e}tale maps $A \to B$
with $A = D$.  For $A \to B$ {\'e}tale implies $D \to B \wedge_A D$
{\'e}tale by the base change formula $\TAQ(B \wedge_A D/D) \simeq
\TAQ(B/A) \wedge_A D$ \cite{Bas99, 4.6} and Lemma~6.2.3.

\proclaim{Lemma 9.5.6}
Let $B \to C$ and $D \to C$ be maps of commutative $A$-algebras,
with pushout $B \wedge_A D \to C$.  The commutative diagram
$$
\xymatrix{
B \ar[r] \ar@{-->}[drr] & B \wedge_A D \ar[r] \ar@{-->}[dr] & C \\
A \ar[r] \ar[u] & D \ar[r]_{=} \ar[u] & D \ar[u]
}
$$
induces a weak equivalence
$$
(\Cal C_A/C)(B, D) \simeq (\Cal C_D/C)(B \wedge_A D, D) \,.
$$
\endproclaim

\demo{Proof}
Dual to the proof of Lemma~9.5.1.
\qed
\enddemo

\proclaim{Proposition 9.5.7}
The class of Henselian maps $D \to C$ contains the square-zero extensions
$C \vee M \to C$ and the singular extensions $\pi \: D \to C$.  It is
closed under weak equivalences, compositions, retracts and
filtered homotopy limits (for diagrams of maps to a fixed $C$).
\endproclaim

\demo{Proof}
The first claims follow from Lemma~9.5.4 and the remark that square-zero
extensions are trivial examples of singular extensions.  The closure
claims are clear, perhaps except for the the last one.  If $\alpha \mapsto
(D_\alpha \to C)$ is a diagram of Henselian maps to $C$, then let $D =
\holim_\alpha D_\alpha$.  For each {\'e}tale map $A \to B$ (mapping to
$D \to C$ as above) there is a weak equivalence
$$
(\Cal C_A/C)(B, D) \simeq \holim_\alpha (\Cal C_A/C)(B, D_\alpha)
\simeq * \,,
$$
since each $D_\alpha \to C$ is Henselian and the limit category is
assumed to be filtering.
\qed
\enddemo

In fact, the Henselian maps that we will encounter in the following
section are sequential homotopy limits of towers of singular extensions,
and thus of a rather special form.  If desired, the reader can view them
as the residue maps of complete local rings, and refer to them as formal
thickenings, rather than as general Henselian maps.

\subhead 9.6. $I$-adic towers \endsubhead

In this section we let $R$ be a commutative $S$-algebra and $R/I$ an
$R$-ring spectrum, i.e., an $R$-module with homotopy unital and
homotopy associative maps $R \to R/I$ and $R/I \wedge_R R/I \to R/I$.
Define the $R$-module $I$ by the cofiber sequence $I \to R \to R/I$,
and let
$$
I^{(s)} = I \wedge_R \dots \wedge_R I
$$
($s$ copies of $I$) be its $s$-fold smash power over $R$, for each
$s\ge1$.  Define the $R$-module $R/I^{(s)}$ by the cofiber sequence
$I^{(s)} \to R \to R/I^{(s)}$.  There is then a tower of $R$-modules
$$
R \to \dots \to R/I^{(s)} \to \dots \to R/I
\tag 9.6.1
$$
that Baker and Lazarev \cite{BL01, \S4} refer to as the {\it external
$I$-adic tower.}

Angeltveit \cite{An:a, 4.2} recently showed that when $R$ is even
graded, i.e., the homotopy ring $\pi_*(R)$ is concentrated in even
degrees, and $R/I$ is a regular quotient, i.e., $\pi_*(I)$ is an ideal
in $\pi_*(R)$ that can be generated by a regular sequence, then each
$R$-ring spectrum multiplication on $R/I$ can be rigidified to an
associative $R$-algebra structure.

Given that $R/I$ is an associative $R$-algebra, Lazarev \cite{La01,
7.1} proved earlier on that the whole $I$-adic tower can be given the
structure of a tower of associative $R$-algebras, and that each cofiber
sequence
$$
I^{(s)}/I^{(s+1)} \to R/I^{(s+1)} @>\pi>> R/I^{(s)}
$$
is a singular extension of associative $R$-algebras.

It remains an open problem to decide when the diagram~(9.6.1) can be
realized as a tower of commutative $R$-algebras, and whether each map
$R/I^{(s+1)} \to R/I^{(s)}$ can be taken to be a singular extension in
the commutative context.  See \cite{La04, 4.5} for a remark on a
similar problem for square-zero extensions.

The homotopy limit
$$
\hat L^R_{R/I} R = \holim_s R/I^{(s)}
$$
of the $I$-adic tower is the Bousfield {\it $R/I$-nilpotent
completion\/} of $R$, formed in the category of $R$-modules (or
$R$-algebras), which we introduced in Definition~8.2.2.  It is in
general not the same as the Bousfield {\it $R/I$-localization\/} of
$R$, formed in the category of $R$-modules, which we denote by
$L^R_{R/I} R$.

However, Baker and Lazarev \cite{BL01, 6.3} use an {\it internal
$I$-adic tower\/} to prove that when $R$ is even graded and $R/I$ is a
homotopy commutative regular quotient $R$-algebra, then the
$R/I$-nilpotent completion has the expected homotopy ring
$$
\pi_* \hat L^R_{R/I} R \cong \pi_*(R)^\wedge_{\pi_*(I)} \,.
$$
If the regular sequence generating $\pi_*(I)$ is finite, then they also
show that the $R/I$-localization and the $R/I$-nilpotent completion of
$R$, both formed in $R$-modules, do in fact agree
$$
L^R_{R/I} R \simeq \hat L^R_{R/I} R \,,
$$
but we shall most be interested in cases when the regular ideal
$\pi_*(I)$ is not finitely generated.

Proposition~9.5.7 therefore has the following consequence, which admits
some fairly obvious algebraically localized generalizations that we
shall also make use of.

\proclaim{Proposition 9.6.2 (Baker--Lazarev)}
Let $R$ be an even graded commutative $S$-algebra, and $R/I$ a homotopy
commutative regular quotient $R$-algebra.  Then the limiting map
$$
\hat L^R_{R/I} R = \holim_s R/I^{(s)} \to R/I
$$
is symmetrically (=thh-)Henselian, and induces the canonical surjection
$$
\pi_*(R)^\wedge_{\pi_*(I)} \to \pi_*(R)/\pi_*(I)
$$
of homotopy rings.  In particular, if $\pi_*(R)$ is already
$\pi_*(I)$-adically complete, so that $R \simeq \hat L^R_{R/I} R$, then
$R \to R/I$ is symmetrically Henselian.
\qed
\endproclaim

We now claim that the complex cobordism spectrum $MU$ can be viewed
as a global model, up to Henselian maps, of each of the commutative
$S$-algebras $\widehat{E(n)} = L_{K(n)} E(n)$ that occur as fixed
$S$-algebras in the $p$-primary $K(n)$-local pro-Galois extensions
$L_{K(n)} S \to E_n \to E_n^{nr}$.  So, even if there is ramification
between the expected maximal unramified Galois extensions (covering
spaces) over the different chromatic strata, reflected in the changing
pro-Galois groups $\G_n$ and $\G_n^{nr}$ for varying $n$ and~$p$,
these can all be compensated for by appropriate Henselian maps (formal
thickenings), and unified into one global model, namely $MU$.

For the sphere spectrum $S$, the chromatic stratification we have in
mind is first branched over the rational primes $p$, and then $S_{(p)}$
is filtered by the Bousfield localizations $L_n S = L_{E(n)} S$ for
each $n\ge0$.  The associated (Zariski) stack has the category $\Cal
M_{S, E(n)}$ of $E(n)$-local $S$-modules over the $n$-th open subobject
in the filtration, and the $n$-th {\it monochromatic category\/} of
$E(n)$-local $E(n{-}1)$-acyclic $S$-modules over the $n$-th half-open
stratum.  The latter category is equivalent to the category $\Cal M_{S,
K(n)}$ of $K(n)$-local $S$-modules, at least in the sense that their
homotopy categories are equivalent \cite{HSt99, 6.19}.

The latter $K(n)$-local module category is in turn equivalent to
the category of $K(n)$-local $L_{K(n)} S$-modules, and we propose
to understand it better by way of Galois descent from the related
categories $\Cal M_{B, K(n)}$ of $K(n)$-local $B$-modules, for the various
$K(n)$-local Galois extensions $L_{K(n)}S \to B$.  The limiting case
of pro-Galois descent from $K(n)$-local modules over $B = E_n^{nr}$,
or over the separable closure $B = \bar E_n$ (cf.~Section~10.3), can
optimistically be hoped to be particularly transparent.

This decomposition of the sphere spectrum, appearing in the lower row
in the diagram below, can be paralleled for $MU$ by applying the same
localization functors in spectra.  However, the proposition above
indicates that it may be more appropriate to nilpotently complete $MU$,
in the category of $MU$-modules.  In other words, we are led to focus
attention on the upper row, rather than the middle row, in the
following commutative diagram.
$$
\xymatrix{
MU \ar[r] & \hat L^{MU}_{E(n)} MU \ar[r] & \hat L^{MU}_{K(n)} MU \\
MU \ar[r] \ar[u]^{=} & L_{E(n)} MU \ar[r] \ar[u] & L_{K(n)} MU \ar[u] \\
S \ar[r] \ar[u] & L_{E(n)} S \ar[r] \ar[u] & L_{K(n)} S \ar[u]
}
\tag 9.6.3
$$
In the middle column we have $L_{E(n)} S \simeq \hat L_{E(n)} S$, since
every $E(n)$-local spectrum is $E(n)$-nilpotent \cite{HSa99, 5.3}.
However, in the right hand column $L_{K(n)} S \not\simeq \hat L_{K(n)}
S$, since $K(n)$-localization is not smashing \cite{Ra92, 8.2.4} and
\cite{HSt99, 8.1}.

The coefficient rings of the various localizations and nilpotent
completions of $MU$ occurring in the diagram above are mostly
understood.  See \cite{Ra92, 8.1.1} for $\pi_* L_{E(n)} MU$ (or rather,
its $BP$-version).  Let $J_n \subset \pi_* MU_{(p)}$ be the kernel of
the ring homomorphism $\pi_* MU_{(p)} \to \pi_* E(n)$, i.e., the
regular ideal generated by the kernel of $\pi_* MU_{(p)} \to \pi_* BP$
and the infinitely many classes $v_k$ for $k > n$.  Let $I_n = (p, v_1,
\dots, v_{n-1})$, also considered as an ideal in $\pi_* MU_{(p)}$, so
that the sum of ideals $I_n + J_n$ is the kernel of the ring
homomorphism $\pi_* MU_{(p)} \to \pi_* K(n)$.  Then
$$
\pi_* L_{K(n)} MU = \pi_* MU_{(p)} [v_n^{-1}]^\wedge_{I_n}
$$
by \cite{HSt99, 7.10(e)}.  By \cite{HSa99, Thm.~B}, $L_{K(n)} BP$
splits as the $K(n)$-localization of an explicit countable wedge sum of
suspensions of $\widehat{E(n)}$.  It follows that $L_{K(n)} MU$ splits
in a similar way.

By Proposition~9.6.2, applied to $R = MU_{(p)}[v_n^{-1}]$ and $R/I =
E(n)$, we find that $\hat L^{MU}_{E(n)} MU \simeq \hat L^R_{R/I} R \to
E(n)$ is symmetrically Henselian, with
$$
\pi_* \hat L^{MU}_{E(n)} MU = \pi_* MU_{(p)}[v_n^{-1}]^\wedge_{J_n} \,.
$$
By the same proposition applied to $R = MU_{(p)}[v_n^{-1}]$ and $R/I =
K(n)$, at least for $p\ne2$ to ensure that $K(n)$ is homotopy
commutative, we also find that $\hat L^{MU}_{K(n)} MU \simeq \hat
L^R_{R/I} R \to K(n)$ is symmetrically Henselian, with
$$
\pi_* \hat L^{MU}_{K(n)} MU = \pi_* MU_{(p)}[v_n^{-1}]^\wedge_{I_n+J_n}
\,.
\tag 9.6.4
$$
This differs from the $K(n)$-localization of $MU$ in $S$-modules by the
additional completion along $J_n$.

This $K(n)$-nilpotently complete part, in $MU$-modules, of the global
commutative $S$-algebra $MU$, can now be related by a symmetrically
Henselian map to the extension $L_{K(n)} S \to \widehat{E(n)}$, which
is closely related to the $K(n)$-locally pro-Galois extension
$L_{K(n)}S \to E_n$.
$$
\xymatrix{
E_n \ar[r]^-{\simeq?} & \hat L^{MU}_{K(n)} E_n \\
{}\widehat{E(n)} \ar[r]^-{\simeq?} \ar[u] & \hat L^{MU}_{K(n)} E(n) \ar[u] \\
L_{K(n)} MU \ar[r] \ar[u] & \hat L^{MU}_{K(n)} MU \ar[u] \ar@{-->}[ul]_{q} \\
L_{K(n)} S \ar[r]_{=} \ar[u] & L_{K(n)} S \ar[u]
}
\tag 9.6.5
$$
Here the horizontal map $\widehat{E(n)} = L_{K(n)} E(n) \to \hat
L^{MU}_{K(n)} E(n)$, and its analog for $E_n$, are both plausibly weak
equivalences.  For instance, the corresponding map of nilpotent
completions of $MU$ induces completion along $J_n$ at the level of
homotopy groups, and $\pi_* \widehat{E(n)}$ and $\pi_* E_n$ are already
$J_n$-adically complete in a trivial way.

We shall now apply Proposition~9.6.2 with $R = \hat L^{MU}_{K(n)} MU$.
Formula~(9.6.4) exhibits $R$ as an even graded commutative
$S$-algebra.  Considering $J_n$ as an ideal in $\pi_* R$, it is still
generated by a regular sequence and $(\pi_* R)/J_n \cong \pi_*
\widehat{E(n)}$.  So we can form $R/I \simeq \widehat{E(n)}$ as a
homotopy commutative regular quotient $R$-algebra.  Then $\pi_*(I) =
J_n$, and $\pi_*(R)$ is $J_n$-adically complete, so by the last clause
of Proposition~9.6.2 the map $q \: R \to R/I$, labeled $q$ in the
diagram~(9.6.5) above, is symmetrically Henselian.

\proclaim{Corollary 9.6.6}
Each $K(n)$-local pro-Galois extension $L_{K(n)}S \to E_n$
factors as the composite map of commutative $S$-algebras
$$
L_{K(n)}S \to \hat L^{MU}_{K(n)} MU @>q>> \widehat{E(n)} \to E_n \,,
$$
where the first map admits the global model $S \to MU$, the second map
is symmetrically (= thh-)Henselian, and the third map is a $K(n)$-local
pro-Galois extension.

In other words, each $K(n)$-local stratum of $S$ is related by a chain
of pro-Galois covers $L_{K(n)} S \to E_n \leftarrow \widehat{E(n)}$ to
a formal thickening $q \: \hat L^{MU}_{K(n)} MU \to \widehat{E(n)}$ of
the corresponding $K(n)$-nilpotently complete stratum of $MU$, formed
in $MU$-modules.
\endproclaim

We shall argue in Section~12.2 that there is a Hopf--Galois structure
on this global model $S \to MU$ that also encapsulates all the known
Galois symmetries over $L_{K(n)} S$, at least up to the adjunction of
roots of unity, i.e., up to the passage from $\widehat{E(n)}$ to $E_n$
(or to $E_n^{nr}$).  The question remains whether $q$ is (commutatively)
Henselian, which would follow if the diagram~(9.6.1) could be realized
by singular extensions of commutative $S$-algebras.

\medskip

After this discussion of $K(n)$-localization and $K(n)$-nilpotent
completion in $MU$-modules, we make some remarks on the chromatic
filtration in $MU$-modules.  The study of the chromatic filtration and
the monochromatic category of $S$-modules relies on the basic fact
\cite{JY80, 0.1} that $E(n)_*(X) = 0$ implies $E(n{-}1)_*(X)$ for
$S$-modules $X$, so that there is a natural map $L_{E(n)} X \to
L_{E(n-1)} X$.  The analogous claim in the context of $MU$-modules is
false, i.e., that $E(n)^{MU}_*(X) = 0$ implies $E(n{-}1)^{MU}_*(X) =
0$, as the easy example $X = MU_{(p)}/(v_n)$ illustrates.  Thus there
is no natural map $L^{MU}_{E(n)} X \to L^{MU}_{E(n-1)} X$.

For brevity, let $K[0,n] = K(0) \vee \dots \vee K(n)$.  It is
well-known that $L_{K[0,n]} = L_{E(n)}$ in the category of $S$-modules
\cite{Ra84, 2.1(d)}.  For any $MU$-module $X$ it is obvious that
$K[0,n]^{MU}_*(X) = 0$ implies $K[0,n{-}1]^{MU}_*(X) = 0$, so that there
is a natural map $L^{MU}_{K[0,n]} X \to L^{MU}_{K[0,n-1]} X$.
Therefore the example above shows that the two localization functors
$L^{MU}_{K[0,n]}$ and $L^{MU}_{E(n)}$ in $MU$-modules cannot be
equivalent.

We therefore think that it will be more appropriate to filter the
category of $MU$-modules by the essential images
$$
\Cal M_{MU} \supset \dots \supset \Cal M^{MU}_{MU, K[0,n]}
\supset \Cal M^{MU}_{MU, K[0,n-1]} \supset \dots
$$
of the Bousfield localization functors $L^{MU}_{K[0,n]}$, i.e., the
full subcategories of $K[0,n]$-local $MU$-modules, within $MU$-modules,
or the corresponding essential images
$$
\Cal M_{MU} \supset \dots \supset \hat\Cal M^{MU}_{MU, K[0,n]}
\supset \hat\Cal M^{MU}_{MU, K[0,n-1]} \supset \dots
$$
of the nilpotent completion functors $\hat L^{MU}_{K[0,n]}$, i.e., the
full subcategories of $K[0,n]$-nilpotently complete $MU$-modules, within
$MU$-modules.  Then we can consider the $MU$-chromatic towers
$$
X \to \dots \to L^{MU}_{K[0,n]} X \to L^{MU}_{K[0,n-1]} X \to \dots
$$
and
$$
X \to \dots \to \hat L^{MU}_{K[0,n]} X \to \hat L^{MU}_{K[0,n-1]} X \to \dots
$$
for each $MU$-module $X$.  We then suspect that $L^{MU}_{K[0,n]}$ is a
smashing localization, and that there is an equivalence of homotopy
categories between the $n$-th monochromatic category of $MU$-modules
and the $K(n)$-local category of $MU$-modules, like that of
\cite{HSt99, 6.19}, but we have not verified this expectation.  To be
precise, the monochromatic category in question has objects the
$MU$-modules that are $L^{MU}_{K[0,n]}$-local and
$L^{MU}_{K[0,n-1]}$-acyclic.  The $K(n)$-local category has objects the
$MU$-modules that are $L^{MU}_{K(n)}$-local.

The thrust of Corollary~9.6.6 is now that the chromatic filtration on
$S$-modules is related to a chromatic filtration on $MU$-modules, by a
chain of pro-Galois extensions and Henselian maps with geometric content.
The chromatic filtration on $MU$-modules is likely to be much easier to
understand algebraically, in terms of the theory of formal group laws.
Taken together, these two points of view may clarify the chromatic
filtration on $S$-modules.

\head 10. Mapping spaces of commutative $S$-algebras \endhead

We turn to the computation of the mapping space $\Cal C_A(B, B)$
for a $G$-Galois extension $A \to B$, and related mapping spaces of
commutative $S$-algebras, using the Hopkins--Miller obstruction theory
in the commutative form presented by Goerss and Hopkins \cite{GH04}.
For the more restricted problem of the classification of commutative
$S$-algebra structures, the related obstruction theory of Alan Robinson
\cite{Rob03} is also relevant.

\subhead 10.1. Obstruction theory \endsubhead

Let $A$ be a cofibrant commutative $S$-algebra and let $E$ be an
$S$-module.  We shall need an extension of the Goerss--Hopkins theory
to the context of (simplicial algebras over simplicial operads in)
the category $\Cal M_{A,E}$ of $E$-local $A$-modules.  The base change
to $A$-modules is harmless, but in working $E$-locally we may loose
the identification of the dualizable $A$-modules with the (homotopy
retracts of) finite cell $A$-modules, recalled in Proposition~3.3.3 above.
It seems clear that only the formal properties of dualizable modules are
important to the Goerss--Hopkins theory, so that the whole extension
can be carried through in full generality.  However, for our specific
purposes the only dualizable $A$-modules we must consider will in fact
be finite cell $A$-modules, so we do not actually need to carry the
generalization through.

Next, consider a fixed (cofibrant, $E$-local) commutative $A$-algebra
$B$.  The Goerss--Hopkins spectral sequence \cite{GH04, Thm.~4.3 and
Thm.~4.5} for the computation of the homotopy groups of commutative
$A$-algebra mapping spaces like $\Cal C_A(C, B)$, for various commutative
$A$-algebras $C$, is based on working with a fixed homology theory given
by a commutative $A$-algebra that they call $E$, but which we will
take to be $B$.  In particular, the target $B$ in the mapping space
is then equivalent to its completion along the given homology theory
(cf.~Definition~8.2.1), as required for the convergence of the spectral
sequence.

This commutative $A$-algebra $B$ is required to satisfy the so-called
{\it Adams conditions} \cite{Ad69, p.~28}, \cite{GH04, Def.~3.1},
which in our notation asks that $B$ is weakly equivalent to a homotopy
colimit of finite cell $A$-module spectra $B_\alpha$, satisfying two
conditions.  For our purposes it will suffice that $B$ itself satisfies
the two conditions, i.e., that there is only a trivial colimit system.
The conditions are then:

\definition{Adams conditions 10.1.1}
The commutative $A$-algebra $B$ is weakly equivalent to a finite cell
$A$-module, such that

(a) $B^A_*(D_AB)$ is finitely generated and projective as a $B_*$-module.

(b) For each $B$-module $M$ the K{\"u}nneth map
$$
[D_AB, M]^A_* \to \Hom_{B_*}(B^A_*(D_AB), M_*)_*
$$
is an isomorphism.
\enddefinition

In the $E$-local situation we expect that it suffices to assume that $B$
is a dualizable $A$-module, but in our applications the stronger finite
cell hypothesis will always be satisfied.

\proclaim{Lemma 10.1.2}
The Adams conditions~(a) and~(b) are satisfied when $A \to B$ is an
$E$-local $G$-Galois extension, with $G$ a finite discrete group.
\endproclaim

\demo{Proof}
From Lemma~6.1.2 we know that $j \: B\langle G\rangle \to F_A(B, B)$
is a weak equivalence, and that $h_M \: B \wedge_A M \to F(G_+, M)$
is a weak equivalence for each $B$-module~$M$.  By Proposition~6.2.1,
$B$ is dualizable over $A$, so $B \wedge_A D_A B \simeq F_A(B, B)$.
So $B^A_*(D_AB) \cong \pi_* F_A(B, B) \cong B_*\langle G\rangle$ is a
finitely generated free $B_*$-module, and $B^A_* M \cong [D_AB, M]^A_*$
is isomorphic to
$$
\Hom_{B_*}(B^A_*(D_AB), M_*) \cong \Hom_{B_*}(B_*\langle G\rangle, M_*)
\cong \prod_G M_* \cong \pi_* F(G_+, M) \,.
$$
A diagram chase verifies that the K{\"u}nneth map equals the composite
of this chain of isomorphisms.
\qed
\enddemo

The more general situation, with $G$ an indiscrete stably dualizable
group, will lead to much more complicated spectral sequence calculations,
which we will not try to address here.

Goerss and Hopkins proceed to consider an {\it $E_2$-} or {\it resolution
model structure\/} on spectra, which is suitably generated by a class
$\Cal P$ of finite cellular spectra.  This class is required to satisfy
a list of conditions \cite{GH04, Def.~3.2.(1)--(5)}.  Following the
proof of \cite{BR:r, 2.2.4}, by Baker and Richter, we take $\Cal P$
to be the smallest set of dualizable $A$-modules that contains $A$
and $B$, and is closed under (de-)suspensions and finite wedge sums.
This immediately takes care of conditions~(3) and~(4).

\proclaim{Lemma 10.1.3}
The resolution model category conditions \cite{GH04, Def.~3.2.(1)--(5)}
are satisfied when $A \to B$ is a finite $E$-local $G$-Galois extension.
\endproclaim

\demo{Proof}
(1) $B^A_*(X)$ is a finite sum of shifted copies of $B^A_*(A) = B_*$
and $B^A_*(B) \cong \prod_G B_*$, for each $A$-module $X \in \Cal P$,
hence is projective as a $B_*$-module.  (2) $D_AB$ is represented in
$\Cal P$, since $B$ is self-dual as an $A$-module by Proposition~6.4.7.
(5) The K{\"u}nneth map
$$
[X, M]^A_* \to \Hom_{B_*}(B^A_*(X), M_*)_*
$$
is an isomorphism for all $B$-module spectra $M$ when $X = D_AB$, by the
Adams condition~(b), and trivially for $X = A$, so the same follows for
all $X \in \Cal P$ by passage to (de-)suspensions and finite wedge sums.
\qed
\enddemo

To sum up, a finite Galois extension $A \to B$ satisfies the Adams
conditions and has an associated resolution model structure on
$A$-modules, as required by \cite{GH04, \S3}, whenever $B$ is weakly
equivalent to a finite cell $A$-module.  It seems likely that the cited
theory also extends to cover all finite Galois extensions, by replacing
all references to finite cell objects by dualizable objects.  However, in
the following applications we shall always make use of the identification
$$
\Cal C_A(C, B) \cong \Cal C_B(B \wedge_A C, B)
$$
and only apply the Goerss--Hopkins spectral sequence in the case of
commutative $B$-algebra maps to $B$.  This is the very special case of
Lemmas~10.1.2 and~10.1.3 when $A = B$ and $G$ is the trivial group, in
which case $B$ is certainly a finite cell $A$-module.  So we are only
using the straightforward extension of \cite{GH04} to a more general
(cofibrant, commutative) ground $S$-algebra, namely~$B$.  Note also
that $B^B_*(B \wedge_A C) \cong B^A_*(C)$, so the two equivalent mapping
spaces above will have the same associated spectral sequences, which we
now review.

Goerss and Hopkins define Andr{\'e}--Quillen cohomology groups $D^s$
of algebras and modules over a simplicially resolved $E_\infty$-operad
\cite{GH04, (4.1)}, as non-abelian right derived functors of algebra
derivations.  They then construct a convergent spectral sequence
of Bousfield--Kan type \cite{GH04, Thm.~4.5}, which in our notation
appears as
$$
E_2^{s,t} \Longrightarrow \pi_{t-s} \Cal C_A(C, B)
\tag 10.1.4
$$
(based at a given commutative $A$-algebra map $C \to B$),
with $E_2$-term
$$
E_2^{0,0} = \Alg_{B_*}(B^A_*(C), B_*)
$$
and
$$
E_2^{s,t} = D^s_{B_*T}(B^A_*(C), \Omega^t B_*)
$$
for $t>0$.  Here $\Omega^t B_*$ is the $t$-th desuspension of the module
$B_*$.  As usual for Bousfield--Kan spectral sequences, this spectral
sequence is concentrated in the wedge-shaped region $0 \le s \le t$.

The subscript $B_*T$ refers to a (Reedy cofibrant, etc.) simplicial
$E_\infty$ operad~$T$ that resolves the commutative algebra operad
in the sense of \cite{GH04, Thm.~2.1}, and $B_*T$ is the associated
simplicial $E_\infty$ operad in the category of $B_*$-modules.
The Goerss--Hopkins Andr{\'e}--Quillen cohomology groups $D^s$ are the
right derived functors of derivations of $B_*T$-algebras in $B_*$-modules,
in the sense of Quillen's homotopical algebra.  As surveyed by Basterra
and Richter \cite{BR04, 2.6}, these groups $D^s$ do not depend on the
choice of resolving simplicial $E_\infty$ operad $T$, and agree with
the Andr{\'e}--Quillen cohomology groups $AQ^s_{sE_\infty}$ defined by
Mandell in \cite{Man03, 1.1} for $E_\infty$ simplicial $B_*$-algebras.
These do in turn agree with the Andr{\'e}--Quillen cohomology groups
$AQ^s_{dgE_\infty}$ defined by Mandell for $E_\infty$ differential
graded $B_*$-algebras \cite{Man03, 1.8}, and with Basterra's topological
Andr{\'e}--Quillen cohomology groups $\TAQ^s$ of the Eilenberg--Mac\,Lane
spectra associated to these algebras and modules \cite{Man03, \S7}.
By the comparison result of Basterra and McCarthy \cite{BMc02, 4.2}, these
are finally isomorphic to the $\Gamma$-cohomology groups $H\Gamma^s$
of Robinson and Sarah Whitehouse \cite{RoW02}, when $B^A_*(C)$ is
projective over $B_*$, or more generally, when $B^A_*(C)$ is flat over
$B_*$ and the universal coefficient spectral sequence from homology to
cohomology collapses.  So in these cases the Goerss--Hopkins groups can
be rewritten as
$$
D^s_{B_*T}(B^A_*(C), \Omega^t B_*) = H\Gamma^{s,-t}(B^A_*(C)|B_*; B_*) \,.
$$

It is not quite obvious from the above references that this chain of
identifications preserves the internal $t$-grading of these cohomology
groups, since this grading could be lost by the passage through
Eilenberg--Mac\,Lane spectra.  However, Birgit Richter has checked that
both gradings are indeed respected, up to the sign indicated above.
In our applications all of these cohomology groups will in fact be zero,
so the finer point about the internal grading is not so important.

If $B^A_*(C)$ is an {\'e}tale commutative $B_*$-algebra (thus flat
over $B_*$), then by \cite{RoW02, 6.8(3)} all $\Gamma$-homology and
$\Gamma$-cohomology groups of $B^A_*(C)$ over $B_*$ are zero, so by
the sequence of comparison results above (and the universal coefficient
spectral sequence for $\TAQ$), all the Goerss--Hopkins Andr{\'e}--Quillen
cohomology groups $D^s_{B_*T}(B^A_*(C), \Omega^t B_*)$ vanish.
Therefore one can conclude:

\proclaim{Corollary 10.1.5}
Let $C \to B$ be a map of commutative $A$-algebras.  If $B^A_*(C)$ is
{\'e}tale over $B_*$, then the Goerss--Hopkins spectral sequence for
$$
\pi_* \Cal C_A(C, B) \cong \pi_* \Cal C_B(B \wedge_A C, B)
$$
collapses to the origin at the $E_2$-term, so $\Cal C_A(C, B)$ is homotopy
discrete (each path component is weakly contractible) with
$$
\pi_0 \Cal C_A(C, B) \cong \Alg_{B_*}(B^A_*(C), B_*) \,.
$$
\endproclaim

\subhead 10.2. Idempotents and connected $S$-algebras \endsubhead

The converse part of the Galois correspondence, begun in Theorem~7.2.3,
should intrinsically characterize the intermediate extensions $A \to C
\to B$ that occur as $K$-fixed $S$-algebras $C = B^{hK}$ by allowable
subgroups $K \subset G$.  Already in the algebraic case of a $G$-Galois
extension $R \to T$ of discrete rings there are additional
complications (compared to the field case) when $T$ admits non-trivial
idempotents, i.e., when the spectrum of $T$ is not connected in the
sense of algebraic geometry.  See \cite{Mag74} for a general treatment
of these complications.  We do not expect that these issues are so
central to the extension of the theory from discrete rings to
$S$-algebras, so we prefer to focus on the analog of the situation when
$T$ is connected.

We can identify the idempotents of a commutative ring $T$ with the
non-unital $T$-algebra endomorphisms $T \to T$, taking an idempotent
$e$ (with $e^2=e$) to the homomorphism $t \mapsto et$.  The forgetful
functor from $T$-algebras to non-unital $T$-algebras has a left
adjoint, taking a non-unital $T$-algebra $N$ to $T \oplus N$, with the
multiplication $(t_1, n_1) \cdot (t_2, n_2) = (t_1t_2, t_1n_2 + n_1t_2
+ n_1n_2)$ and unit $(1,0)$.  In particular, we can identify the set of
idempotents $E(T) = \{e \in T \mid e^2=e \}$ with the set of
$T$-algebra maps
$$
E(T) \cong \Alg_T(T \oplus T, T) \,.
$$
Here $T \oplus T \cong T[x]/(x^2-x)$ is finitely generated and free as a
$T$-module.  It is {\'e}tale as a commutative $T$-algebra by \cite{Mil80,
I.3.4}, since $(x^2-x)' = 2x-1$ is its own multiplicative inverse in
$T[x]/(x^2-x)$.

This leads us to the following definitions.

\definition{Definition 10.2.1}
Let $B$ be a (cofibrant) commutative $S$-algebra.  Let the {\it space
of idempotents}
$$
\Cal E(B) = \Cal N_B(B, B)
$$
be the mapping space of non-unital commutative $B$-algebra \cite{Bas99,
\S1} endomorphisms $B \to B$.  We say that $B$ is {\it connected\/}
if the map $\{0,1\} \to \Cal E(B)$ taking $0$ and $1$ to the constant
map and the identity map $B \to B$, respectively, is a weak equivalence.
\enddefinition

We shall not have need to do so, but if we wanted to express that the
spectrum~$B$ has the property that $\pi_*(B) = 0$ for all $* \le 0$, we
would say that $B$ is $0$-connected, reserving the term ``connected'' for
the algebro-geometric interpretation just introduced.  A spectrum $B$ with
$\pi_*(B) = 0$ for $* < 0$ will be called $(-1)$-connected or connective.

There is a homeomorphism
$$
\Cal E(B) \cong \Cal C_B(B \vee B, B) \,,
$$
where $B \vee B$ is defined as the split commutative $S$-algebra extension
of $B$ with fiber the underlying non-unital commutative $S$-algebra
of $B$.  Its unit $B \to B \vee B$ is the inclusion on the first wedge
summand, and its multiplication is the composite
$$
(B \vee B) \wedge_B (B \vee B) \cong B \vee (B \vee B \vee B)
@>1\vee\nabla>> B \vee B
$$
where $\nabla$ folds the last three wedge summands together.

\proclaim{Proposition 10.2.2}
Let $B$ be any commutative $S$-algebra.  The space of idempotents $\Cal
E(B)$ is homotopy discrete, with $\pi_0 \Cal E(B) \cong E(\pi_0(B))$.
In particular, the commutative $S$-algebra $B$ is connected if and only
if the commutative ring $\pi_0(B)$ is connected.
\endproclaim

\demo{Proof}
We compute the homotopy groups of $\Cal E(B) \cong \Cal C_B(B \vee B,
B)$ by means of the Goerss--Hopkins spectral sequence~(10.1.4), in the
almost degenerate case when $A = B$ and $C = B \vee B$.  Here $A \to B$
is of course a $G$-Galois extension, in the trivial case $G = 1$, so our
discussion in Section~10.1 justifies the use of this spectral sequence.
It specializes to
$$
E_2^{s,t} \Longrightarrow \pi_{t-s} \Cal E(B)
$$
with
$$
E_2^{0,0} = \Alg_{B_*}(B_* \oplus B_*, B_*) = E(B_*)
$$
and
$$
E_2^{s,t} = D^s_{B_*T}(B_* \oplus B_*, \Omega^t B_*)
$$
for $t>0$.  Here $B_* \oplus B_* = B_*[x]/(x^2-x)$ is {\'e}tale
over $B_*$, so all the Andr{\'e}--Quillen cohomology groups $D^s =
H\Gamma^s$ vanish \cite{RoW02, 6.8(3)}, and we deduce that $\Cal E(B)$
is homotopy discrete, with $\pi_0 \Cal E(B) \cong E(B_*)$ equal to the
set of idempotents in the graded ring $B_*$, which of course are the
same as the idempotents in the ring $\pi_0(B)$.  In short, we have
applied Corollary~10.1.5.
\qed
\enddemo

The following argument, explained by Neil Strickland, shows that the above
definition of connectedness for structured ring spectra is equivalent to
another definition originally proposed by the author.  We say that an
$S$-algebra $B$ is {\it trivial\/} if it is weakly contractible, i.e.,
if $\pi_*(B) = B_* = 0$, and {\it non-trivial\/} otherwise.

\proclaim{Lemma 10.2.3}
A non-trivial commutative $S$-algebra $B$ is either connected, or weakly
equivalent to a product $B_1 \times B_2$ of non-trivial commutative
$B$-algebras, but not both.
\endproclaim

\demo{Proof}
If $B$ is non-trivial and not connected then there exists an idempotent
$e \in \pi_0(B)$ different from $0$ and $1$.  Let $f_1$ and $f_2 \: B \to
B$ be $B$-module maps inducing multiplication by $e$ and $1-e$ on
$\pi_*(B)$, respectively.  (These could also be taken to be non-unital
commutative $B$-algebra maps by the previous proposition.)  For $i=1,
2$ let $B[f_i^{-1}]$ be the mapping telescope for the iterated self-map
$f_i$, and let
$$
B_i = L^B_{B[f_i^{-1}]} B
$$
be the Bousfield $B[f_i^{-1}]$-localization of $B$ in the category of
$B$-modules.  Then there are commutative $B$-algebra maps $B \to B_1$
and $B \to B_2$ inducing isomorphisms $e \pi_*(B) \cong \pi_*(B_1)$
and $(1-e) \pi_*(B) \cong \pi_*(B_2)$, of nontrivial groups, and their
product $B \to B_1 \times B_2$ is the asserted weak equivalence.

Conversely, if $B \simeq B_1 \times B_2$ as commutative $B$-algebras
(or even just as ring spectra), with $B_1$ and $B_2$ non-trivial,
then $\pi_0(B)$ is not connected as a commutative ring, so $B$ is not
connected as a commutative $S$-algebra.
\qed
\enddemo

\subhead 10.3. Separable closure \endsubhead

The following terminology presumes, in some sense, that each finite
separable extension can be embedded in a finite Galois extension, i.e.,
a kind of normal closure.  We will not prove this in our context, but
keep the terminology, nonetheless.

\definition{Definition 10.3.1}
Let $A$ be a connected commutative $S$-algebra.  We say that $A$ is
{\it separably closed\/} if there are no $G$-Galois extensions $A \to B$
with $G$ finite and non-trivial and $B$ connected, i.e., if each finite
$G$-Galois extension $A \to B$ has $G = \{e\}$ or $B$ not connected.

A {\it separable closure\/} of $A$ is a pro-$G_A$-Galois extension
$A \to \bar A$ such that $\bar A$ is connected and separably closed.
The pro-finite Galois group $G_A$ of $\bar A$ over $A$ is the {\it
absolute Galois group\/} of $A$.
\enddefinition

The existence of a separable closure follows from Zorn's lemma.  However,
we have not yet proved that two separable closures of $A$ are weakly
equivalent, so talking of ``the'' absolute Galois group is also a bit
presumptive.

By Minkowski's theorem on the discriminant \cite{Ne99, III.2.17},
for every number field $K$ different from $\Q$ the inclusion $\Z \to
\Cal O_K$ is ramified at one or more primes.  In particular, there
are no Galois extensions $\Z \to \Cal O_K$ other than the identity.
The following inference appears to be well-known.

\proclaim{Proposition 10.3.2}
The only connected Galois extension of the integers is $\Z$ itself,
so $\Z = \bar\Z$ is separably closed.
\endproclaim

\demo{Proof}
Let $\Z \to T$ be a $G$-Galois extension of commutative rings, so $T$
is a finitely generated free $\Z$-module.  Then $\Q \to \Q \otimes T$
is also a $G$-Galois extension, so $\Q \otimes T \cong \prod_i K_i$ is
a product of number fields \cite{KO74, III.4.1}.  Then $T$ is contained
in the integral closure of $\Z$ in $\Q \otimes T$, which is a product
$\prod_i \Cal O_{K_i}$ of number rings.  The condition $T \otimes T \cong
\prod_G T$ and an index count imply, in combination, that $T = \prod_i
\Cal O_{K_i}$ and that each $\Cal O_{K_i}$ is unramified over $\Z$.
By Minkowski's theorem, this only happens when each $K_i =\Q$, so $T =
\prod_G \Z$.  If $T$ is connected, this implies that
$G$ is the trivial group and that $T = \Z$.
\qed
\enddemo

In other words, to have interesting Galois extensions of $\Z$ one must
localize away from one or more primes.  We have the following analog in
the context of commutative $S$-algebras.  The examples in Section~5.4
demonstrate that after localization there are indeed interesting
examples of (local) Galois extensions of $S$.

\proclaim{Theorem 10.3.3}
The only (global, finite) connected Galois extension of the sphere
spectrum $S$ is $S$ itself, so $S = \bar S$ is separably closed.
\endproclaim

\demo{Proof}
Let $S \to B$ be any finite $G$-Galois extension of global, i.e.,
unlocalized, commutative $S$-algebras (Definition~4.1.3).  Then $B$ is a
dualizable $S$-module (Proposition~6.2.1), hence of the homotopy type of
(a retract of) a finite CW spectrum (Proposition~3.3.3).  Thus $H_*(B)
= H_*(B; \Z)$ is finitely generated in each degree, and non-trivial only
in finitely many degrees.

Let $k$ be minimal such that $H_k(B) \ne 0$ and let $\ell$ be maximal
such that $H_\ell(B) \ne 0$.  The condition $B \wedge B \simeq \prod_G
B$ implies that $k = \ell = 0$.  For if $k < 0$ then $H_k(B) \otimes
H_k(B)$ is isomorphic to $H_{2k}(B \wedge B) \cong \prod_G H_{2k}(B) =
0$, which contradicts $H_k(B) \ne 0$ and finitely generated.  If $\ell >
0$ then $H_\ell(B) \otimes H_\ell(B)$ injects into $H_{2\ell}(B \wedge B)
\cong \prod_G H_{2\ell}(B) = 0$, which again contradicts $H_\ell(B) \ne 0$
and finitely generated.  Thus $H_*(B) = T$ is concentrated in degree~$0$.

By the Hurewicz theorem, $B$ is a connective spectrum with $\pi_0(B)
\cong H_0(B) = T$.  The K{\"u}nneth formula then implies that $T \otimes
T \cong \prod_G T$ and $\Tor_1^{\Z}(T, T) = 0$, so the unit map $\Z \to T$
makes $T$ a free abelian $\Z$-module of rank equal to the order of $G$.
In particular, $T$ is a faithfully flat $\Z$-module.

The result of inducing $B$ up along the Hurewicz map $S \to H\Z$ has
homotopy $\pi_*(H\Z \wedge B) = H_*(B) = T$ concentrated in degree~$0$,
so there is a pushout square
$$
\xymatrix{
B \ar[r] & HT \\
S \ar[u] \ar[r] & H\Z \ar[u]
}
$$
of commutative $S$-algebras.  By a variation on the proof of Lemma~7.1.1,
we shall now show that $H\Z \to HT$ is $G$-Galois.

The map $HT \wedge_{H\Z} HT \to \prod_G HT$ is induced up from the weak
equivalence $B \wedge B \to \prod_G B$, cf.~diagram~(7.1.2), and is
therefore a weak equivalence.  Next, $S \to B$ is dualizable, so $H\Z
\to HT$ is dualizable (Lemma~6.2.3).  Finally, $T$ is faithfully flat
over $\Z$ and so $HT$ is faithful over $H\Z$ by the proof of Lemma~4.3.5.
Thus $H\Z \to HT$ is a faithful $G$-Galois extension (Proposition~6.3.2).

From Proposition~4.2.1 we deduce that $\Z \to T$ is a $G$-Galois extension
of commutative rings.  By the classical theorem of Minkowski, this is
only possible if $G = \{e\}$ is the trivial group or $T$ is not connected.
And $\pi_0(B) \cong T$, so either $G$ is trivial or $B$ is not connected
(Proposition~10.2.2).  Thus $S$ is separably closed.
\qed
\enddemo

Note that we did not have to (possibly) restrict attention to faithful
$G$-Galois extensions $S \to B$ in this proof.

\example{Question 10.3.4}
Can the absolute Galois group $G_A$, or its maximal abelian quotient
$G_A^{ab}$, be expressed in terms of arithmetic invariants of $A$, such
as its algebraic $K$-theory $K(A)$?  This would constitute a form of
class field theory for commutative $S$-algebras.  The author expects that
there is a better hope for a simple answer in the maximally localized
category of $K(n)$-local commutative $S$-algebras, than for general
commutative $S$-algebras.
\endexample

\example{Question 10.3.5}
If an $E$-local commutative $S$-algebra $A$ is an even periodic Landweber
exact spectrum, and $A \to B$ is a finite $E$-local $G$-Galois extension,
does it then follow that $B$ is also an even periodic Landweber exact
spectrum, and that $\pi_0(A) \to \pi_0(B)$ is a $G$-Galois extension of
commutative rings?
\endexample

In the case of $E = K(n)$ and $A = E_n^{nr}$, for which $\pi_0(A) =
\W(\bar\F_p)[[u_1, \dots, u_{n-1}]]$ is separably closed, there are no
non-trivial such algebraic extensions to a connected ring $\pi_0(B)$,
so it would follow that $E_n^{nr}$ is $K(n)$-locally separably closed.
In particular, $E_n^{nr}$ would be the $K(n)$-local separable closure
of $L_{K(n)} S$, with absolute Galois group $\G_n^{nr} = \SS_n \rtimes
\hat\Z$.  This amounts to Conjecture~1.3 in the introduction.

Baker and Richter \cite{BR:r} have partial results in this direction,
in the global category.  They are able to show that $E_n^{nr}$ does not
admit any non-trivial connected faithful abelian $G$-Galois
extensions.  So $E_n^{nr} = E_n^{ab}$ is the maximal global faithful
abelian extension of $E_n$.

\head 11. Galois theory II \endhead

As before, we are implicitly working $E$-locally, for some spectrum $E$.

\subhead 11.1. Recovering the Galois group \endsubhead

The space of commutative $A$-algebra endomorphisms of $B$ in a $G$-Galois
extension $A \to B$ can be rewritten as
$$
\Cal C_A(B, B) \cong \Cal C_B(B \wedge_A B, B)
\simeq \Cal C_B(F(G_+, B), B) \,,
$$
in view of the weak equivalence $h \: B \wedge_A B \to F(G_+, B)$.
When $G$ is finite and discrete, and $B$ admits no non-trivial
idempotents, we can compute the homotopy groups of this mapping space
by the Goerss--Hopkins spectral sequence.

When $G$ is not discrete, these spectral sequence computations appear
to be much harder, and we will not attempt them.  We are therefore
principally working in the context of the separable/{\'e}tale extensions
from Chapter~9.

\proclaim{Theorem 11.1.1}
Let $A \to B$ be a finite $G$-Galois extension of commutative
$S$-algebras, with $B$ connected.  Then the natural map
$$
G \to \Cal C_A(B, B) \,,
$$
giving the action of $G$ on $B$ through commutative $A$-algebra maps,
is a weak equivalence.  In particular, $\Cal C_A(B, B)$ is a homotopy
discrete grouplike monoid, so each commutative $A$-algebra endomorphism
of $B$ is an automorphism, up to a contractible choice.
\endproclaim

\demo{Proof}
This time we compute the homotopy groups of $\Cal C_A(B, B) \simeq
\Cal C_B(\prod_G B, B)$ by means of~(10.1.4), once again in the
almost degenerate case when $A = B$ and $C = F(G_+, B) = \prod_G B$.
The $E_2$-term has
$$
E_2^{0,0} = \Alg_{B_*}(\prod_G B_*, B_*) \cong G
$$
since $B_* = \pi_*(B)$ is connected in the graded sense, or equivalently,
$\pi_0(B)$ has no non-trivial idempotents.  The remainder of the
$E_2$-term is
$$
E_2^{s,t} = D^s_{B_*T}(\prod_G B_*, \Omega^t B_*) = 0 \,,
$$
since $\prod_G B_*$ is {\'e}tale over $B_*$.  We are therefore in the
collapsing situation of Corollary~10.1.5, and $\Cal C_A(B, B) \simeq
G$ follows.
\qed
\enddemo

The extension to profinite pro-Galois extensions is straightforward.

\proclaim{Proposition 11.1.2}
Let $A \to B = \colim_\alpha B_\alpha$ be a pro-$G$-Galois extension,
with each $A \to B_\alpha$ a finite $G_\alpha$-Galois extension and $G =
\lim_\alpha G_\alpha$.  Suppose that $B$ is connected.  Then $\Cal C_A(B,
B)$ is homotopy discrete, and the natural map $G \to \pi_0 \Cal C_A(B,
B)$ is a group isomorphism.
\endproclaim

\demo{Proof}
Using~(8.1.2), we rewrite the commutative $A$-algebra mapping space as
$$
\multline
\Cal C_A(B, B) \cong \Cal C_B(B \wedge_A B, B) \\
\simeq \Cal C_B(\colim_\alpha F(G_{\alpha+}, B), B)
\simeq \holim_\alpha \Cal C_B(F(G_{\alpha+}, B), B) \,.
\endmultline
$$
By the finite case, each $\Cal C_B(F(G_{\alpha+}, B), B)$ is homotopy
discrete with
$$
\pi_0 \Cal C_B(F(G_{\alpha+}, B), B) \cong \Alg_{B_*}(\prod_{G_\alpha}
B_*, B_*) \cong G_\alpha \,,
$$
when $B$ is connected.  So $\Cal C_A(B, B) \simeq \holim_\alpha G_\alpha$
is homotopy discrete, with $\pi_0 \Cal C_A(B, B) \cong \lim_\alpha
G_\alpha \cong G$.
\qed
\enddemo

\subhead 11.2. The brave new Galois correspondence \endsubhead

We now turn to the converse part of the Galois correspondence.  The proper
role of the separability condition in the following result was found
in a conversation with Birgit Richter.

\proclaim{Proposition 11.2.1}
Let $A \to B$ be a $G$-Galois extension, with $B$ connected and $G$
finite and discrete, and let
$$
A \to C \to B
$$
be a factorization of this map through a separable commutative $A$-algebra
$C$.  Then $\Cal C_C(B, B)$ is homotopy discrete, and the natural map
$\Cal C_C(B, B) \to \Cal C_A(B, B)$ identifies $K = \pi_0 \Cal C_C(B,
B)$ with a subgroup of $G = \pi_0 \Cal C_A(B, B)$.  Furthermore, the
action of $\Cal C_C(B, B) \simeq K$ on $B$ induces a weak equivalence
$$
h \: B \wedge_C B \to \prod_K B \,.
$$
\endproclaim

\demo{Proof}
By assumption $A \to C$ is separable, so there are maps
$$
C' @>\sigma>> C \wedge_A C @>\mu>> C
$$
of $C$-bimodules relative to $A$ such that $\mu \sigma \: C' \to C$ is
a weak equivalence.  Inducing these maps and modules up along $C \to B$,
both as left and right modules, we get maps
$$
B \wedge_C C' \wedge_C B @>\bar\sigma>> B \wedge_A B
	@>\bar\mu>> B \wedge_C B
$$
of $B$-bimodules relative to $A$, such that the composite is a weak
equivalence.  We consider $C \wedge_A C$ as a commutative $C$-algebra via
the left unit $C \cong C \wedge_A A \to C \wedge_A C$, and similarly for
$B \wedge_A B$ over $B$.  Then $\mu$ is a map of commutative $C$-algebras
and $\bar\mu$ is a map of commutative $B$-algebras.

At the level of homotopy groups, we get a diagram
$$
B^C_*(B) @>\bar\sigma_*>> B^A_*(B) @>\bar\mu_*>> B^C_*(B)
$$
of $B^A_*(B)$-module homomorphisms, whose composite is the identity.
Furthermore, $\bar\mu_*$ is a $B_*$-algebra homomorphism.  It follows
from the $B^A_*(B)$-linearity of the homomorphism $\bar\sigma_*$ that
it is also a $B_*$-algebra homomorphism.  For if $x, y \in B^C_*(B)$
then $\bar\sigma_* x \in B^A_*(B)$ acts on $y$ through multiplication
by its image $\bar\mu_* \bar\sigma_*x = x$ and on $\bar\sigma_* y \in
B^A_*(B)$ by multiplication by $\bar\sigma_* x$.  The
$B^A_*(B)$-linearity of $\bar\sigma_*$ now provides the left hand
equality below:
$$
\bar\sigma_* x \cdot \bar\sigma_* y
	= \bar\sigma_*(\bar\mu_* \bar\sigma_* x \cdot y)
	= \bar\sigma_* (x \cdot y) \,.
$$
Thus $\bar\sigma_*$ is a $B_*$-algebra homomorphism, so that $B^C_*(B)$
is a retract of $B^A_*(B)$, both in the category of $B^A_*(B)$-modules
and, more importantly to us, in the category of commutative
$B_*$-algebras.

Recall that $B^A_*(B) \cong \prod_G B_*$, since $A \to B$ is $G$-Galois.
Here $\prod_G B_* \cong \bigoplus_G B_*$ is a finitely generated free
$B_*$-module, since $G$ is finite, so the retraction above implies
that $B^C_*(B)$ is a finitely generated projective $B_*$-module.
We may therefore once more consider the Goerss--Hopkins spectral
sequence~(10.1.4), now for the mapping space
$$
\Cal C_C(B, B) \cong \Cal C_B(B \wedge_C B, B) \,.
$$

The $E_2$-term has
$$
E_2^{s,t} = D^s_{B_*T}(B^C_*(B), \Omega^t B_*)
$$
for $t>0$.  The commutative $B_*$-algebra retraction $\bar\mu_* \:
B^A_*(B) \to B^C_*(B)$ induces a split injection from each of these
cohomology groups to
$$
D^s_{B_*T}(B^A_*(B), \Omega^t B_*) \,,
$$
which we saw was zero in the proof of Theorem~11.1.1, since $B^A_*(B)
= \prod_G B_*$ is {\'e}tale over $B_*$.  We therefore have $E_2^{s,t} =
0$ away from the origin, also in the Goerss--Hopkins spectral sequence
for $\pi_{t-s} \Cal C_C(B, B)$.

Thus $\Cal C_C(B, B)$ is homotopy discrete, in the sense that each path
component is weakly contractible, with set of path components
$$
K = \pi_0 \Cal C_C(B, B) \cong \Alg_{B_*}(B^C_*(B), B_*) \,.
$$
The $B_*$-algebra retraction $\bar\mu_*$ induces a split injection from
this set to
$$
\Alg_{B_*}(B^A_*(B), B_*) \cong G \,.
$$
It is clear that the natural map $\Cal C_C(B, B) \to \Cal C_A(B,
B)$, viewing a map $B \to B$ of commutative $C$-algebras as a map of
commutative $A$-algebras, is a monoid map with respect to the composition
of maps.  Therefore the injection $K \to G$ identifies $K$ as a sub-monoid
of $G$.  But $G$ is a finite group, so $K$ is in fact a subgroup of $G$.
This completes the proof of the first claims of the proposition.

The tautological action by $\Cal C_C(B, B)$ on $B$ through commutative
$C$-algebra maps can be converted to an action by $K$ on a commutative
$C$-algebra $B'$ weakly equivalent to $B$.  We hereafter implicitly make
this replacement, so as to have $K$ acting directly on $B$ over $C$,
and turn to the proof of the final claim.

The composite $\bar\sigma_* \bar\mu_* \: B^A_*(B) \to B^A_*(B)$ is an
idempotent $B_*$-algebra map.  Under the isomorphism $B^A_*(B) \cong
\prod_G B_*$ it corresponds to an idempotent $B_*$-algebra endomorphism
of $\prod_G B_*$.  Since $B_*$ is connected, it must be the retraction
of $\prod_G B_*$ onto the subalgebra $\prod_{K'} B_*$, for some subset
$K' \subset G$ (containing $e \in G$).  Thus $B^C_*(B) \cong \prod_{K'}
B_*$, which implies
$$
K = \Alg_{B_*}(B^C_*(B), B_*) \cong \Alg_{B_*}(\prod_{K'} B_*, B_*)
\cong K' \,.
$$
Thus $K = K'$ as subsets of $G$, and the weak equivalence $B \wedge_A B
\simeq \prod_G B$ retracts to a weak equivalence
$$
h \: B \wedge_C B \to \prod_K B
\,.
$$
It is quite clearly given by the action of $K$ on $B$ through commutative
$C$-algebra maps, as in~(4.1.2).
\qed
\enddemo

This leads us to the converse part of the Galois correspondence for
$E$-local commutative $S$-algebras, in the case of finite, faithful
Galois extensions.

\proclaim{Theorem 11.2.2}
Let $A \to B$ be a $G$-Galois extension, with $B$ connected and $G$
finite and discrete.  Furthermore, let
$$
A \to C \to B
$$
be a factorization of this map through a separable commutative $A$-algebra
$C$ such that $C \to B$ is faithful, and let $K = \pi_0 \Cal C_C(B, B)
\subset G$.

If $A \to B$ is faithful, or more generally, if $B$ is dualizable over
$C$, then $C \simeq B^{hK}$ as commutative $C$-algebras, and $C \to B$
is a faithful $K$-Galois extension.
\endproclaim

\demo{Proof}
We first prove that $A \to B$ faithful implies that $B$ is dualizable
over $C$.

By hypothesis, $A \to C$ is separable, so the multiplication map $\mu$
and its weak section $\sigma$ make $C$ a retract up to homotopy of $C
\wedge_A C$, as a $C \wedge_A C$-module.  Therefore $\mu$ makes $C$
a dualizable $C \wedge_A C$-module, by Lemma~3.3.2(c).  Similarly,
for each $g \in G$ the twisted multiplication map
$$
\mu (1 \wedge g) \: C \wedge_A C \to C
$$
and its weak section $(1 \wedge g^{-1})\sigma$ make $C$ a dualizable $C
\wedge_A C$-module.  Inducing up along $C \to B$, Lemma~6.2.3 implies
that each map $B \wedge_A C \to B$, given algebraically as $b \wedge
c \mapsto b \cdot g(c)$, makes $B$ a dualizable $B \wedge_A C$-module.
By Lemma~3.3.2(c) again, it follows that the natural map $B \wedge_A C
\to B \wedge_A B$ makes $B \wedge_A B \simeq \prod_G B \simeq \bigvee_G B$
a dualizable $B \wedge_A C$-module.  By Proposition~6.2.1 and hypothesis, $B$
is dualizable and faithful over $A$, so by Lemma~6.2.3 and Lemma~4.3.3
we know that $B \wedge_A C$ is faithful and dualizable over $C$.
Thus by Lemma~6.2.4 it follows that the natural map $C \to B$ makes $B$
a dualizable $C$-module.

By Proposition~11.2.1, $K = \pi_0 \Cal C_C(B, B) \subset G$ acts on (a
weakly equivalent replacement for) $B$ through $C$-algebra maps, so that
$h \: B \wedge_C B \to \prod_K B$ is a weak equivalence.  By hypothesis
(and the argument above), $C \to B$ is faithful and dualizable.
Then by Proposition~6.3.2 the natural map $i \: C \to B^{hK}$ is a weak
equivalence, and so $C \to B$ is a faithful $K$-Galois extension.
\qed
\enddemo

\head 12. Hopf--Galois extensions in topology \endhead

In this final chapter we work globally, i.e., not implicitly localized
at any spectrum (other than at $E = S$).

\subhead 12.1.  Hopf--Galois extensions of commutative $S$-algebras
\endsubhead

Let $A \to B$ be a $G$-Galois extension of commutative $S$-algebras,
with $G$ stably dualizable, as usual.  The right adjoint $\tilde\alpha \:
B \to F(G_+, B)$ of the group action map $\alpha \: G_+ \wedge B \to B$
can be lifted up to homotopy through the weak equivalence $\nu\gamma \:
B \wedge DG_+ \to F(G_+, B)$, to a map $\beta \: B \to B \wedge DG_+$.
The group multiplication $G \times G \to G$ induces a functionally
dual map $DG_+ \to D(G \times G)_+$, which likewise can be lifted up to
homotopy through the weak equivalence $\wedge \: DG_+ \wedge DG_+ \to
D(G \times G)_+$ to a coproduct $\psi \: DG_+ \to DG_+ \wedge DG_+$.
We shall require rigid forms of these structure maps.

\definition{Definition 12.1.1}
A {\it commutative Hopf $S$-algebra\/} is a cofibrant commutative
$S$-algebra $H$ equipped with a counit $\epsilon \: H \to S$ and a
coassociative and counital coproduct $\psi \: H \to H \wedge_S H$,
in the category of commutative $S$-algebras.
\enddefinition

Note that we are not assuming that the coproduct $\psi$ is (strictly)
cocommutative, nor that it admits a strict antipode/conjugation $\chi \:
H \to H$.  This would severely limit the number of interesting examples.

\example{Example 12.1.2}
Let $X$ be an infinite loop space.  The $E_\infty$ structure on $X$
makes $S[X] = S \wedge X_+$ an $E_\infty$ ring spectrum.  The diagonal
map $\Delta \: X \to X \times X$ and $X \to *$ induce a coproduct $\psi
\: S[X] \to S[X \times X] \cong S[X] \wedge S[X]$ and counit $\epsilon
\: S[X] \to S$, which altogether can be rigidified to make $H \simeq
S[X]$ a commutative Hopf $S$-algebra.  The rigidification takes the
coassociative and counital coproduct and counit on $S[X]$ to a
corresponding co-$A_\infty$ structure on $H$, which in turn can be
rigidified to strictly coassociative and counital operations, by
working entirely within commutative $S$-algebras.  It is in general
not possible to make a similar rigidification of co-$E_\infty$ structures,
within commutative $S$-algebras.
\endexample

\definition{Definition 12.1.3}
Let $A$ be a cofibrant commutative $S$-algebra, let $B$ be a cofibrant
commutative $A$-algebra and let $H$ be a commutative Hopf $S$-algebra.
We say that $H$ {\it coacts\/} on $B$ over $A$ if there is a coassociative
and counital map
$$
\beta \: B \to B \wedge H
$$
of commutative $A$-algebras.  In this situation, let
$$
h \: B \wedge_A B \to B \wedge H
$$
be the composite map $(\mu \wedge 1)(1 \wedge \beta)$ of commutative
$B$-algebras.
\enddefinition

\definition{Definition 12.1.4}
The (Hopf) {\it cobar complex} $C^\bullet(H; B)$, for $H$ coacting on $B$
over~$A$, is the cosimplicial commutative $A$-algebra with
$$
C^q(H; B) = B \wedge H \wedge \dots \wedge H
$$
($q$ copies of $H$) in codegree~$q$.  The coface maps are $d_0 =
\beta \wedge 1^{\wedge q}$, $d_i = 1^{\wedge i} \wedge \psi \wedge
1^{\wedge(q-i)}$ for $0 < i < q$ and $d_q = 1^{\wedge q} \wedge \eta$,
where $\eta \: S \to H$ is the unit map.  The codegeneracy maps involve
the counit $\epsilon \: H \to S$.  Let $C(H; B) = \Tot C^\bullet(H; B)$
be its totalization.  The algebra unit $A \to B$ induces a coaugmentation
$A \to C^\bullet(H; B)$, and a map
$$
i \: A \to C(H; B) \,.
$$
\enddefinition

\definition{Definition 12.1.5}
A map $A \to B$ of commutative $S$-algebras is an {\it $H$-Hopf--Galois
extension\/} if $H$ is a commutative Hopf $S$-algebra that coacts on $B$
over $A$, so that the maps $i \: A \to C(H; B)$ and $h \: B \wedge_A B
\to B \wedge H$ are both weak equivalences.
\enddefinition

Note that there is no finiteness/dualizability condition on $H$ in
this definition.  See \cite{Chi00} for a recent text on Hopf--Galois
extensions in the algebraic setting.

\example{Example 12.1.6}
Let $G$ be a stably dualizable topological group.  The weak coproduct
on $DG_+ = F(G_+, S)$, derived from the group multiplication, can be
rigidified to give $H \simeq DG_+$ the structure of a commutative Hopf
$S$-algebra.  If $G$ acts on~$B$ over $A$, then the weak coaction of
$DG_+$ on $B$ can be rigidified to a coaction of~$H$ on $B$ over $A$.
Then the (Hopf) cobar complex $C^\bullet(H; B)$ maps by a degreewise
weak equivalence to the (group) cobar complex $C^\bullet(G; B)$ from
Definition~8.2.5.  In codegree~$q$ it is weakly equivalent to the
composite natural map
$$
B \wedge DG_+ \wedge \dots \wedge DG_+ @>\wedge>\simeq>
	B \wedge DG^q_+ @>\gamma>\cong>
	DG^q_+ \wedge B @>\nu>\simeq> F(G^q_+, B) \,.
$$
On totalizations, we obtain a weak equivalence $C(H; B) \simeq B^{hG}$.
In this case, the definition of an $H$-Hopf--Galois extension $A \to B$
generalizes that of a $G$-Galois extension $A \to B$, since $i \: A \to
B^{hG}$ factors as
$$
A @>i>> C(H; B) @>\simeq>> B^{hG} \,,
$$
and $h \: B \wedge_A B \to \prod_G B$ factors as
$$
B \wedge_A B @>h>> B \wedge H @>\simeq>> F(G_+, B) \,.
$$
\endexample

Recall the Amitsur complex $C^\bullet(B/A)$ from Definition~8.2.1.

\definition{Definition 12.1.7}
There is a natural map of cosimplicial commutative $A$-algebras
$h^\bullet \: C^\bullet(B/A) \to C^\bullet(H; B)$ given in
codegree~$q$ by the map
$$
h^q \: B \wedge_A B \wedge_A \dots \wedge_A B
	\to B \wedge H \wedge \dots \wedge H
$$
that is the composite of the maps
$$
\multline
B^{\wedge_A (i+1)} \wedge H^{\wedge j} \cong
B^{\wedge_A (i-1)} \wedge_A (B \wedge_A B) \wedge H^{\wedge j} \\
@>1^{\wedge (i-1)} \wedge h \wedge 1^{\wedge j}>>
B^{\wedge_A (i-1)} \wedge_A (B \wedge H) \wedge H^{\wedge j}
\cong B^{\wedge_A i} \wedge H^{\wedge (j+1)}
\endmultline
$$
for $j = 0, \dots, q-1$ and $i+j=q$.
Upon totalization, it induces a map $h' \: A^\wedge_B \to C(H; B)$
of commutative $A$-algebras.
\enddefinition

The diagram chase needed to verify that $h^\bullet$ indeed is cosimplicial
uses the strict coassociativity and counitality of the Hopf $S$-algebra
structure on $H$.

\proclaim{Proposition 12.1.8}
Suppose that $H$ coacts on $B$ over $A$, as above, and that $h \: B
\wedge_A B \to B \wedge H$ is a weak equivalence.  Then $h' \: A^\wedge_B
\to C(H; B)$ is a weak equivalence.  As a consequence, $A \to B$ is an
$H$-Hopf--Galois extension if and only if $A$ is complete along $B$.
\endproclaim

\demo{Proof}
The cosimplicial map $h^\bullet$ is a weak equivalence in each codegree,
so the induced map of totalizations $h'$ is a weak equivalence.
Therefore the composite $i = h' \circ \eta$, of the two maps
$$
A @>\eta>> A^\wedge_B @>h'>> C(H; B) \,,
$$
is a weak equivalence if and only if $\eta$ is one.
\qed
\enddemo

\subhead 12.2. Complex cobordism \endsubhead

Let $A = S$ be the sphere spectrum, $B = MU$ the complex cobordism
spectrum and $H = S[BU] = \Sigma^\infty BU_+$ the unreduced suspension
spectrum of $BU$.  Bott's infinite loop space structure on $BU$ makes
$H$ a commutative $S$-algebra, and the diagonal map $\Delta \: BU
\to BU \times BU$ induces the Hopf coproduct $\psi \: S[BU] \to S[BU]
\wedge S[BU]$.  The Thom diagonal
$$
\beta \: MU \to MU \wedge BU_+
$$
defines a coaction by $S[BU]$ on $MU$ over $S$.  The induced map
$$
h \: MU \wedge MU \to MU \wedge BU_+
$$
is the weak equivalence known as the Thom isomorphism.  The Bousfield--Kan
spectral sequence associated to the cosimplicial commutative $S$-algebra
$C^\bullet(MU/S)$ is the Adams--Novikov spectral sequence
$$
E_2^{s,t} = \Ext^{s,t}_{MU_*MU}(MU_*, MU_*)
\Longrightarrow \pi_{t-s}(S) \,.
$$
The convergence of this spectral sequence is the assertion that the
coaugmentation
$$
i \: S \to S^\wedge_{MU} = \Tot C^\bullet(MU/S)
$$
is a weak equivalence.  In view of Proposition~12.1.8 we can summarize
these facts as follows:

\proclaim{Proposition 12.2.1}
The unit map $S \to MU$ is an $S[BU]$-Hopf--Galois extension of
commutative $S$-algebras. \qed
\endproclaim

\remark{Remark 12.2.2}
There is no topological group $G$ such that $S \to MU$ is a $G$-Galois
extension, but $S[BU]$ is taking on the role of its functional dual
$DG_+$, as in Example~12.1.6.  So the commutative Hopf $S$-algebra
$S[BU]$ is trying to be the ring of functions on the non-existent
Galois group of $MU$ over $S$.  Note that there is no bimodule section
to the multiplication map $\mu \: MU \wedge MU \to MU$, since the left
and right units $\eta_L, \eta_R \: MU_* \to MU_* MU$ are really
different, so $S \to MU$ is not separable in the sense of Section~9.1.
\endremark

\remark{Remark 12.2.3}
There are similar $S[X]$-Hopf--Galois extensions $S \to Th(\gamma)$
to the Thom spectrum induced by any infinite loop map $\gamma \: X \to
BGL_1(S)$.  For example, there is such an extension $S \to MUP$ to the
even periodic version $MUP$ of $MU$, which is the Thom spectrum of the
tautological virtual bundle over $X = \Z \times BU = \Omega^\infty ku$.
More generally, for any commutative $S$-algebra $R$ and infinite loop
map $\gamma \: X \to BGL_1(R)$ there is an $R$-based Thom spectrum
$Th^R(\gamma)$, i.e., a $\gamma$-twisted form of $R[X] = R \wedge X_+$,
and an $R[X]$-Hopf--Galois extension $R \to Th^R(\gamma)$.
\endremark

\remark{Remark 12.2.4}
The extension $S \to MU$ is known not to be faithful, since by \cite{Ra84,
\S3} or \cite{Ra92, 7.4.2} $MU_*(cY) = 0$ for every finite complex $Y$
with trivial rational cohomology.  Here $cY$ denotes the Brown--Comenetz
dual of $Y$.  This faithlessness leaves the telescope conjecture
\cite{Ra84, 10.5} or \cite{Ra92, 7.5.5} a significant chance to be false.
Recall that if $F(n)$ is a finite complex of type~$n$ (with a $v_n$-self
map), and $T(n) = v_n^{-1} F(n)$ its mapping telescope, the conjecture is
that the natural map $\lambda \: T(n) \to L_n F(n)$ is a weak equivalence.
After inducing up to $MU$, $1 \wedge \lambda \: MU \wedge T(n) \to MU
\wedge L_n F(n)$ is an equivalence, by the localization theorem $v_n^{-1}
MU \wedge F(n) \simeq L_n MU \wedge F(n)$ \cite{Ra92, 7.5.2}.  Positive
information about the faithfulness of Galois- or Hopf--Galois extensions
(Question~4.3.6) might conceivably reflect back on this conjecture.
\endremark
\medskip

To conclude this paper, we wish to discuss how the Hopf--Galois
extension $S \to MU$ provides a global, integral object whose
$p$-primary $K(n)$-localization and nilpotent completion $L_{K(n)} S
\to \hat L^{MU}_{K(n)} MU$ governs the pro-Galois extensions $L_{K(n)}
S \to E_n$, for each rational prime~$p$ and integer $n\ge0$.
$$
\xymatrix@!C{
MU & \hat L^{MU}_{K(n)} MU \ar[r]^-t & E_n \\
S \ar[u]^{S[BU]} & L_{K(n)} S \ar[u] \ar[ur]_{\G_n}
}
\tag 12.2.5
$$
This suggests that $S \to MU$ is a kind of near-maximal ramified Galois
extension, and that its weak Galois group (``weak'' in the analytic
sense that it is only realized through its functional dual $DG_+ =
S[BU]$ that coacts on $MU$) is a kind of near-absolute ramified Galois
group of the sphere.  More precisely, the maximal extension may be the
one obtained from the even periodic theory $MUP$ by tensoring with the
ring $\Cal O_{\bar\Q}$ of algebraic integers.

Even if $S \to MU$ does not admit many Galois automorphisms, the Hopf
coaction $\beta \: MU \to MU \wedge BU_+$ still determines the Galois
action of each element $g \in \G_n$ on~$E_n$.  By the Hopkins--Miller
theory, each commutative $S$-algebra map $g \: E_n \to E_n$ is uniquely
determined by the underlying map of (commutative) ring spectra, so it
is the description of the latter that we shall review.

Recall from~5.4.2 that $\Gamma_n$ is the Honda formal group law over
$\F_{p^n}$ and $\widetilde\Gamma_n$ its universal deformation, defined
over $\pi_0(E_n)$.  By the Lubin--Tate theorem \cite{LT66, 3.1}, each
automorphism $g \in \SS_n \subset \G_n$ of $\Gamma_n$ determines a
unique pair $(\phi, \tilde g)$, where $\phi \: \pi_0(E_n) \to
\pi_0(E_n)$ is a ring automorphism and $\tilde g \: \widetilde\Gamma_n
\to \phi_* \widetilde\Gamma_n$ is an isomorphism of formal group laws
over $\pi_0(E_n)$, whose expansion $\tilde g(x) \in \pi_0(E_n)[[x]]
\cong E_n^0(\C P^\infty)$ reduces modulo $(p, u_1, \dots, u_{n-1})$ to
the expansion $g(x) \in \F_{p^n}[[x]]$ of $g$.  Then $\phi = \pi_0(g)$,
when $g$ is considered as a self-map of $E_n$.  Furthermore,
$$
\align
\tilde g(x) \in E_n^0(\C P^\infty)
	& \cong \Hom_{E_{n*}}(E_{n*}(\C P^\infty), E_{n*}) \\
	& \cong \Alg_{E_{n*}}(E_{n*}(BU), E_{n*})
	\subset E_n^0(BU)
\endalign
$$
corresponds to a unique map of ring spectra $\tilde g \: S[BU] \to E_n$ .
Let $t \: MU \to E_n$ be the usual complex orientation, corresponding to
the graded version of $\widetilde\Gamma_n$.  Then the following diagram
commutes up to homotopy:
$$
\xymatrix{
MU \ar[r]^-{\beta} \ar[d]_t
	& MU \wedge BU_+ \ar[d]^{\mu(t \wedge \tilde g)} \\
E_n \ar[r]^g & E_n
}
\tag 12.2.6
$$
The composite $g \circ t = \mu (t \wedge \tilde g) \beta$ determines $g$,
in view of $t^* \: E_n^*(E_n) \to E_n^*(MU)$ being nearly injective.
Only the Galois automorphisms in $\Gal \subset \G_n$ are missing, but
these may be ignored if we are focusing on $\widehat{E(n)}$, or can be
detected by passing to $MUP$ and adjoining some roots of unity.

By analogy, for number fields $K \subset L$ and primes $\goth p \in
\Cal O_K$, a factorization $\goth p \Cal O_L = \goth P_1^{e_1} \cdot
\dots \cdot \goth P_r^{e_r}$ leads to a splitting of completions
$K_{\goth p} \to L \otimes_K K_{\goth p} \cong \prod_i L_{\goth P_i}$.
If the field extension $K \to L$ is $G$-Galois, then each local
extension $K_{\goth p} \to L_{\goth P_i}$ is $G_{\goth P_i}$-Galois,
where $G_{\goth P_i} \subset G$ is the decomposition group of $\goth
P_i$, and $G$ acts transitively on the finite set of primes over $\goth
p$.  Thus when the global extension $K \to L$ is localized (i.e.,
completed), it splits as a product of smaller local extensions, in a
way that depends on the place of localization.
$$
\xymatrix{
L & L \otimes_K K_{\goth p} \ar[r]^-{pr_i} & L_{\goth P_i} \\
K \ar[u]^G & K_{\goth p} \ar[u] \ar[ur]_{G_{\goth P_i}}
}
\tag 12.2.7
$$
In the algebraic case of a pro-Galois extension $K \to \bar K$ there
is a profinite set of places over each prime~$\goth p$, still forming
a single orbit for the action by the absolute Galois group $G_K$.

In the topological setting of $S \to MU$ there is likewise a single
orbit of chromatic primes of $MU$ over the one of $S$ that corresponds
to the localization functor $L_{K(n)}$ on $\Cal M_S$, namely those
corresponding to the nilpotent completion functors $\hat L^{MU}_{K(n)}$
on $\Cal M_{MU}$, for all the various possible complex orientations $MU
\to K(n)$.  In the absence of a real Galois group of automorphisms of
$MU$ these do not form a geometric orbit of places, but the next-best
thing is available, namely the $S$-algebraic coaction by $S[BU]$ via
the Thom diagonal, a sub-coaction of which indeed links the various
complex orientations of $K(n)$ into one ``weak'' orbit.

Jack Morava \cite{Mo05} has developed this Galois theoretic perspective
on the stable homotopy category further.

\enddocument